\documentclass[11pt,a4paper]{article}
\pdfoutput=1

\usepackage{a4wide}
\usepackage{tikz}
\usetikzlibrary{arrows}
\usetikzlibrary{positioning}
\usetikzlibrary{arrows.meta}

\usepackage{enumitem}
\usepackage{amsmath}
\usepackage{subcaption}
\usepackage{xr}
\RequirePackage{hyperref}
\usepackage{epsfig}
\usepackage{amssymb}
\usepackage{natbib}
\usepackage{graphicx}
\usepackage{color}
\usepackage{setspace}
\usepackage{xspace}
\usepackage{tabularx}
\usepackage[T1]{fontenc}
\usepackage{blkarray}
\newcommand\bigzero{\makebox(0,0){\text{\huge0}}}

\newtheorem{definition}{Definition}
\newtheorem{theorem}[definition]{Theorem} 
\newtheorem{lemma}[definition]{Lemma}

\newtheorem{corollary}[definition]{Corollary}

\newtheorem{example}{Example}

\newcommand{\ddo}{\textrm{do}}
\newcommand{\bX}{\mathbf{X}}
\newcommand{\bY}{\mathbf{Y}}

\newcommand{\bZ}{\mathbf{Z}}

\newcommand\given[1][]{\:#1\vert\:}

\newcommand{\dsepp}{\perp_{d}}

\newtheorem{claim}{Claim}

      \newenvironment{proofof}[1][]{\begin{trivlist}
   \item[\hskip \labelsep {\bfseries Proof of #1.}]}{\hfill{}$\square$\end{trivlist}}

\definecolor{byzantium}{rgb}{0.44, 0.16, 0.39}
\definecolor{burgundy}{rgb}{0.5, 0.0, 0.13}
\definecolor{blue-violet}{rgb}{0.54, 0.17, 0.89}
\definecolor{antiquefuchsia}{rgb}{0.57, 0.36, 0.51}
\definecolor{amethyst}{rgb}{0.6, 0.4, 0.8}
\definecolor{blue-violet}{rgb}{0.54, 0.17, 0.89}
\definecolor{ao}{rgb}{0.0, 0.5, 0.0}
\definecolor{blue(ncs)}{rgb}{0.0, 0.53, 0.74}

\newcommand{\tild}{\raise.17ex\hbox{ $\scriptstyle\sim$ }}

\newcommand{\specialcell}[2][l]{%
  \begin{tabular}[#1]{@{}l@{}}#2\end{tabular}}

\DeclareMathOperator{\Det}{Det}
\DeclareMathOperator{\Cov}{Cov}
\DeclareMathOperator{\Var}{Var}
\DeclareMathOperator{\CPDAG}{CPDAG}
\DeclareMathOperator{\DAG}{DAG}
\DeclareMathOperator{\MAG}{MAG}
\DeclareMathOperator{\PAG}{PAG}
\DeclareMathOperator{\PossDe}{PossDe}

\DeclareMathOperator{\De}{De}
\DeclareMathOperator{\An}{An}
\DeclareMathOperator{\PossAn}{PossAn}
\DeclareMathOperator{\Pa}{Pa}
\DeclareMathOperator{\Adj}{Adj}
\DeclareMathOperator{\Adjust}{Adjust}
\DeclareMathOperator{\distancefrom}{distance-from-\!}

\DeclareMathOperator{\Forbi}{Forb}

\newcommand{\amen}{the amenability condition\xspace}
\newcommand{\forb}{the forbidden set condition\xspace}
\newcommand{\blck}{the blocking condition\xspace}
\newcommand{\condtwoprime}{the separation condition\xspace}

\newcommand{\vars}[1][V]{\mathbf{#1}}
\newcommand{\e}[1][E]{\mathbf{#1}}

\newcommand{\g}[1][G]{\mathcal{#1}}
\newcommand{\gout}[2][G]{\mathcal{#1}_{\underline{#2}}}

\newcommand{\gpbd}[2][G]{\mathcal{#1}_{\mathbf{#2}}^{pbd}}

\newcommand{\dsep}[2][X,Y]{DSEP(#1,#2)}

\newcommand{\adjust}[2][X,Y]{\Adjust(#1,#2)}
\newcommand{\adjustb}[2][X,Y]{\Adjust(\mathbf{#1},#2)}

\newcommand{\f}[2][X,Y]{\Forbi(#1,#2)}
\newcommand{\fb}[2][X,Y]{\Forbi(\mathbf{#1},#2)}

\newcommand{\pstar}[1][p]{{#1}^{*}}

\newcommand{\bulletbullet}{
  \setlength{\unitlength}{1mm}
  \begin{picture}(5,1)(0,0)
    \put(0.2,0){$\bullet$}
    \put(1.3,1){\line(1,0){2.4}}
    \put(2.8,0){$\bullet$}
  \end{picture}
}

\newcommand{\bulletcirc}{
  \setlength{\unitlength}{1mm}
  \begin{picture}(5,1)(0,0)
    \put(0.2,0){$\bullet$}
    \put(1.1, 1){\line(1,0){2.4}}
    \put(4, 1){\circle{1}}
  \end{picture}
}

\newcommand{\circbullet}{
  \setlength{\unitlength}{1mm}
  \begin{picture}(5,1)(0,0)
    \put(1,1){\circle{1}}
    \put(1.5,1){\line(1,0){2.4}}
    \put(2.9,0){$\bullet$}
  \end{picture}
}

\newcommand{\bulletarrow}{
  \setlength{\unitlength}{1mm}
  \begin{picture}(5,1)(0,0)
    \put(0.2,0){$\bullet$}
    \put(1,0){$\rightarrow$}
  \end{picture}
}

\newcommand{\arrowbullet}{
  \setlength{\unitlength}{1mm}
  \begin{picture}(5,1)(0,0)
    \put(0.2,0){$\leftarrow$}
    \put(3,0){$\bullet$}
  \end{picture}
}
\newcommand{\tailbullet}{
  \setlength{\unitlength}{1mm}
  \begin{picture}(5,1)(0,0)
    \put(0.4,1){\line(1,0){3.2}}
    \put(2.8,0){$\bullet$}
  \end{picture}
}

\newcommand{\circarrow}{
  \setlength{\unitlength}{1mm}
  \begin{picture}(5,1)(0,0)
    \put(1,1){\circle{1}}
    \put(1.2,0){$\rightarrow$}
  \end{picture}
}

\newcommand{\arrowcirc}{
  \setlength{\unitlength}{1mm}
  \begin{picture}(5,1)(0,0)
    \put(0.2,0){$\leftarrow$}
    \put(4.3,1){\circle{1}}
  \end{picture}
}

\newcommand{\circcirc}{
  \setlength{\unitlength}{1mm}
  \begin{picture}(5,1)(0,0)
    \put(1,1){\circle{1}}
    \put(1.5,1){\line(1,0){2}}
    \put(4,1){\circle{1}}
  \end{picture}
}

\makeatletter
\def\ctext#1{\expandafter\@ctext\csname c@#1\endcsname}
\def\@ctext#1{\ifcase#1\or 0\or 1\or 2\or
3\or 4\or 5\fi}
\makeatother
\AddEnumerateCounter{\ctext}{\@ctext}{0}

\makeatletter
\def\btext#1{\expandafter\@btext\csname c@#1\endcsname}
\def\@btext#1{\ifcase#1\or B-i\or B-ii\or B-iii\fi}
\makeatother
\AddEnumerateCounter{\btext}{\@btext}{B-i}

\makeatletter
\def\cctext#1{\expandafter\@cctext\csname c@#1\endcsname}
\def\@cctext#1{\ifcase#1\or \textbf{Amenability}\or \textbf{Forbidden set}\or \textbf{Blocking}\fi}
\makeatother
\AddEnumerateCounter{\cctext}{\@cctext}{Amenability}

\makeatletter
\def\ccctext#1{\expandafter\@ccctext\csname c@#1\endcsname}
\def\@ccctext#1{\ifcase#1\or \textbf{Forbidden set}\or \textbf{Blocking}\fi}
\makeatother
\AddEnumerateCounter{\ccctext}{\@ccctext}{Amenability}

\makeatletter
\def\dtext#1{\expandafter\@dtext\csname c@#1\endcsname}
\def\@dtext#1{\ifcase#1\or \textbf{m-separation in~$\gpbd{\mathbf{XY}}$}\or \textbf{m-separation}\or \textbf{m-separation}\fi}
\makeatother
\AddEnumerateCounter{\dtext}{\@dtext}{m-separation}

\begin{document}

\title{Complete Graphical Characterization and Construction of Adjustment Sets in Markov Equivalence Classes of Ancestral Graphs}
\author{Emilija Perkovi\'c, Johannes Textor, Markus Kalisch and Marloes H. Maathuis}
\date{}

\maketitle

\begin{abstract}
We present a graphical criterion for covariate adjustment that is sound and complete for four
different classes of causal graphical models: directed acyclic graphs (DAGs), maximum ancestral graphs (MAGs),
completed partially directed
acyclic graphs (CPDAGs), and partial ancestral
graphs (PAGs). Our criterion unifies covariate adjustment for a large set of graph classes.
Moreover, we define an explicit set that satisfies our criterion, if there is any set that satisfies our criterion. We also give
efficient algorithms for constructing all sets that fulfill our criterion, implemented in the \texttt{R} package \texttt{dagitty}.
Finally, we discuss the relationship between our criterion and other criteria for adjustment, and we provide new soundness and completeness proofs
for the adjustment criterion for DAGs.
\end{abstract}

\noindent{}\textbf{Keywords:} causal effects, graphical models, covariate adjustment, latent variables, confounding

\section{Introduction} \label{sec:intro}

Covariate adjustment is a well-known method to estimate causal effects from observational data. 
There are, however, still common misconceptions about what variables one should or should not adjust for.
 For example,
it is sometimes thought that adjusting for more variables will lead to a more precise estimate, as long as the added variables are not affected by the exposure variable.  While this is true for a randomized exposure, in observational data even adjustment for pre-exposure variables may lead to so-called collider bias as described in the ``M-bias graph'' \citep{Shrier2008,Rubin2008}. 
Another example is the ``Table 2 fallacy'' \citep{Westreich2013}. In observational research papers, Table 1 often describes the data, and Table 2
shows a multiple regression analysis. By presenting all estimated coefficients in one table, it is implicitly suggested that all estimates
can be interpreted similarly. This is usually not the case: some coefficients may be interpreted as a total causal effect, some may be
interpreted as a direct causal effect, and some do not have any causal interpretation at all.

The practical importance of covariate adjustment has inspired
a growing body of theoretical work on graphical criteria for adjustment.
Pearl's back-door criterion \citep{pearl1993bayesian} is probably the most well-known, and is
sound but not complete for adjustment in $\DAG$s.
\citet{shpitser2012validity} and \citet{shpitser2012avalidity} refined the back-door criterion to a sound and complete graphical criterion
for adjustment in $\DAG$s.
Others considered more general graph classes, which can represent structural uncertainty.
\Citet{vanconstructing,van2018separators} gave sound and complete graphical criteria for $\MAG$s that
 allow for unobserved variables (latent confounding). \citet{maathuis2013generalized}
presented a generalized back-door criterion for $\DAG$s, $\CPDAG$s, $\MAG$s and $\PAG$s,
where $\CPDAG$s and $\PAG$s represent Markov equivalence classes of $\DAG$s or $\MAG$s, respectively,
and can be inferred directly from data (see, for example, \citealp{spirtes2000causation, Chickering02-optimal, Colombo2012, ClaassenMooijHeskes13, colombo2014order, nandyMaathuisArges, frot2017learning, heinze2017causal}). The generalized back-door criterion
is sound but not complete for adjustment. Another line of work explores data driven covariate adjustment that does not require knowing the graph \citep{vanderweele2011new,DeLunaEtAl11,EntnerHoyerSpirtes13}. Some of these data driven results are sound and complete for adjustment, but they all rely on some additional assumptions. We will not explore this direction in our paper.

In \cite{perkovic15_uai}, the preliminary conference version of the present paper,
we extended the results of
\citet{shpitser2012validity,shpitser2012avalidity}, \citet{vanconstructing} and \citet{maathuis2013generalized} to derive a single sound and complete adjustment criterion for $\DAG$s,
$\CPDAG$s, $\MAG$s and $\PAG$s. The different adjustment criteria are summarized in Table~\ref{table:relation}. Additionally, we note that \citet{vanDerZander16} showed that the generalized adjustment criterion can also be applied to the more general class of restricted chain graphs (representing a subset of a Markov equivalence class of $\DAG$s). Furthermore, in \cite{perkovic17}, we extend the generalized adjustment criterion to maximally oriented partially directed acyclic graphs (PDAGs), which represent $\CPDAG$s with added background knowledge. 
\begin{table}
\centering
 \begin{tabular}{| l | c | c | c | c |}
       \hline
     & DAG & MAG & CPDAG & PAG \\ \hline
     \specialcell{\small Back-door criterion \\ \small \cite{pearl1993bayesian} \vspace{0.05cm}} & {\Large $\Rightarrow$} & & &  \\ \hline
     \specialcell{\small Adjustment criterion \\ \small \cite{shpitser2012validity}, \cite{shpitser2012avalidity} \vspace{0.06cm}}  & {\Large $\Leftrightarrow$} & & &  \\ \hline
     \specialcell{\small Adjustment criterion \\ \small \cite{vanconstructing} \vspace{0.05cm}} & {\Large $\Leftrightarrow$} & {\Large $\Leftrightarrow$} & &\\ \hline
     \specialcell{\small Generalized back-door criterion \\  \small \cite{maathuis2013generalized} \vspace{0.05cm}} & {\Large $\Rightarrow$} & {\Large $\Rightarrow$} & {\Large $\Rightarrow$} & {\Large $\Rightarrow$} \\ \hline
     \specialcell{\vspace{0.02cm} \small \textbf{Generalized adjustment criterion} \\  \small \cite{perkovic15_uai} \vspace{0.05cm}}  & {\Large $\Leftrightarrow$} & {\Large $\Leftrightarrow$} & {\Large $\Leftrightarrow$} & {\Large $\Leftrightarrow$} \\ \hline
    \end{tabular}
    \caption{Graphical criteria for covariate adjustment:  {\Large $\Rightarrow$} - sound,
{\Large $\Leftrightarrow$} - \normalsize sound and complete.}
     \label{table:relation}

\end{table}

\begin{figure}
\tikzstyle{every edge}=[draw,>=stealth',->]
\newcommand\dagvariant[1]{\begin{tikzpicture}[xscale=.6,yscale=0.6,main node/.style={font=\small}]
\node[main node] (a) at (0,2) {A};
\node[main node] (b) at (2,2) {B};
\node[main node,inner sep=1pt] (i) at (-1,1) {I};
\node[main node,inner sep=1pt] (z) at (1,1) {Z};
\node[main node] (x) at (0,0) {X};
\node[main node] (y) at (2,0) {Y};
\draw (a) edge (x);
\draw (b) edge (y);
\draw (z) edge (x) edge (y);
\draw (x) edge (y);
\draw (i) edge (x);
\draw #1;
\end{tikzpicture}}
\hspace{1cm}
\begin{subfigure}{.3\textwidth}
\center
\begin{tikzpicture}[xscale=1,yscale=1]
\node (a) at (0,2) {A};
\node (b) at (2,2) {B};
\node (i) at (-1,1) {I};
\node (z) at (1,1) {Z};
\node (x) at (0,0) {X};
\node (y) at (2,0) {Y};
\draw (i) edge [o-o] (a) edge  (x);
\draw (z) edge [o-o] (a) edge [o-o] (b) (a) edge [o-o] (b);
\draw (a) edge (x);
\draw (b) edge (y);
\draw (z) edge (x) edge (y);
\draw (x) edge (y);
\end{tikzpicture}
\caption{}
\end{subfigure}
\vrule
\hspace{1cm}
\begin{subfigure}{.44\textwidth}
\dagvariant{
(i) edge [->] (a)
(z) edge [<-] (a) edge [<-] (b)
(a) edge [->] (b)
}
\dagvariant{
(i) edge [<-] (a)
(z) edge [<-] (a) edge [<-] (b) (a) edge [->] (b)
}
\dagvariant{
(i) edge [<-] (a)
(z) edge [<-] (a) edge [<-] (b) (a) edge [<-] (b)
}
\dagvariant{
(i) edge [->] (a)
(z) edge [<-] (a) edge [->] (b) (a) edge [->] (b)
}
\dagvariant{
(i) edge [<-] (a)
(z) edge [<-] (a) edge [->] (b) (a) edge [->] (b)
}
\dagvariant{
(i) edge [<-] (a)
(z) edge [->] (a) edge [<-] (b) (a) edge [<-] (b)
}
\dagvariant{
(i) edge [<-] (a)
(z) edge [->] (a) edge [->] (b) (a) edge [->] (b)
}
\dagvariant{
(i) edge [<-] (a)
(z) edge [->] (a) edge [->] (b) (a) edge [<-] (b)
}
\caption{}
\end{subfigure}
\caption{\small (a) A $\CPDAG$ in which, according to our
criterion, $\{A,Z\}$ is an  adjustment set
for the total causal effect of $X$ on $Y$. (b)
The Markov equivalence class of $\DAG$s represented by the $\CPDAG$. An adjustment
set for a $\CPDAG$ $(\PAG)$ is one that
is valid for all $\DAG$s ($\MAG$s) in the Markov equivalence class.}
\label{fig:cpdagexample}
\end{figure}
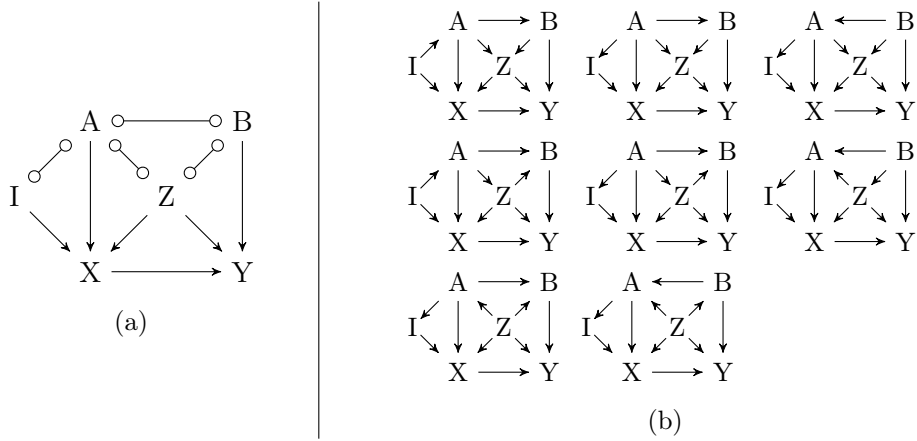

To illustrate the use of our generalized adjustment
criterion, suppose we are given the $\CPDAG$
in Figure~\ref{fig:cpdagexample}a and we want to
estimate the total causal effect of $X$ on $Y$.
Our criterion will inform us that the set $\{A,Z\}$ is an
adjustment set for this $\CPDAG$, meaning that it
is an adjustment set in every $\DAG$ that the $\CPDAG$
represents (Figure~\ref{fig:cpdagexample}b). Hence, we can estimate the causal
effect without knowledge of the full causal structure.
In a similar manner, by applying our criterion to a $\MAG$
or a $\PAG$, we find adjustment sets that are
valid for all $\DAG$s represented by this $\MAG$ or $\PAG$.
Our criterion finds such
adjustment sets whenever they exist;
 else, the causal effect is not identifiable by covariate adjustment. 
We hope that this ability to allow for incomplete
structural knowledge, latent confounding, or
both will help address concerns
that graphical causal modelling ``assumes that all [...] $\DAG$s
have been properly specified'' \citep{West2014}. Moreover, our criterion for
$\CPDAG$s and $\PAG$s can be combined with causal structure learning algorithms.

In the current paper, we give full proofs of the results in \cite{perkovic15_uai}. In addition, we provide several new results that allow us to construct sets $\mathbf{Z}$ that fulfill the generalized adjustment criterion for given sets $\mathbf{X}$ of exposures and $\mathbf{Y}$ of response variables in a $\DAG$, $\CPDAG$, $\MAG$ or $\PAG$ $\g$. 
In Corollary~\ref{cor:adjust xy gac general} we define a specific set that satisfies our criterion, if any set does. We refer to this set as a ``constructive set''. In Theorem~\ref{theorem:gac-alt}, we show how one can express adjustment sets in terms of $m$-separating sets in a certain subgraph of $\g$. This theorem reduces the problem of finding adjustment sets to the problem of finding $m$-separating sets, which has been studied in detail by \cite{vanconstructing}. In Lemma~\ref{lemmaequivalenceofcondb}, we prove that all adjustment sets for a $\CPDAG$ $(\PAG)$ $\g$ can be found in an arbitrary orientation of $\g$ into a valid $\DAG$ ($\MAG$). This allows us to leverage existing implementations \citep{vanconstructing}. We implemented the criterion itself and the construction of all adjustment sets in the software {\tt dagitty} \citep{textor2016robust}, available as a web-based GUI and an R package, and in the R package {\tt pcalg} \citep{kalischpcalg}.

Furthermore, we explore the relationships between our generalized adjustment criterion and the previously suggested generalized back-door criterion and Pearl's back-door criterion.
For both Pearl's back-door criterion and the generalized back-door criterion, a constructive set was given only in the case when the number of exposures is limited to one ($|\mathbf{X}|=1$).
We give constructive sets for each of these criteria for general $\mathbf{X}$ in Corollary~\ref{cor:constr-bc} and Corollary~\ref{cor:constr-set-gbc}.
Moreover, in Theorem~\ref{theorem:unified-nope} we identify cases in which there exist sets satisfying two, or all three of these criteria, as well as cases in which there are only sets satisfying the generalized adjustment criterion.

Another important contribution, included in the appendix, are new soundness and completeness proofs of the adjustment criterion for $\DAG$s as defined in \cite{shpitser2012validity} and in the unpublished addendum \cite{shpitser2012avalidity}, where the adjustment criterion in \cite{shpitser2012avalidity} is a revised version of the criterion in \cite{shpitser2012validity} (see Definition~\ref{def:johannes ac shpitser} in Appendix~\ref{subsec:proofs shpitser}). Since there are no published soundness and completeness proofs for the revised criterion and since we build on this work, we felt it was important to provide these proofs. The proofs are non-trivial, but rely only on elementary concepts.

We note that, although we can find all causal effects that are identifiable by covariate adjustment, we generally do not find
all identifiable causal effects, since some effects may be identifiable only by other means, using for example IDA approaches \citep{MaathuisKalischBuehlmann09, MaathuisColomboKalischBuehlmann10, nandy2014estimating, malinsky2017estimating}, 
Pearl's front-door criterion \citep[Section 3.3.2]{Pearl2009} or the ID algorithm \citep{tian2002general,shpitser2006identification}.

We also point out that $\MAG$s and $\PAG$s are in principle
not only able to represent unobserved confounding, but
can also account for unobserved selection variables.
In this paper, however, we assume that there are no
unobserved selection variables, since
selection bias often rules out causal effect identification using
just covariate adjustment. \citet{Bareinboim2014} discuss
these problems and present creative approaches to work
around them, for example by combining data from different sources.
The question whether our adjustment criterion
could be combined with such auxiliary methods
is left for future research.

\section{Preliminaries} \label{sec:prelim}

Throughout the paper we denote sets in bold (for example $\mathbf{X}$), graphs in calligraphic font (for example $\g$) and nodes in a graph in uppercase letters (for example $X$). All omitted proofs are given in the appendix.

\textbf{Nodes and edges.} A graph $\g= (\vars,\e) $ consists of a set of nodes (variables) $ \vars=\left\lbrace X_{1},\dots,X_{p}\right\rbrace$ and a set of edges $ \e $. We consider simple graphs, meaning that there is at most one edge between any pair of nodes. Two nodes are called \textit{adjacent} if they are connected by an edge. Every edge has two edge marks that can be arrowheads, tails or circles. Edges can be \emph{directed} $\rightarrow$,  \emph{bi-directed}  $\leftrightarrow$, \emph{non-directed} $\circcirc$, or \emph{partially directed} $\circarrow$. We use $\bullet$ as a stand in for any of the allowed edge marks. An edge is \textit{into} (\textit{out of}) a node $X$ if the edge has an arrowhead (tail) at $X$. A \emph{directed graph} contains only directed edges. A \emph{mixed graph} may contain directed and bi-directed edges. A \emph{partial mixed graph} may contain any of the described edges. Unless stated otherwise, definitions apply to partial mixed graphs.

\textbf{Paths.} A \textit{path} $p$ from $X$ to $Y$ in~$\g$ is a sequence of distinct nodes $\langle X, \dots,Y \rangle$ in which every pair of successive nodes is adjacent in~$\g$. If $p = \langle X_1, X_2, \dots , X_k, \rangle, k \ge 2$, then with $-p$ we denote the path $\langle X_k, \dots , X_2, X_1 \rangle$. A node $V$ \emph{lies on a path} $p$ if $V$ occurs in the sequence of nodes. If $p = \langle X_1, X_2, \dots , X_k, \rangle, k \ge 2$, then $X_1$ and $X_k$ are \textit{endpoints} of $p$, and any other node $X_i, 1 <i<k,$ is a \textit{non-endpoint} node on $p$.
The \textit{length} of a path equals the number of edges on the path. A \textit{directed path} from $X$ to $Y$ is a path from $X$ to $Y$ in which all edges are directed towards $Y$, that is, $X \to \dots\to Y$. We also refer to this as a \textit{causal path}. A \textit{possibly directed path} or \textit{possibly causal path} from $X$ to $Y$ is a path from $X$ to $Y$ that does not contain an arrowhead pointing in the direction of $X$. 
A path from $X$ to $Y$ that is not possibly causal is called a \textit{non-causal path} from $X$ to $Y$. A directed path from $X$ to $Y$ together with $Y\to X$ forms a \emph{directed cycle}. A directed path from $X$ to $Y$ together with $Y\leftrightarrow X$ forms an \emph{almost directed cycle}.
For two disjoint subsets $\mathbf{X}$ and $\mathbf{Y}$ of $\mathbf{V}$, a path from $\mathbf{X}$ to $\mathbf{Y}$
is a path from some $X \in \mathbf{X}$ to some $Y \in \mathbf{Y}$.
A path from $\mathbf{X}$ to $\mathbf{Y}$ is \textit{proper} (wrt $\mathbf{X}$) if only its first node is in $\mathbf{X}$. If $\g$ and $\g^*$ are two graphs with identical adjacencies and $p$ is a path in $\g$, then the \textit{corresponding path} $\pstar$ is the path in $\g^*$ constituted by the same sequence of nodes as $p$.

\textbf{Subsequences, subpaths and concatenation.} A \textit{subsequence} of a path $p$ is a sequence of nodes obtained by deleting some nodes from $p$ without changing the order of the remaining nodes. A subsequence of a path is not necessarily a path. For a path $p = \langle X_1,X_2,\dots,X_m \rangle$, the \textit{subpath} from $X_i$ to $X_k$ ($1\le i\le k\le m)$ is the path $p(X_i,X_k) = \langle X_i,X_{i+1},\dots,X_{k}\rangle$. We denote the concatenation of paths by $\oplus$, so
that for example $p = p(X_1,X_{k}) \oplus p(X_{k},X_{m})$. In this paper, we only concatenate paths if the result of the concatenation is again a path. 

\textbf{Ancestral relationships.} If $X\to Y$, then $X$ is a \textit{parent} of $Y$. If there is a directed (possibly directed) path from $X$ to $Y$, then $X$ is a \textit{ancestor} (\textit{possible ancestor}) of $Y$, and $Y$ is a \textit{descendant} (\textit{possible descendant}) of $X$. We also use the convention that every node is a descendant, possible descendant, ancestor and possible ancestor of itself. The sets of parents, descendants and ancestors of $X$ in~$\g$ are denoted by $\Pa(X,\g)$, $\De(X,\g)$ and $\An(X,\g)$ respectively. The sets of possible descendants and possible ancestors of $X$ in $\g$ are denoted by $\PossDe(X,\g)$ and $\PossAn(X,\g)$ respectively. For a set of nodes $\mathbf{X} \subseteq \mathbf{V}$, we let $\Pa(\mathbf{X},\g) = \cup_{X \in \mathbf{X}} \Pa(X,\g)$, with analogous definitions for $\De(\mathbf{X},\g)$, $\An(\mathbf{X},\g)$, $\PossDe(\mathbf{X},\g)$ and $\PossAn(\mathbf{X},\g)$.

\textbf{Colliders, shields and definite status paths.} If a path $p$ contains $X_i \bulletarrow X_j \arrowbullet X_k$ as a subpath, then $X_j$ is a \textit{collider} on $p$. A \textit{collider path} is a path on which every non-endpoint node is a collider. A  path of length one is a trivial collider path. A path $\langle X_{i},X_{j},X_{k} \rangle$ is an \emph{(un)shielded triple} if $ X_{i} $ and $ X_{k}$ are (not) adjacent. A path is \textit{unshielded} if all successive triples on the path are unshielded. A node $X_{j}$ is a \textit{definite non-collider} \citep{zhang2008causal} on a path $p$ if there is at least one edge out of $X_{j}$ on $p$, or if $X_{i} \bulletcirc X_j \circbullet X_k$ is a subpath of $p$ and $\langle X_i,X_j,X_k\rangle$ is an unshielded triple. Any collider on a path is always of definite status and hence, a \textit{definite collider}. In a $\DAG$ ($\MAG$) we refer to definite non-colliders as \emph{non-colliders}. A node is of \textit{definite status} on a path if it is a collider or a definite non-collider on the path. A path $p$ is of definite status if every non-endpoint node on $p$ is of definite status.

\textbf{m-separation and m-connection.} A definite status path \textit{p} between nodes $X$ and $Y$ is \textit{m-connecting} given a set of nodes $\mathbf{Z}$ ($X,Y \notin \mathbf{Z}$) if every definite non-collider on $p$ is not in $\mathbf{Z}$, and every collider on $p$ has a descendant in $\mathbf{Z}$ \citep{richardson2003markov}. Otherwise $\mathbf{Z}$ blocks $p$. If $\g$ is a $\DAG$ or $\MAG$ (defined later) and if $\mathbf{Z}$ blocks all paths between $X$ and $Y$, we say that $X$ and $Y$ are m-separated given $\mathbf{Z}$ in $\g$. Otherwise, $X$ and $Y$ are m-connected given $\mathbf{Z}$ in $\g$. For pairwise disjoint subsets $\mathbf{X}$, $\mathbf{Y}$ and $\mathbf{Z}$ of $\vars$ in $\g$, we say that $\mathbf{X}$ and $\mathbf{Y}$ are m-separated given $\mathbf{Z}$ in $\g$ if $X$ and $Y$ are m-separated given $\mathbf{Z}$ in $\g$ for any $X \in \mathbf{X}$ and $Y\in \mathbf{Y}$. Otherwise, $\mathbf{X}$ and $\mathbf{Y}$ are m-connected given $\mathbf{Z}$ in $\g$. In a $\DAG$, m-separation and m-connection simplify to d-separation and d-connection \citep[][]{Pearl2009}.

\textbf{Causal Bayesian networks.} A directed graph without directed cycles is a \emph{directed acyclic graph $(\DAG)$}. A Bayesian network for a set of variables $\vars=\{X_{1},\dots,X_{p}\}$ is a pair ($\g,f$), where $\g$ is a $\DAG$, and $f$ is a joint density for $\vars$ that factorizes as $f(\vars)= \prod_{i=1}^{p}f(X_{i}|\Pa(X_{i},\g))$ \citep{Pearl2009}. We call a $\DAG$ \emph{causal} if every edge $X_{i} \rightarrow X_{j}$ in~$\g$ represents a direct causal effect of $X_{i}$ on $X_{j}$. A Bayesian network ($\g,f$) is a \emph{causal Bayesian network} if $\g$ is a causal $\DAG$. If a causal Bayesian network is given and all variables are observed, one can easily derive post-intervention densities. In particular, we consider interventions $do(\mathbf{X=x})$, or shorthand $do(\mathbf{x})$, ($\mathbf{X}\subseteq \mathbf{V}$), which represent outside interventions that set $\mathbf{X}$ to $\mathbf{x}$, uniformly in the population \citep[see][]{Pearl2009}:
\begin{equation} \label{eq:fact-form}
f(\mathbf{v}\mid do(\mathbf{x})) =
\begin{cases}
\prod_{\{i \mid X_{i} \in \vars \setminus \mathbf{X}\}}f(x_{i}\mid\Pa(x_{i},\g)), &  \text{if }\mathbf{v} \text{ is consistent with }\mathbf{x,} \\
0, & \text{otherwise.}
\end{cases}
\end{equation}
Equation \eqref{eq:fact-form} is known as the truncated factorization formula \citep{Pearl2009}, the g-formula \citep{robins1986new} or the manipulated density \citep{spirtes2000causation}.

\textbf{Maximal ancestral graphs.}
A mixed graph $\g$ without directed cycles and almost directed cycles is called \textit{ancestral}. A \emph{maximal ancestral graph $(\MAG)$} is an ancestral graph $\g = (\vars, \e)$ where every pair of non-adjacent nodes $X$ and $Y$ in~$\g$ can be m-separated by a set $\mathbf{Z} \subseteq \vars \setminus \{X,Y\}$. A $\DAG$ with unobserved variables can be uniquely represented by a $\MAG$ on the observed variables that preserves the ancestral and m-separation relationships among the observed variables \citep[page~981~in][]{richardson2002ancestral}.  
Since we consider $\MAG$s that do not encode selection bias, the $\MAG$s in this paper can only contain directed ($\rightarrow$) and bi-directed ($\leftrightarrow$) edges.
The $\MAG$ of a causal $\DAG$ is a \emph{causal $\MAG$}.

\textbf{Markov equivalence.}
 Several $\DAG$s can encode the same conditional independencies via d-separation. Such $\DAG$s form a \emph{Markov equivalence class} which can be described uniquely by a \emph{completed partially directed acyclic graph $(\CPDAG)$} \citep{meek1995causal}. A $\CPDAG$ $\g[C]$ has the same adjacencies as any $\DAG$ in the Markov equivalence class described by $\g[C]$. A directed edge $X \rightarrow Y$ in a $\CPDAG$ $\g[C]$ corresponds to a directed edge $X \rightarrow Y$ in every $\DAG$ in the Markov equivalence class described by $\g[C]$. For any non-directed edge $X \circcirc Y$ in a $\CPDAG$ $\g[C]$, the Markov equivalence class described by $\g[C]$ contains a $\DAG$ with $X \rightarrow Y$ and a $\DAG$ with $X \leftarrow Y$\footnote{The non-directed edges in a $\CPDAG$, which we denote as $\circcirc$, are often denoted as $-$ in the relevant literature, see for example \citep{meek1995causal}. We use $\circcirc$ instead of $-$ for the sake of consistency among different graph classes.}. 
 Thus, $\CPDAG$s only contain directed ($\rightarrow$) and non-directed ($\circcirc$) edges.

Several $\MAG$s can also encode the same conditional independencies via m-separation. Such $\MAG$s form a Markov equivalence class which can be described uniquely by a \emph{partial ancestral graph $(\PAG)$} \citep{richardson2002ancestral, ali2009markov}. A $\PAG$ $\g[P]$ has the same adjacencies as any $\MAG$ in the Markov equivalence class described by $\g[P]$. Any non-circle edge-mark in a $\PAG$ $\g[P]$ corresponds to that same non-circle edge-mark in every $\MAG$ in the Markov equivalence class described by $\g[P]$. We only consider maximally informative $\PAG$s \citep{zhang2008completeness}, that is, for any circle mark $X \circbullet $ in a $\PAG$ $\g[P]$, the Markov equivalence class described by $\g[P]$ contains a $\MAG$ with $X \arrowbullet Y$ and a $\MAG$ with $X \tailbullet Y$. We denote all $\DAG$s ($\MAG$s) in the Markov equivalence class described by a $\CPDAG$ $(\PAG)$ $\g$ by $[\g]$. The $\CPDAG$ ($\PAG$) of a causal $\DAG$ ($\MAG$) is a \emph{causal} $\CPDAG$ ($\PAG$).

\textbf{Consistent densities.} A density $f$ is \textit{consistent} with a causal $\DAG$ $\g[D]$ if the pair ($\g[D],f$) forms a causal Bayesian network. A density $f$ is consistent with a causal $\MAG$ $\g[M]$ if there exists a causal Bayesian network ($\g[D]',f'$) such that $\g[M]$ represents $\g[D]'$ and $f$ is the observed marginal of $f'$. A density $f$ is consistent with a causal $\CPDAG$ $(\PAG)$ $\g$ if it is consistent with a causal $\DAG$ $(\MAG)$ in $[\g]$.

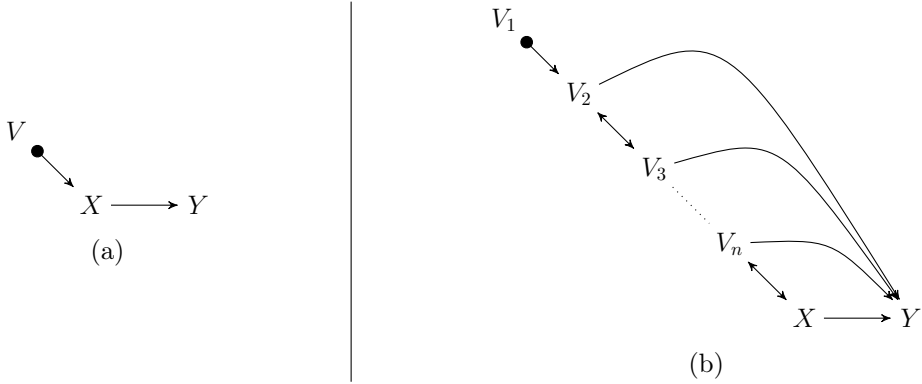
\begin{figure}
\centering
\begin{subfigure}{.25\textwidth}
  \centering
  \begin{tikzpicture}[>=stealth',shorten >=1pt,auto,node distance=2cm,main node/.style={minimum size=0.6cm,font=\sffamily\Large\bfseries},scale=0.7,transform shape]
\node[main node]         (X)                        {$X$};
\node[main node]         (V) [above left of= X]  	{$V$};
\node[main node]       	 (Y) [right of= X]         	{$Y$};
\draw[*->] (V) edge    (X);
\draw[->] (X) edge    (Y);
\end{tikzpicture}
  \caption{}
  \label{visible:sub1}
\end{subfigure}
\hspace{1cm}
\vrule
\begin{subfigure}{.6\textwidth}
  \centering
  \begin{tikzpicture}[>=stealth',shorten >=1pt,auto,node distance=2cm,main node/.style={minimum size=0.6cm,font=\sffamily\Large\bfseries},scale=.7,transform shape]
\node[main node]         (X)                        {$X$};
\node[main node]         (Vn) [above left of= X]  	{$V_n$};
\node[main node]         (V3) [above left of= Vn]  {$V_3$};
\node[main node]         (V2) [above left of= V3]  	{$V_2$};
\node[main node]         (V1) [above left of= V2]  	{$V_1$};
\node[main node]       	 (Y) [right of= X]         	{$Y$};
\draw[*->] (V1) edge    (V2);
\draw[<->] (V2) edge    (V3);
\draw[dotted] (V3) edge    (Vn);
\draw[<->] (Vn) edge    (X);
\draw[->] (X) edge    (Y);
\draw[->] (Vn) ..  controls (0.5,1.5) ..  (Y);
\draw[->] (V3) ..  controls (-0.5,3.5) ..  (Y);
\draw[->] (V2) ..  controls (-1.5,5.6) ..  (Y);
\end{tikzpicture}
  \caption{}
  \label{visible:sub2}
\end{subfigure}
\caption{\small Two configurations where the edge $X \rightarrow Y$ is visible. Nodes $V$ and $Y$ must be nonadjacent in~\ref{visible:sub1}, and $V_1$ and $Y$ must be nonadjacent in~\ref{visible:sub2}. }
\label{fig:visible}
\end{figure}

\textbf{Visible and invisible edges.} All directed edges in $\DAG$s and $\CPDAG$s are said to be visible. Given a $\MAG$ $\g[M]$ or a $\PAG$ $\g$, a directed edge $X \rightarrow Y$ is \textit{visible} if there is a node $V$ not adjacent to $Y$ such that there is an edge $V \bulletarrow X$, or if there is a collider path between $V$ and $X$ that is into $X$ and every non-endpoint node on the path is a parent of $Y$, see Figure~\ref{fig:visible} \citep{zhang2006causal}. A visible edge $X\rightarrow Y$ means that there are no latent confounders between $X$ and $Y$. A directed edge $X \rightarrow Y$ that is not visible in a $\MAG$ $\g[M]$ or a $\PAG$ $\g$ is said to be \textit{invisible}. In the FCI algorithm, invisible edges can occur due to orientation rules \textbf{R5} - \textbf{R10} of \cite{zhang2008completeness}.  
When considering $\MAG$s and $\PAG$s that do not encode selection bias, invisible edges occur as a consequence of the orientation rules \textbf{R8} - \textbf{R10} of \cite{zhang2008completeness}. 

\section{The Generalized Adjustment Criterion} \label{sec:main}

Throughout, let $\g=(\vars,\e)$ represent a $\DAG$, $\CPDAG$, $\MAG$ or $\PAG$, and let $\mathbf{X}$, $\mathbf{Y}$ and $\mathbf{Z}$ be pairwise disjoint subsets of $\vars$, with $\mathbf{X}\neq \emptyset$ and $\mathbf{Y}\neq \emptyset$. Here, $\mathbf{X}$ represents the set of exposures and $\mathbf{Y}$ represents the set of response variables. 

We will define sound and complete graphical conditions for adjustment sets relative to $(\mathbf{X,Y})$ in~$\g$. Thus, if a set $\mathbf{Z}$ satisfies our conditions  relative to ($\mathbf{X},\mathbf{Y}$) in~$\g$ (see Definition~\ref{def:gac}), then it is a valid adjustment set for calculating the causal effect of $\mathbf{X}$ on $\mathbf{Y}$ (see Definition~\ref{defadjustment}), and every existing valid adjustment set satisfies our conditions (see Theorem~\ref{theorem:gac}).
First, we define what we mean by an adjustment set.
\begin{definition}{(\textbf{Adjustment set}; \citealp{maathuis2013generalized})}
   Let $\mathbf{X,Y}$ and $\mathbf{Z}$ be pairwise disjoint node sets in a causal $\DAG$, $\CPDAG$, $\MAG$ or $\PAG$ $\g$. Then $\mathbf{Z}$ is an adjustment set relative to $(\mathbf{X,Y})$ in~$\g$  if for any density\footnote{We use the notation for continuous random variables throughout. The discrete analogues should be obvious.} $f$ consistent with $\g$ we have
   \begin{equation}
   f(\mathbf{y}\mid do(\mathbf{x}))=
   \begin{cases}
   f(\mathbf{y}\mid\mathbf{x}) & \text{if }\mathbf{Z} = \emptyset,\\
   \int_{\mathbf{z}}f(\mathbf{y}\mid\mathbf{x,z})f(\mathbf{z})d\mathbf{z}  & \text{otherwise.}
   \end{cases}
   \label{lab}
   \end{equation}
   \label{defadjustment}
\end{definition}
\noindent{}Thus, adjustment sets allow post-intervention densities involving the do-operator (left-hand side of Equation \ref{lab}) 
to be identified as specific functions of conditional densities (right-hand side of Equation \ref{lab}). The latter can be estimated from observational data. As a result, adjustment sets are important for the computation of causal effects. This is illustrated in Example~\ref{ex:gaussian-linear} for the special case of multivariate Gaussian densities.

\begin{example}
Suppose $f$ is a multivariate Gaussian density that is consistent with a causal $\DAG$ $\g[D]$. Let $\mathbf{Z} \neq \emptyset$ be an adjustment set relative to two distinct variables $X$ and $Y$ in~$\g[D]$ such that $\mathbf{Z} \cap \{ X \cup Y \} = \emptyset$. Then
\begin{align}
&E(Y \given do(x)) =  \int_{y} y f(y \given do(x))dy  =   \int_{y} y \int_\mathbf{z} f(y \given x,\mathbf{z})f(\mathbf{z})d\mathbf{z}dy  \nonumber \\
& =  \int_\mathbf{z} \int_{y} y f(y \given x,\mathbf{z})dyf(\mathbf{z})d\mathbf{z} = \int_\mathbf{z} E(Y \given x,\mathbf{z})f(\mathbf{z})d\mathbf{z}  \nonumber \\
&=  \int_\mathbf{z}(\alpha + \gamma x + \beta^T \mathbf{z})f(\mathbf{z})d\mathbf{z}=\alpha + \gamma x + \beta^T E(\mathbf{Z}), \nonumber
\end{align}

\noindent{}where we use the fact that all conditional expectations in a multivariate Gaussian distribution are linear, so that $E(Y \given x,\mathbf{z}) = \alpha + \gamma x + \beta^{T} \mathbf{z}$, for some $\alpha, \gamma \in \mathbb{R}$ and $\beta \in \mathbb{R}^{|\mathbf{z}|}$. Defining the total causal effect of $X$ on $Y$ as $\frac{\partial }{\partial x}E(Y \given do(x))$, we obtain that the total causal effect of $X$ on $Y$ is $\gamma$, that is, the regression coefficient of $X$ in the regression of $Y$ on $X$ and $\mathbf{Z}$.
\label{ex:gaussian-linear}
\end{example}

Our first goal in this paper is to give a graphical criterion (see Definition~\ref{def:gac}) that is equivalent to Definition~\ref{defadjustment}. To this end, we introduce some additional terminology.

\begin{definition}{(\textbf{Amenability})}
   Let $\mathbf{X}$ and $\mathbf{Y}$ be disjoint node sets in a $\DAG$, $\CPDAG$, $\MAG$ or $\PAG$ $\g$. Then $\g$ is said to be amenable relative to $(\mathbf{X},\mathbf{Y})$ if every proper possibly directed path from $\mathbf{X}$ to $\mathbf{Y}$ in~$\g$ starts with a visible edge out of $\mathbf{X}$.
  \label{defadjamengen}
\end{definition}

\noindent{}If $\g$ is a $\MAG$, then Definition \ref{defadjamengen} reduces to the notion of amenability as introduced in \cite{vanconstructing}. The intuition behind the concept of amenability is the following: In $\MAG$s and $\PAG$s, directed edges $X\to Y$ can represent causal effects, but also mixtures of causal effects and latent confounding. For instance, when the graph $X \to Y$ is interpreted as a $\DAG$, the empty set is a valid adjustment set with respect to $(X,Y)$. When the same graph is interpreted as a $\MAG$, it can still represent the $\DAG$ $X \to Y$, but also the $\DAG$ $X\to Y$ with an additional non-causal path $X \leftarrow L \to Y$ where $L$ is latent.

In $\CPDAG$s and $\PAG$s, there are edges with unknown direction. This complicates adjustment because paths containing such edges can correspond to causal paths in some represented $\DAG$s and to non-causal paths in others.
For example, the $\CPDAG$ $X \circcirc Y$ represents the $\DAG$s $X \rightarrow Y$ and $X \leftarrow Y$. Amenable graphs are graphs where these problems do not occur.
\begin{definition}(\textbf{Forbidden set; $\fb{\g}$})  Let $\mathbf{X}$ and $\mathbf{Y}$ be disjoint node sets in a $\DAG$, $\CPDAG$, $\MAG$ or $\PAG$ $\g$. Then the forbidden set relative to $(\mathbf{X,Y})$ is defined as
\begin{align}
   \fb{\g} = \{ & W' \in \vars: W' \in \PossDe(W,\g) \text{, for some } W \notin \mathbf{X} \, \notag \\
   & \text{which lies on a proper possibly} \text{ directed path from } \, \mathbf{X} \,\text{to}\, \mathbf{Y}\text{in } \g\}.\notag
\end{align} \label{def:forbidden nodes}
\end{definition}
\begin{definition}{(\textbf{Generalized adjustment criterion})}
   Let $\mathbf{X,Y}$ and $\mathbf{Z}$ be pairwise disjoint node sets in a $\DAG$, $\CPDAG$, $\MAG$ or $\PAG$ $\g$. Then $\mathbf{Z}$ satisfies the generalized adjustment criterion
      relative to $(\mathbf{X,Y})$ in~$\g$ if the following three conditions hold:
   \begin{enumerate}[label = (\cctext*), leftmargin=0.5cm,align=left]
   \item\label{cond0} $\g$ is adjustment amenable relative to $(\mathbf{X,Y})$, and
   \item\label{cond1} $\mathbf{Z} \cap \fb{\g} = \emptyset$, and
   \item\label{cond2} all proper definite status non-causal paths from $\mathbf{X}$ to $\mathbf{Y}$ are blocked by $\mathbf{Z}$~in~$\g$.
   \end{enumerate}  \label{def:gac}
\end{definition}

\noindent{}If $\g$ is a $\DAG$ ($\MAG$), our criterion reduces to the adjustment criterion of \cite{shpitser2012avalidity}, \citep{vanconstructing} (see Definition~\ref{def:johannes ac shpitser}~in~Appendix~\ref{subsec:proofs shpitser}). For consistency, however, we will refer to the generalized adjustment criterion for all graph types.

We note that \amen{} does not depend on $\mathbf{Z}$. In other words, if \amen{} is violated, then no set satisfies the generalized adjustment criterion relative to $(\mathbf{X,Y})$ in~$\g$. The forbidden set contains nodes that cannot be used for adjustment. We will try to give some intuition. For simplicity, we consider $\mathbf{X} = \{X\}$ and $\mathbf{Y} = \{Y\}$ in a $\DAG$ $\g[D]$, and we are interested in estimating the total causal effect of $X$ on $Y$ in $\g[D]$. It is clear that nodes on any causal path from $X$ to $Y$  in $\g[D]$ should not be included in the set used for adjustment, since including such nodes would block the causal path.

To understand why we cannot include descendants of nodes on a causal path from $X$ to $Y$ in $\g[D]$ (except for descendants of $X$), it is useful to consider \textit{walks} from $X$ to $Y$ in $\g[D]$, that is, sequences of nodes $\langle X, \dots,Y \rangle$ in which every pair of successive nodes is adjacent in~$\g[D]$ (but the nodes are not necessarily distinct). A walk $r = \langle  X =V_0,V_1,\dots,V_k=Y \rangle$ is non-causal if $V_i \arrowbullet V_{i+1}$ for at least one $i \in \{0,\dots,k-1\}$. A walk $r$ from $X$ to $Y$ in $\g[D]$ is \textit{connecting} given a set of nodes $\mathbf{Z}$ if $\mathbf{Z}$ contains all colliders on $r$ and no non-collider on $r$ is in $\mathbf{Z}$. If a walk $r$ is not connecting given $\mathbf{Z}$, then $r$ is \textit{blocked} by $\mathbf{Z}$. \cite{koster2002marginalizing} proved that there is a walk from $X$ to $Y$ that is connecting given $\mathbf{Z}$ in $\g[D]$ if and only if there is path from $X$ to $Y$ that is d-connecting given $\mathbf{Z}$ in $\g[D]$.

Intuitively, all non-causal walks from $X$ to $Y$ should be blocked in order to estimate the total causal effect of $X$ on $Y$. Now, consider a path $p$ of the form $X \rightarrow V_1 \rightarrow \dots \rightarrow V_k \rightarrow Y$ in $\g$. Assume $V_i \notin \mathbf{Z}$, for all $i \in \{1,\dots,k\}$. Including a descendant $A$ of $V_i$ in the set $\mathbf{Z}$ leads to walk of the form $X \rightarrow \dots \rightarrow V_i \rightarrow \dots \rightarrow A \leftarrow \dots \leftarrow V_i \rightarrow \dots \rightarrow V_k \rightarrow Y$ being connecting given $\mathbf{Z}$ in $\g$. Hence, including $A$ in the adjustment set opens a non-causal walk from $X$ to $Y$ in $\g$. 

We now give the main theorem of this section. Corresponding examples can be found in Section~\ref{subsec:examples1} and the proof of the theorem is given in Section~\ref{subsec:proof-theorem-gac}.

\begin{theorem}
     Let $\mathbf{X,Y}$ and $\mathbf{Z}$ be pairwise disjoint node sets in a causal $\DAG$, $\CPDAG$, $\MAG$ or $\PAG$ $\g$. Then $\mathbf{Z}$ is an adjustment set relative to $(\mathbf{X},\mathbf{Y})$ in~$\g$ (see Definition~\ref{defadjustment}) if and only if $\mathbf{Z}$ satisfies the generalized adjustment criterion relative to $(\mathbf{X},\mathbf{Y})$ in~$\g$ (see Definition~\ref{def:gac}).
   \label{theorem:gac}
\end{theorem}

Verifying \blck by checking all paths requires keeping track of which paths are non-causal and hence, scales poorly to larger graphs. We therefore give an alternative definition of this condition which relies on m-separation in a so-called proper back-door graph.

\begin{definition}(\textbf{Proper back-door graph}; $\gpbd{\mathbf{XY}}$)
  Let $\mathbf{X}$ and $\mathbf{Y}$ be disjoint node sets in a $\DAG$, $\CPDAG$, $\MAG$ or $\PAG$ $\g$. The proper back-door graph $\gpbd{\mathbf{XY}}$ is obtained from $\g$ by removing all visible edges out of $\mathbf{X}$ that are on proper possibly directed paths from $\mathbf{X}$ to $\mathbf{Y}$ in~$\g$.
\label{def:gpbd}
\end{definition}

\noindent{}If $\g$ is a $\DAG$ or $\MAG$, then Definition~\ref{def:gpbd} reduces to the definition of proper back-door graphs as introduced in \cite{vanconstructing}.

\begin{theorem}
Replacing the Blocking condition in Definition~\ref{def:gac} with:

\begin{enumerate}[label = (\roman*), leftmargin=0.5cm,align=left]
\item[({\bfseries Separation})] \label{cond2prime} $\mathbf{Z}$ m-separates $\mathbf{X}$ and $\mathbf{Y}$ in~$\gpbd{\mathbf{XY}}$,
\end{enumerate}
results in a criterion that is equivalent to the generalized adjustment criterion.
\label{theorem:gac-alt}
\end{theorem}

\noindent{}If $\g$ is a $\DAG$ or $\MAG$, then Theorem~\ref{theorem:gac-alt} reduces to Theorem~4.6 in \cite{vanconstructing}.

\subsection{Examples} \label{subsec:examples1}

We now provide some examples that illustrate how the generalized adjustment criterion can be applied.

\begin{example}
   We return to the $\CPDAG$ $\g[C]$ in Figure~\ref{fig:cpdagexample}(a). $\g[C]$ is amenable relative to $(X,Y)$ and $\f{\g[C]} = \{ Y \}$.
One can easily verify that any superset of $\{Z,A\}$ or of $\{Z,B\}$ that does not contain $X$ or $Y$ satisfies the generalized adjustment criterion relative to $(X,Y)$ in~$\g[C]$.
  \label{ex:everything-works}
\end{example}

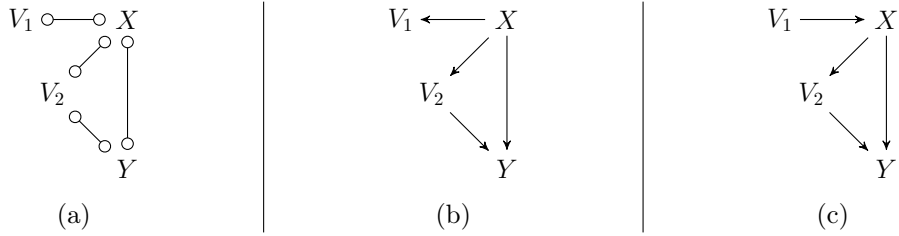
\begin{figure}
   \centering
   \begin{subfigure}{.32\textwidth}
     \centering
     \begin{tikzpicture}[>=stealth',shorten >=1pt,auto,node distance=2cm,main node/.style={minimum size=0.6cm,font=\sffamily\Large\bfseries},scale=0.7,transform shape]
     \node[main node]         (X)                        {$X$};
   \node[main node]         (V1) [left of= X]  		{$V_{1}$};
   \node[main node]         (V2) [below left of = X] 	{$V_{2}$};
   \node[main node]       	 (Y)  [below right of= V2] 	{$Y$};
   \draw[o-o] (V2) edge    (X);
   \draw[o-o] (V1) edge    (X);
   \draw[o-o] (X) edge    (Y);
   \draw[o-o] (V2) edge    (Y);
   \end{tikzpicture}
     \caption{}
   \end{subfigure}
   \unskip
   \vrule
   \begin{subfigure}{.32\textwidth}
     \centering
     \begin{tikzpicture}[>=stealth',shorten >=1pt,auto,node distance=2cm,main node/.style={minimum size=0.6cm,font=\sffamily\Large\bfseries},scale=.7,transform shape]
   \node[main node]         (X)                        {$X$};
   \node[main node]         (V1) [left of= X]  		{$V_{1}$};
   \node[main node]         (V2) [below left of = X] 	{$V_{2}$};
   \node[main node]       	 (Y)  [below right of= V2] 	{$Y$};
   \draw[<-] (V2) edge    (X);
   \draw[<-] (V1) edge    (X);
   \draw[->] (X) edge    (Y);
   \draw[->] (V2) edge    (Y);
   \end{tikzpicture}
     \caption{}
   \end{subfigure}
   \unskip
   \vrule
   \begin{subfigure}{.32\textwidth}
     \centering
     \begin{tikzpicture}[>=stealth',shorten >=1pt,auto,node distance=2cm,main node/.style={minimum size=0.6cm,font=\sffamily\Large\bfseries},scale=.7,transform shape]
   \node[main node]         (X)                        {$X$};
   \node[main node]         (V1) [left of= X]  		{$V_{1}$};
   \node[main node]         (V2) [below left of = X] 	{$V_{2}$};
   \node[main node]       	 (Y)  [below right of= V2] 	{$Y$};
   \draw[<-] (V2) edge    (X);
   \draw[->] (V1) edge    (X);
   \draw[->] (X) edge    (Y);
   \draw[->] (V2) edge    (Y);
   \end{tikzpicture}
     \caption{}
   \end{subfigure}
   \caption{\small (a) $\PAG$ $\g[P]$, (b) $\MAG$ $\g[M]_1$, (c) $\MAG$ $\g[M]_2$ used in Example~\ref{ex2}.}
   \label{figex2}
\end{figure}

\begin{example}
   To illustrate the concept of amenability, consider Figure~\ref{figex2} with a $\PAG$ $\g[P]$ in (a), and two $\MAG$s $\g[M]_1$ and $\g[M]_2$ in $[\g[P]]$ in (b) and (c). The graphs $\g[P]$ and $\g[M]_1$ are not amenable relative to $(X,Y)$. For $\g[P]$ this is due to the path $X \circcirc Y$, and for $\g[M]_1$ this is due to the invisible edge $X \to Y$  
   (which implies that $\g[M]_1$ also represents a DAG that contains a hidden confounder that is an ancestor of both $X$ and $Y$). On the other hand, $\g[M]_2$ is amenable relative to $(X,Y)$, since the edges $X\to Y$ and $X \to V_2$ are visible due to the edge $V_1\to X$, with $V_1$ not adjacent to $Y$ or $V_2$.
   Since there are no proper definite status non-causal paths from $X$ to $Y$ in~$\g[M]_2$, it follows that the empty set satisfies the generalized adjustment criterion relative to $(X,Y)$ in~$\g[M]_2$.
   Finally, note that $\g[M]_1$ could also be interpreted as a $\DAG$. In that case, it would be amenable relative to $(X,Y)$. This shows that amenability depends crucially on the interpretation of the graph.
   \label{ex2}
\end{example}

\begin{figure}
   \centering
   \begin{subfigure}{.48\textwidth}
     \centering
     \begin{tikzpicture}[>=stealth',shorten >=1pt,auto,node distance=2cm,main node/.style={minimum size=0.6cm,font=\sffamily\Large\bfseries},scale=0.7,transform shape]
   \node[main node]         (X)                        {$X$};
   \node[main node]         (V2) [above left of= X]  	{$V_{2}$};
   \node[main node]         (V1) [above right of= X]  	{$V_{1}$};
   \node[main node]         (V3) [below right of = X]  			{$V_{3}$};
   \node[main node]         (V4) [left of= V3]  		{$V_{4}$};
   \node[main node]       	 (Y)  [left of= V4]         	{$Y$};
   \draw[o->] (V1) edge    (X);
   \draw[o->] (V2) edge    (X);
   \draw[o->] (V3) edge    (X);
   \draw[->] (V3) edge    (V4);
   \draw[o->] (V4) edge    (Y);
   \draw[->] (X) edge   (Y);
   \draw[->] (X) edge   (V4);
   \draw[->] (V3) ..  controls (-0.6,-2.7) ..   (Y);
   \end{tikzpicture}
     \caption{}
   \end{subfigure}
   \unskip
   \vrule
   \hspace{-0.5cm}
   \begin{subfigure}{.48\textwidth}
     \centering
     \begin{tikzpicture}[>=stealth',shorten >=1pt,auto,node distance=2cm,main node/.style={minimum size=0.6cm,font=\sffamily\Large\bfseries},scale=.7,transform shape]
   \node[main node]         (X)                        {$X$};
   \node[main node]         (V2) [above left of= X]  	{$V_{2}$};
   \node[main node]         (V1) [above right of= X]  	{$V_{1}$};
   \node[main node]         (V3) [below right of = X]  {$V_{3}$};
   \node[main node]         (V4) [left of = V3]  		{$V_{4}$};
   \node[main node]       	 (Y)  [left of= V4]         	{$Y$};
   \draw[o->] (V1) edge    (X);
   \draw[o->] (V2) edge    (X);
   \draw[<->] (V3) edge    (X);
   \draw[<->] (V3) edge    (V4);
   \draw[->] (V4) edge    (Y);
   \draw[->] (X) edge   (V4);
   \draw[->] (V3) ..  controls (-0.6,-2.7) ..   (Y);
   \end{tikzpicture}
     \caption{}
   \end{subfigure}
   \caption{\small (a) $\PAG$ $\g[P]_1$, (b) $\PAG$ $\g[P]_{2}$ used in Example~\ref{ex:amenable-noset} and Example~\ref{ex:constr-amenable-noset}.}
   \label{fig:amenable-noset}
\end{figure}
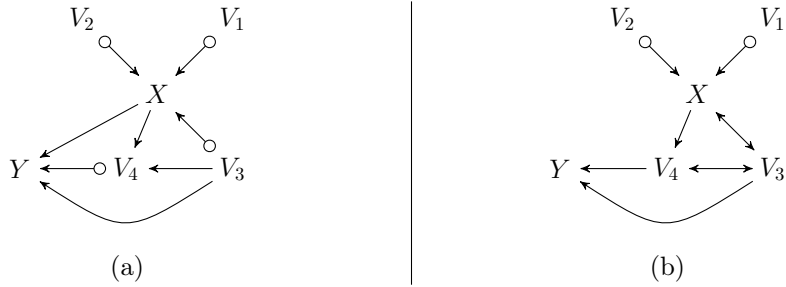

\begin{example}
  Let $\g[P]_1$ and $\g[P]_{2}$ be the $\PAG$s in Figure~\ref{fig:amenable-noset}(a) and Figure~\ref{fig:amenable-noset}(b), respectively.
  Both $\PAG$s are amenable relative to $(X,Y)$. We will show that there is an adjustment set relative to $(X,Y)$ in~$\g[P]_1$ but not in~$\g[P]_{2}$. This illustrates that amenability is not a sufficient
  criterion for the existence of an adjustment set.

  We first consider $\g[P]_1$. Note that $\f{\g[P]_1} = \{ V_4, Y \}$ and that there are two proper definite status non-causal paths from $X$ to $Y$: $X \arrowcirc V_3 \rightarrow Y$ and $X \rightarrow V_4 \leftarrow V_3 \rightarrow Y$. Path $X \arrowcirc V_3 \rightarrow V_4 \circarrow Y$ is not of definite status, as node $V_4$ is not of definite status on this path. Both proper definite status non-causal paths from $X$ to $Y$ are blocked by any set containing $V_3$. Hence, all sets satisfying the generalized adjustment criterion relative to $(X,Y)$ in~$\g[P]_1$ are: $\{ V_3 \}$, $\{V_1, V_3 \}, \{V_2, V_3 \}$ and $\{ V_1, V_2, V_3 \}$.

  In~$\g[P]_{2}$, we have $\f{\g[P]_{2}} = \f{\g[P]_1} = \{ V_4, Y \}$, and there are three proper definite status non-causal paths from $X$ to $Y$ in~$\g[P]_{2}$:  $p_1$ of the form $X \leftrightarrow V_3 \rightarrow Y$, $p_2$ of the form $ X \leftrightarrow V_3 \leftrightarrow V_4 \to Y$ and  $p_3$ of the form $X \rightarrow V_4 \leftrightarrow V_3 \rightarrow Y$. To block $p_1$, we must use $V_3$, and this implies that we must use $V_4$ to block $p_2$. But $V_4 \in \f{\g[P]_{2}}$. Hence, no set satisfies the generalized adjustment criterion $\mathbf{Z}$ relative to $(X,Y)$ in~$\g[P]_{2}$.
  \label{ex:amenable-noset}
\end{example}

\subsection{Proof of Theorem~\ref{theorem:gac}} \label{subsec:proof-theorem-gac}

To prove that the generalized adjustment criterion is sound and complete for adjustment (Theorem~\ref{theorem:gac}), we build on the fact that the adjustment criterion for $\DAG$s and $\MAG$s is sound and complete for adjustment. The adjustment criterion for $\DAG$s was first presented in \cite{shpitser2012validity} and was modified in the unpublished addendum \citep{shpitser2012avalidity}. In Appendix~\ref{subsec:proofs shpitser} we give the revised version of the criterion, as well as new soundness and completeness proofs, relying only on basic probability calculus, linear algebra and the do-calculus rules. 
The adjustment criterion for $\MAG$s was presented and proved to be sound and
complete for adjustment in \citet[][see Theorem~5.8]{vanconstructing}. However, \cite{vanconstructing}'s proof assumed the soundness and completeness of the adjustment criterion for DAGs given in \cite{shpitser2012avalidity}. This latter claim is
proved here in Theorem \ref{theorem:sound and complete ac} in Appendix~\ref{subsec:proofs shpitser}.

Lastly, our proof of Theorem~\ref{theorem:gac} heavily relies on the three lemmas given below. Their proofs can be found in Appendix~\ref{subsec:proofs GAC}.

\begin{lemma}
        Let $\mathbf{X}$ and $\mathbf{Y}$ be disjoint node sets in a $\CPDAG$ $(\PAG)$ $\g$. If $\g$ is amenable (see Definition~\ref{def:gac}) relative to $(\mathbf{X},\mathbf{Y})$, then every $\DAG$ $(\MAG)$ in $[\g]$ is amenable relative to $(\mathbf{X},\mathbf{Y})$. On the other hand, if $\g$ violates \amen{} relative to $(\mathbf{X},\mathbf{Y})$, then there is no adjustment set relative to $(\mathbf{X},\mathbf{Y})$ in~$\g$ (see Definition~\ref{defadjustment}).
   \label{lemma:adj amen gen}
\end{lemma}

\begin{lemma}
  Let $\mathbf{X,Y}$ and $\mathbf{Z}$ be pairwise disjoint node sets in a $\CPDAG$ $(\PAG)$ $\g$. If $\g$ is amenable relative to $(\mathbf{X},\mathbf{Y})$, then the following statements are equivalent:
\begin{enumerate}[label = (\roman*)]
\item\label{l:eqa1} $\mathbf{Z}$ satisfies \forb (see Definition~\ref{def:gac}) relative to $(\mathbf{X},\mathbf{Y})$ in~$\g$.
\item\label{l:eqa2}  $\mathbf{Z}$ satisfies \forb relative to $(\mathbf{X},\mathbf{Y})$ in every $\DAG$ $(\MAG)$ in $[\g]$.
\end{enumerate}
   \label{lemmaprovingeqofacconditiona}
\end{lemma}

\begin{lemma}
     Let $\mathbf{X,Y}$ and $\mathbf{Z}$ be pairwise disjoint node sets in a $\CPDAG$ $(\PAG)$ $\g$. If $\g$ is amenable relative to $(\mathbf{X},\mathbf{Y})$, and $\mathbf{Z}$ satisfies \forb{} relative to $(\mathbf{X},\mathbf{Y})$, then the following statements are equivalent:
   \begin{enumerate}[label = (\roman*)]
\item\label{l:eqb1} $\mathbf{Z}$ satisfies \blck{} (see Definition~\ref{def:gac}) relative to $(\mathbf{X},\mathbf{Y})$ in~$\g$.
\item\label{l:eqb2}  $\mathbf{Z}$ satisfies \blck{} relative to $(\mathbf{X},\mathbf{Y})$ in every $\DAG$ $(\MAG)$ in $[\g]$.
\item\label{l:eqb3} $\mathbf{Z}$ satisfies \blck{} relative to $(\mathbf{X},\mathbf{Y})$ in a $\DAG$ $\g[D]$ $(\MAG \g[M])$ in $[\g]$.
\end{enumerate}
   \label{lemmaequivalenceofcondb}
\end{lemma}

\begin{proofof}[Theorem~\ref{theorem:gac}]
 If $\g$ is a $\DAG$ ($\MAG$), then our criterion reduces to the adjustment criterion from \cite{shpitser2012avalidity} \citep{vanconstructing} which is sound and complete for adjustment (Theorem~\ref{theorem:sound and complete ac} in Appendix~\ref{subsec:proofs shpitser}, \citealp[Theorem~5.8 in][]{vanconstructing}). Hence, we only consider the case that $\g$ is a $\CPDAG$ $(\PAG)$.

   Suppose first that $\mathbf{Z}$ satisfies the generalized adjustment criterion relative to $(\mathbf{X,Y})$ in the $\CPDAG$ $(\PAG)$ $\g$. We need to show that $\mathbf{Z}$ is an adjustment set (see Definition~\ref{defadjustment}) relative to $(\mathbf{X,Y})$ in every $\DAG$ $\g[D]$ ($\MAG$ $\g[M]$) in $[\g]$. By applying Lemmas~\ref{lemma:adj amen gen},~\ref{lemmaprovingeqofacconditiona} and~\ref{lemmaequivalenceofcondb} in turn, it directly follows that $\mathbf{Z}$ satisfies the generalized adjustment criterion relative to $(\mathbf{X,Y})$ in every $\DAG$ $\g[D]$ ($\MAG$ $\g[M]$) in $[\g]$.
  Since the generalized adjustment criterion is sound for adjustment in $\DAG$s ($\MAG$s) 
  (see Theorem~\ref{theorem:soundness ac} in Appendix \ref{subsec:proofs shpitser} and Theorem 5.8 in \citealp{vanconstructing}), $\mathbf{Z}$ is an adjustment set relative to $(\mathbf{X,Y})$ in~every $\g[D]$ ($\g[M]$) in $[\g]$.

   To prove the other direction, suppose that $\mathbf{Z}$ does not satisfy the generalized adjustment criterion relative to $(\mathbf{X,Y})$ in~$\g$. First, suppose that $\g$ violates \amen{} relative to $(\mathbf{X},\mathbf{Y})$. Then by Lemma~\ref{lemma:adj amen gen}, there is no adjustment set relative to $(\mathbf{X},\mathbf{Y})$ in~$\g$.
   Otherwise, $\g$ is amenable relative to $(\mathbf{X,Y})$, but $\mathbf{Z}$ violates \forb or \expandafter\ignorespaces\blck.
   We need to show $\mathbf{Z}$ is not an adjustment set in at least one $\DAG$ $\g[D]$ ($\MAG$ $\g[M]$) in $[\g]$.
   Suppose $\mathbf{Z}$ violates \expandafter\ignorespaces\forb. Then by Lemma~\ref{lemmaprovingeqofacconditiona}, it follows that there exists a $\DAG$ $\g[D]$ ($\MAG$ $\g[M]$) in $[\g]$ such that $\mathbf{Z}$ does not satisfy the generalized adjustment criterion relative to $(\mathbf{X,Y})$ in~$\g[D]$ ($\g[M]$). Since the generalized adjustment criterion is complete for adjustment in $\DAG$s ($\MAG$s)
   (see Theorem~\ref{theorem:completeness ac} in Appendix \ref{subsec:proofs shpitser} and Theorem 5.8 in \citealp{vanconstructing}), it follows that $\mathbf{Z}$ is not an adjustment set relative to $(\mathbf{X,Y})$ in~$\g[D]$ ($\g[M]$).
   Otherwise, suppose $\mathbf{Z}$ satisfies \expandafter\ignorespaces\forb, but violates \expandafter\ignorespaces\blck. Then by Lemma~\ref{lemmaequivalenceofcondb}, it follows that there is a $\DAG$ $\g[D]$ ($\MAG$ $\g[M]$) in $[\g]$ such that $\mathbf{Z}$ does not satisfy the generalized adjustment criterion relative to $(\mathbf{X,Y})$ in~$\g[D]$ ($\g[M]$). Since the generalized adjustment criterion is complete for adjustment in $\DAG$s ($\MAG$s), it follows that $\mathbf{Z}$ is not an adjustment set relative to $(\mathbf{X,Y})$ in~$\g[D]$ ($\g[M]$).
\end{proofof}

\section{Constructing Adjustment Sets} \label{subsec:constructive set gac}

We now present approaches to construct adjustment sets.
First, in Theorem \ref{thm:preprocessing}, we discuss a pre-processing of the node set $\mathbf{X}$ that in conjunction with our generalized adjustment criterion can help identify $f(\mathbf{y}|do(\mathbf{x}))$ in DAG, CPDAG, MAG or PAG $\g$ via adjustment.

As mentioned before, if $f(\mathbf{y}|do(\mathbf{x}))$ is not identifiable via adjustment in $\g$, it may be identifiable through other means. In particular, if $f(\mathbf{y}|do(\mathbf{x}))$ is not identifiable via adjustment in $\g$, there may be a set $\mathbf{X'} \subseteq \mathbf{X}$ such that $f(\mathbf{y}|do(\mathbf{x})) = f(\mathbf{y}| do(\mathbf{x'}))$, and $f(\mathbf{y}| do(\mathbf{x'}))$ is identifiable via adjustment in $\g$. One such example is given in Theorem~\ref{thm:preprocessing}. 

Next, we introduce Theorem~\ref{theorem:general set proof} that will allow us to easily construct adjustment sets that do not contain certain nodes, if any such adjustment set exists. We illustrate the results of Theorem~\ref{theorem:general set proof} with examples in Section~\ref{subsec:examples2} and give the proof of this theorem in Section~\ref{subsec:proof-thm-construct}.  In Section~\ref{sec:implementation} we explain how to leverage previous results of \cite{vanconstructing} to enumerate all (minimal) adjustment sets, and discuss how to implement this procedure efficiently.

\begin{theorem} \label{thm:preprocessing}
Let $\mathbf{X}$ and $\mathbf{Y}$ be disjoint node sets in a causal DAG, CPDAG, MAG or PAG $\g$. Let $\mathbf{X'} \subseteq \mathbf{X}$ such that there is no possibly directed path from $\mathbf{X} \setminus \mathbf{X'}$ to $\mathbf{Y}$ that is proper with respect to $\mathbf{X}$. 
Then 
\begin{align*}
f(\mathbf{y}|do(\mathbf{x})) = \begin{cases}
   f(\mathbf{y}) & \text{if }\mathbf{X'} = \emptyset,\\ f(\mathbf{y}|do(\mathbf{x'})) & \text{otherwise.} \end{cases}
   \end{align*} 
  Furthermore, if $\mathbf{X'} \neq \emptyset$ and if $\mathbf{Z}$ is an adjustment set relative to $(\mathbf{X,Y})$ in $\g$, then $\mathbf{Z}$ is an adjustment set relative to  $(\mathbf{X',Y})$ in $\g$.
\end{theorem}

\noindent{}Following Theorem \ref{thm:preprocessing}, we  recommend pre-processing the set $\mathbf{X}$ as follows: remove all nodes $X \in \mathbf{X}$ that do not have a possibly directed path to $\mathbf{Y}$ which is proper with respect to $\mathbf{X}$ in $\g$. In other words, if $\mathbf{W}$ is the set of all nodes that have a possibly directed path to $\mathbf{Y}$ which is proper with respect to $\mathbf{X}$ in $\g$, then $\mathbf{X'} =  \mathbf{X} \cap \mathbf{W}$. Then by choice of $\mathbf{W}$ and the proof of Theorem \ref{thm:preprocessing}, all nodes in $\mathbf{X'}$ have a possibly directed path to $\mathbf{Y}$ that is proper with respect to $\mathbf{X'}$.

By Theorem~\ref{thm:preprocessing}, this pre-processing of $\mathbf{X}$ cannot hurt in identifying $f(\mathbf{y}|do(\mathbf{x}))$ via adjustment. 
Moreover, there are cases when this pre-processing helps to identify  $f(\mathbf{y}|do(\mathbf{x}))$ via adjustment. For example, in the MAG $\g[M]$ in Figure~\ref{fig:magnoset}, there is no adjustment set relative to $(\{X_1,X_2\}, Y)$. However, there is no possibly directed path from $X_2$ to $Y$ that is proper with respect to $\{X_1,X_2\}$. Furthermore, $\{V_2\}$ is an adjustment set relative to $(X_1,Y)$. Hence, by Theorem~\ref{thm:preprocessing} and Theorem~\ref{theorem:gac}, $f(y|do(x_1,x_2)) = f(y|do(x_1)) = \int_{v_2} f(y|x_1,v_2)f(v_2)dv_2$.

\begin{figure}
     \centering
     \begin{tikzpicture}[>=stealth',shorten >=1pt,auto,node distance=2cm,main node/.style={minimum size=0.6cm,font=\sffamily\Large\bfseries},scale=0.8,transform shape]
   \node[main node]         (X1)                       		{$X_1$};
   \node[main node]       	(V1) [left of= X1]	    		{$V_1$};
   \node[main node]       	(V2) [above right of= X1]	 	 {$V_2$};
   \node[main node]         (Y)  [below right of= V2]  			{$Y$};
   \node[main node]         (X2) [right of= Y]				{$X_{2}$};

   \draw[->] 	(V1)	edge    (X1);
   \draw[->] 	(X1)	edge    (Y);
   \draw[->] 	(V2) 	edge    (X1);
   \draw[->] 	(V2)	edge    (Y);
   \draw[<->] 	(X2) 	edge    (Y);
   \end{tikzpicture}
   \caption{ $\MAG$ $\g[M]$.}
   \label{fig:magnoset}
\end{figure}

We now introduce two definitions that will be used in Theorem \ref{theorem:general set proof}. First, we define the set $\adjustb{\g}$ relative to disjoint node sets $\mathbf{X}$ and $\mathbf{Y}$ in a $\DAG$, $\CPDAG$, $\MAG$ or $\PAG$ $\g$.

\begin{definition}{($\adjustb{\g}$)}
Let $\mathbf{X}$ and $\mathbf{Y}$ be disjoint node sets in a $\DAG$, $\CPDAG$, $\MAG$ or $\PAG$ $\g$. We define
\begin{equation}
\adjustb{\g}=\PossAn(\mathbf{X} \cup \mathbf{Y}, \g)\setminus(\mathbf{X} \cup \mathbf{Y} \cup \fb{\g}).
\end{equation}
\label{def:adjust xy}
\end{definition}
\noindent{}If $\g$ is a $\DAG$ or $\MAG$, then Definition~\ref{def:adjust xy} reduces to the definition of $\adjustb{\g}$ in \citet{vanconstructing}, that is, $\adjustb{\g} = \An(\mathbf{X} \cup \mathbf{Y},\g) \setminus (\mathbf{X} \cup \mathbf{Y} \cup \fb{\g})$.

\begin{definition} \textbf{(Descendral set)}
Let $\mathbf{I}$ be a node set in a $\DAG$, $\CPDAG$, $\MAG$ or $\PAG$ $\g$. Then $\mathbf{I}$ is called descendral in~$\g$ if $\mathbf{I} = \PossDe(\mathbf{I},\g)$.
\label{def:descendral}
\end{definition}

\noindent{}A descendral set is in a sense analogous to an ancestral set, which is a set containing all ancestors of itself. Note that $\fb{\g}$ and $\PossDe(\mathbf{X},\g)$ are both descendral sets. This property will be used throughout the proofs. 
Note that if $\mathbf{A}$ is an ancestral set and $\mathbf{B}$ is a descendral set, then $\mathbf{A} \setminus \mathbf{B}$ is an ancestral set and $\mathbf{B} \setminus \mathbf{A}$ is a descendral set.

\begin{theorem}(\textbf{Constructive set})
Let $\mathbf{X}$ and $\mathbf{Y}$ be disjoint node sets in a $\DAG$, $\CPDAG$, $\MAG$ or $\PAG$ $\g$. Let $\mathbf{I} \supseteq \fb{\g}$ be a descendral set in~$\g$. Then there exists a set $\mathbf{Z}$ that satisfies the generalized adjustment criterion relative to $(\mathbf{X,Y})$ in~$\g$ such that $\mathbf{Z} \cap \mathbf{I} = \emptyset$ if and only if $\adjustb{\g} \setminus \mathbf{I}$ satisfies the generalized adjustment criterion relative to $(\mathbf{X,Y})$ in~$\g$.
\label{theorem:general set proof}
\end{theorem}

\noindent{}The smallest set we can take for $\mathbf{I}$ in Theorem~\ref{theorem:general set proof} is $\fb{\g}$. This leads to Corollary~\ref{cor:adjust xy gac general}. Pearl's back-door criterion does not allow using descendants of $\mathbf{X}$ in a $\DAG$ $\g[D]$. Moreover, the generalized back-door criterion does not allow using possible descendants of $\mathbf{X}$ in a $\DAG$, $\CPDAG$, $\MAG$ or $\PAG$ $\g$. 
Thus, another natural set to consider for $\mathbf{I}$ is $\PossDe(\mathbf{X},\g)$. We will use $\mathbf{I} = \PossDe(\mathbf{X},\g)$ and Theorem~\ref{theorem:general set proof} in Section~\ref{sec:relation} to define sets that satisfy generalized back-door criterion and Pearl's back-door criterion.

\begin{corollary}
Let $\mathbf{X}$ and $\mathbf{Y}$ be disjoint node sets in a $\DAG$, $\CPDAG$, $\MAG$ or $\PAG$ $\g$. The following statements are equivalent:
\begin{enumerate} [label = (\roman*)]
\item\label{cor-adjust-1} There exists an adjustment set relative to $(\mathbf{X,Y})$ in~$\g$.
\item\label{cor-adjust-2} $\adjustb{\g}$ satisfies the generalized adjustment criterion relative to $(\mathbf{X,Y})$ in~$\g$.
\item\label{cor-adjust-3} $\g$ is amenable relative to $(\mathbf{X,Y})$ in~$\g$ and $\adjustb{\g}$ satisfies \blck{} relative to $(\mathbf{X,Y})$ in~$\g$.
\end{enumerate}
\label{cor:adjust xy gac general}
\end{corollary}

\subsection{Examples} \label{subsec:examples2}
We now provide some examples that illustrate the construction of adjustment sets.
\begin{example}
Consider again the $\CPDAG$ $\g[C]$ in Figure~\ref{fig:cpdagexample}(a). As previously discussed in Example~\ref{ex:everything-works}, $\g[C]$ is amenable relative to $(X,Y)$.
The set $\adjust{\g[C]} = \{ X, Y, I, A, Z, B \} \setminus \{X,Y\} = \{I,A,Z,B \}$ satisfies \blck{} relative to $(X,Y)$ in~$\g[C]$. Hence, by Corollary~\ref{cor:adjust xy gac general}, $\{I,A,Z,B \}$ satisfies the generalized adjustment criterion relative to $(X,Y)$ in~$\g[C]$.
\label{ex:constr-everything-works}
\end{example}
\begin{example}
Consider again the $\PAG$s $\g[P]_1$ and $\g[P]_{2}$ in Figure~\ref{fig:amenable-noset}(a) and Figure~\ref{fig:amenable-noset}(b), respectively. As previously discussed in Example~\ref{ex:amenable-noset}, both $\g[P]_{1}$ and $\g[P]_{2}$ are amenable relative to $(X,Y)$.

In~$\g[P]_{1}$, $\f{\g[P]_{1}} = \{ V_4, Y \}$, so $\adjust{\g[P]_{1}} =$ $\{ X, Y, V_1, V_2, V_3, V_4\}$ $\setminus \{X,Y,V_4\} $ $ = \{V_1,V_2,V_3\}$. Since $\{V_1,V_2,V_3\}$ satisfies \blck{} relative to $(X,Y)$ in~$\g$, it follows that $\{V_1,V_2,V_3\}$ satisfies the generalized adjustment criterion relative to $(X,Y)$ in~$\g$.

In~$\g[P]_{2}$, again $\f{\g[P]_{2}} =\{ V_4, Y \}$, so $\adjust{\g[P]_{2}} =$ $\{ X, Y, V_1, V_2, V_3, V_4\}$ $\setminus \{X,Y,V_4\} $ $ = \{V_1,V_2,V_3\}$. Since $\{V_1,V_2,V_3\}$ does not block the path $X \leftrightarrow V_3 \leftrightarrow V_4 \rightarrow Y$ it does not satisfy \blck{} relative to $(X,Y)$ in~$\g$. Hence, Corollary~\ref{cor:adjust xy gac general} implies that there is no adjustment set relative to $(X,Y)$ in~$\g[P]_{2}$.
\label{ex:constr-amenable-noset}
\end{example}

\subsection{Proof of Theorem~\ref{theorem:general set proof}} \label{subsec:proof-thm-construct}

To prove Theorem~\ref{theorem:general set proof} we heavily rely on Lemma~\ref{lemma:general path proof} and Lemma~\ref{lemma:rich prime} given below. Their proofs are given in Appendix~\ref{subsec:proofs constructive}. Lemma~\ref{lemma:rich prime} is related to Lemma~1 from \citet{richardson2003markov} (see Lemma~\ref{lemma:rich} in Appendix~\ref{subsec:additional}).

\begin{lemma}
Let $\mathbf{X}$ and $\mathbf{Y}$ be disjoint node sets in a $\DAG$, $\CPDAG$, $\MAG$ or $\PAG$ $\g$. Let $\mathbf{I} \supseteq \fb{\g}$ be a descendral set in~$\g$ (see Definition~\ref{def:descendral}). If there is a proper definite status non-causal path from $\mathbf{X}$ to $\mathbf{Y}$ in~$\g$ that is m-connecting given $\adjustb{\g} \setminus \mathbf{I}$, then there is a path $p$ from $\textbf{X}$ to $\textbf{Y}$ in $\g$ such that:
\begin{enumerate}[label=(\roman*)]
\item\label{l:gpp0} $p$ is a proper definite status non-causal path from $\mathbf{X}$ to $\mathbf{Y}$ in~$\g$, and
\item\label{l:gpp1} all colliders on $p$ are in $\adjustb{\g} \setminus \mathbf{I}$, and
\item\label{l:gpp2} all definite non-colliders on $p$ are in $\mathbf{I}$, and
\item\label{l:gpp3} for any collider $C$ on $p$, there is an unshielded possibly directed path from $C$ to $\mathbf{X} \cup \mathbf{Y}$, that starts with $\circarrow$ or $\rightarrow$.
\end{enumerate}
\label{lemma:general path proof}
\end{lemma}
\begin{lemma}
Let $\mathbf{X,Y}$ and $\mathbf{Z}$ be pairwise disjoint node sets in a $\DAG$, $\CPDAG$, $\MAG$ or $\PAG$ $\g$. Let $\mathbf{I} \supseteq \fb{\g}$ be a node set in~$\g$ such that $\mathbf{Z} \cap \mathbf{I} = \emptyset$. Let $p$ be a path from $\mathbf{X}$ to $\mathbf{Y}$ in~$\g$ such that:
\begin{enumerate}[label=(\roman*)]
\item\label{l:rp0} $p$ is a proper definite status non-causal path from $\mathbf{X}$ to $\mathbf{Y}$ in~$\g$, and
\item\label{l:rp2} all colliders on $p$ are in $\An(\mathbf{X} \cup \mathbf{Y} \cup \mathbf{Z}, \g) \setminus \mathbf{I}$, and
\item\label{l:rp1} no definite non-collider on $p$ is in $\mathbf{Z}$.
\end{enumerate}
Then there is a proper definite status non-causal path from $\mathbf{X}$ to $\mathbf{Y}$ that is m-connecting given $\mathbf{Z}$ in~$\g$.
\label{lemma:rich prime}
\end{lemma}
\begin{proofof}[Theorem~\ref{theorem:general set proof}]
We only prove the non-trivial direction. Thus, assume there is a set $\mathbf{Z}$ satisfying the generalized adjustment criterion relative to $(\mathbf{X,Y})$ in~$\g$ such that $\mathbf{Z} \cap \mathbf{I} = \emptyset$. We will prove that $\adjustb{\g} \setminus \mathbf{I}$ satisfies the generalized adjustment criterion relative to $(\mathbf{X,Y})$ in~$\g$.

Since $\mathbf{Z}$ satisfies the generalized adjustment criterion relative to $(\mathbf{X,Y})$ in~$\g$, $\g$ is amenable relative to $(\mathbf{X,Y})$. Additionally, since $\fb{\g} \subseteq \mathbf{I}$, $\adjustb{\g} \setminus \mathbf{I}$ satisfies \forb{} relative to $(\mathbf{X,Y})$ in~$\g$. It is only left to prove that $\adjustb{\g} \setminus \mathbf{I}$ satisfies \blck{} relative to $(\mathbf{X,Y})$ in~$\g$.

Suppose for a contradiction that there is a proper definite status non-causal path from $\mathbf{X}$ to $\mathbf{Y}$ that is m-connecting given $\adjustb{\g} \setminus \mathbf{I}$.
Then we can choose a path $\pstar$ in $\g$ that satisfies~\ref{l:gpp0}$-$\ref{l:gpp3} in Lemma~\ref{lemma:general path proof}. Then $\pstar$ also satisfies \ref{l:rp0}~in~Lemma~\ref{lemma:rich prime}.
By \ref{l:gpp2}~in~Lemma~\ref{lemma:general path proof}, every definite non-collider on $\pstar$ is in $\mathbf{I}$. Since $\mathbf{Z} \cap \mathbf{I} = \emptyset$, no definite non-collider on $\pstar$ is in $\mathbf{Z}$. So $\pstar$ satisfies \ref{l:rp1}~in~Lemma~\ref{lemma:rich prime}. Also, since by \ref{l:gpp3}~in~Lemma~\ref{lemma:general path proof} there is a possibly directed unshielded path $\pstar[q]$ from every collider $C$ on $\pstar$ to $\mathbf{X} \cup \mathbf{Y}$ that starts with $C \circarrow$, Lemma~\ref{lemma:unshielded edges} implies that any other edge on $\pstar[q]$ (if there is any) is directed in $\g$.

Then if $\g$ is a $\DAG$, $\CPDAG$ or $\MAG$, it follows from \ref{l:gpp3}~in~Lemma~\ref{lemma:general path proof} that all colliders on $\pstar$ are in $\An(\mathbf{X} \cup \mathbf{Y},\g)$. Combining this with \ref{l:gpp1}~in~Lemma~\ref{lemma:general path proof} implies that all colliders on $\pstar$ are in $\An(\mathbf{X} \cup \mathbf{Y},\g) \setminus \mathbf{I}$, so that $\pstar$ satisfies \ref{l:rp2}~in~Lemma~\ref{lemma:rich prime}. Hence, if $\g$ is a $\DAG$, $\CPDAG$ or $\MAG$, all conditions of Lemma~\ref{lemma:rich prime} are satisfied, which implies that $\mathbf{Z}$ does not satisfy \blck{} relative to $(\mathbf{X,Y})$ in $\g$.

Thus, assume $\g$ is a $\PAG$.
 Let $\g[M] \in [\g]$ be a $\MAG$ obtained from $\g$ by first replacing all partially directed edges $\circarrow$ by directed edges $\rightarrow$, and then orienting all non-directed edges $\circcirc$ as a $\DAG$ without unshielded colliders (see Lemma~\ref{lemmamarlrx} in Appendix~\ref{subsec:additional}).
  Let $p$ be the path in~$\g[M]$ corresponding to  $\pstar$ in~$\g$.
 Then $p$ satisfies \ref{l:rp0}~and~\ref{l:rp1}~in~Lemma~\ref{lemma:rich prime}. By the choice of $\g[M]$ and Lemma~\ref{lemma:unshielded edges}, any possibly directed unshielded path $\pstar[q]$ in $\g$ that starts with a partially directed edge $\circarrow$, corresponds to a directed path $q$ in $\g[M]$. Hence, by \ref{l:gpp3}~in~Lemma~\ref{lemma:general path proof}, every collider on $p$ is in $\An(\mathbf{X} \cup \mathbf{Y},\g[M])$. Since $\pstar$ satisfies \ref{l:gpp1}~in~Lemma~\ref{lemma:general path proof} in $\g$, no collider on $\pstar$ is in $\mathbf{I}$. Hence, also no collider on $p$ is in $\mathbf{I}$. Then all colliders on $p$ are in $\An(\mathbf{X} \cup \mathbf{Y},\g[M]) \setminus \mathbf{I}$.
Additionally, since $\mathbf{I} \supseteq \fb{\g}$, and $\fb{\g} \supseteq \fb{\g[M]}$, it follows that $p$ satisfies \ref{l:rp2}~in~Lemma~\ref{lemma:rich prime}. Thus, all conditions of Lemma~\ref{lemma:rich prime} are satisfied, which implies that $\mathbf{Z}$ does not satisfy \blck{} relative to $(\mathbf{X,Y})$ in $\g[M]$. This contradicts Lemma~\ref{lemmaequivalenceofcondb}.
\end{proofof}

\subsection{Implementation} \label{sec:implementation}

We now discuss how one can implement the generalized adjustment criterion in an algorithmically efficient manner, and describe our implementation in the software \texttt{dagitty} and the R package \texttt{pcalg}.

\textbf{Verification of the criterion.} Given a $\DAG$, $\CPDAG$, $\MAG$ or $\PAG$ $\g$ and three disjoint node sets $\bX,\bY$, and $\bZ$, we wish to test whether $\bZ$ fulfills the generalized adjustment criterion with respect to $(\bX,\bY)$. Of course, we could do this simply by verifying the three conditions of the generalized adjustment criterion (see Definition~\ref{def:gac}). However, \blck is a statement about individual paths, which can pose problems for large graphs. 

As a worst-case example, consider a $\DAG$ with $X$, $Y$ and $p$ remaining variables. Let every pair of variables be connected by an edge. Then the $\DAG$ contains $\sum_{i=0}^p p! / i! \approx p!\,e$ paths from $X$ to $Y$. Thus, a direct implementation of the generalized adjustment criterion has an exponential runtime in $p$. Still, a verbatim implementation of the criterion can be useful for verification and  didactic purposes, as well as for sparse graphs with few paths, and we therefore provide one in the function \texttt{gac} of the R package \texttt{pcalg}.

The key result for implementing the criterion in an efficient manner is Theorem~\ref{theorem:gac-alt}, which replaces the path blocking condition by an $m$-separation condition in a subgraph of $\g$, the proper back-door graph $\gpbd{\mathbf{XY}}$. This condition can be checked efficiently by a simple depth-first or breadth-first graph traversal, known as the ``Bayes-Ball algorithm'' \citep{Shachter1998}. 
Specifically, for graphs represented as adjacency lists, the runtime is $O(|p|+|\mathbf{E}|)$ where $|\mathbf{E}|$ is the number of edges. Our implementation of this method can be accessed via the function \texttt{isAdjustmentSet} of the R package \texttt{dagitty}.

\textbf{Constructing adjustment sets.} Given a $\DAG$, $\CPDAG$, $\MAG$ or $\PAG$ $\g$ and two disjoint variable sets $\bX,\bY$, we wish to find one or several sets $\bZ$ that fulfill the generalized adjustment criterion relative to $(\bX,\bY)$. If a single set is sufficient, we can directly apply the main result of Section~\ref{subsec:constructive set gac} and construct $\adjustb{\g}$ (see Definition~\ref{def:adjust xy}) and verify whether it satisfies the generalized adjustment criterion relative to $(\bX,\bY)$. Since this set is defined in terms of (possible) ancestors of $\bX$ and $\bY$, it can be constructed by graph traversal in linear time. However, $\adjustb{\g}$ can be a large set, since it contains all (possible) ancestors of $\bX$ and $\bY$ except the forbidden nodes. This may result in a loss of statistical precision.

To avoid this, it is of interest to construct all possible adjustment sets. 
Again, Theorem~\ref{theorem:gac-alt} is key to achieving this: it allows us to use the algorithmic framework developed by \citet{vanconstructing} for constructing and enumerating $m$-separating sets in $\DAG$s and $\MAG$s. For $\DAG$s and $\MAG$s, this can be directly applied. We propose the following procedure for $\CPDAG$s and $\PAG$s:

For a given $\CPDAG$ or $\PAG$ $\g$ and disjoint node sets $\mathbf{X}$ and $\mathbf{Y}$,

\begin{enumerate}[label=(\arabic*)]

\item Check if $\g$ is amenable relative to $(\mathbf{X},\mathbf{Y})$. If not, stop
because there is no adjustment set relative to $(\mathbf{X},\mathbf{Y})$ in~$\g$. Otherwise,
continue.

\item Find $\fb{\g}$.

\item If $\g$ is a $\CPDAG$, orient $\g$ into a $\DAG$ $\g[D]$ in $[\g]$. If $\g$ is a $\PAG$, orient $\g$ into a $\MAG$ $\g[M]$ in $[\g]$ according to Theorem~2 from \cite{zhang2008completeness} (see Lemma~\ref{lemmamarlrx} in Appendix~\ref{subsec:additional}).

\item By Lemma~\ref{lemmaequivalenceofcondb}, finding all sets satisfying the generalized adjustment criterion relative to $(\mathbf{X},\mathbf{Y})$ in~$\g$ is equivalent to finding all sets $\mathbf{Z}$ satisfying \condtwoprime relative to $(\mathbf{X},\mathbf{Y})$ in~$\g[D]$ ($\g[M]$) such that $\mathbf{Z} \cap (\mathbf{X} \cup \mathbf{Y} \cup \fb{\g}) = \emptyset$. Thus, we can apply the algorithms from \cite{vanconstructing} on $\g[D]$ ($\g[M]$). These algorithms are able to deal with the additional restriction that the resulting set must not contain nodes in $\fb{\g}$.

\end{enumerate}

Thus, through this simple procedure we gain complete access to all functions in the algorithmic framework by \cite{vanconstructing}. These include listing all adjustment sets and all minimal adjustment sets (in polynomial time per set that is listed). We have implemented these features in the function \texttt{adjustmentSets} of the R package \texttt{dagitty}.

\section{Relationship to (Generalized) Back-door Criteria} \label{sec:relation}

We now discuss the relationship between our generalized adjustment criterion and some other existing graphical criteria for covariate adjustment.
In particular, we discuss  Pearl's back-door criterion (see Definition~\ref{def:bc}) and the generalized back-door criterion (see Definition~\ref{def:gbc}) and give constructive sets for both in Section~\ref{subsec:GBC}.
We use the results from Section~\ref{subsec:GBC} to precisely characterize the differences between our generalized adjustment criterion, Pearl's back-door criterion and the generalized back-door criterion in Theorem~\ref{theorem:unified-nope} of Section~\ref{subsubsec:equivalence gac gbc}.
We illustrate the results of Sections~\ref{subsec:GBC} and~\ref{subsubsec:equivalence gac gbc} with examples in Section~\ref{subsec:ex3}.
\begin{definition} (\textbf{Back-door criterion}; \citealp{pearl1993bayesian})
 Let $X$ and $Y$ be distinct nodes in a $\DAG$ $\g[D]$. A set of nodes $\mathbf{Z}$ not containing $X$ or $Y$ satisfies the back-door criterion relative to $(X,Y)$ in~$\g[D]$ if:
\begin{enumerate}[label=(\roman*)]
\item\label{d:bc1} no node in $\mathbf{Z}$ is a descendant of $X$, and
\item\label{d:bc2} $\mathbf{Z}$ blocks every path between $X$ and $Y$ that contains an arrow into $X$.
\end{enumerate}
If $\mathbf{X}$ and $\mathbf{Y}$ are two disjoint sets of nodes in~$\g[D]$, then $\mathbf{Z}$ is said to satisfy the back-door criterion relative to $(\mathbf{X,Y})$ if it satisfies the criterion relative to any pair $(X,Y)$ such that $X \in \mathbf{X}$, and $Y \in \mathbf{Y}$.
A set $\mathbf{Z}$ that satisfies the back-door criterion relative to $(\mathbf{X,Y})$ in~$\g[D]$ is called a back-door set relative to $(\mathbf{X,Y})$ in~$\g[D]$.
\label{def:bc}
\end{definition}
\begin{definition}{(\textbf{Back-door path}; \citealp{maathuis2013generalized})}
Let $X$ and $Y$ be distinct nodes in a $\DAG$, $\CPDAG$, $\MAG$ or $\PAG$ $\g$. A path from $X$ to $Y$ in~$\g$ is a back-door path if it does not start with a visible edge out of $X$.
\label{def:backdoor path}
\end{definition}
\begin{definition}{(\textbf{Generalized back-door criterion}; \citealp{maathuis2013generalized})} Let $\mathbf{X,Y}$ and $\mathbf{Z}$ be pairwise disjoint node sets in a $\DAG$, $\CPDAG$, $\MAG$ or $\PAG$ $\g$. Then $\mathbf{Z}$ satisfies the generalized back-door criterion relative to $(\mathbf{X,Y})$ in~$\g$ if:
\begin{enumerate}[label=(\roman*)]
\item\label{gbc:cond1} $\mathbf{Z}$ does not contain possible descendants of $\mathbf{X}$ in~$\g$, and
\item\label{gbc:cond2} for every $X \in \mathbf{X}$, the set $\mathbf{Z} \cup \mathbf{X} \setminus \{X\}$ blocks every definite status back-door path from $X$ to any member of $\mathbf{Y}$, if any, in~$\g$.
\end{enumerate}
A set $\mathbf{Z}$ that satisfies the generalized back-door criterion relative to $(\mathbf{X,Y})$ in~$\g$ is called a generalized back-door set relative to $(\mathbf{X,Y})$ in~$\g$.
\label{def:gbc}
\end{definition}

\subsection{Constructing (Generalized) Back-door Sets} \label{subsec:GBC}

We first focus on Pearl's back-door criterion. If $|\mathbf{X}| = 1$, the existence and construction of a back-door set in a causal $\DAG$ $\g[D]$ is well understood.
If $Y \in \Pa(X,\g[D])$, then there is no back-door set relative to $(X,Y)$ in~$\g[D]$, but it is obvious that $f(y \given \ddo(x)) = f(y)$. If $Y \notin \Pa(X,\g[D])$, then $\Pa(X,\g[D])$ is a back-door set relative  to $(X,Y)$ in~$\g[D]$.

If $|\mathbf{X}| \geq 1$, the construction of a back-door set relative to $(\mathbf{X,Y})$ in~$\g[D]$ is less obvious. One could perhaps think that any set $\mathbf{Z}$ that satisfies our generalized adjustment criterion relative to $(\mathbf{X,Y})$ in~$\g[D]$ such that $\mathbf{Z} \cap \De(\mathbf{X},\g[D])=\emptyset$ satisfies Pearl's back-door criterion. This is not true, as shown in Lemma~\ref{lemma:gbc bc} that describes the graphical pattern that appears when there is an adjustment set but no back-door set. An example of a $\DAG$ that satisfies \ref{l:gbc-bc1}-\ref{l:gbc-bc4} in Lemma~\ref{lemma:gbc bc} is given in Figure~\ref{fig:constr-gbc-nodsep}. 
Using this result and Theorem~\ref{theorem:general set proof} we are able to define a specific set that satisfies Pearl's back-door criterion, when such a set exists. This result is given in Corollary~\ref{cor:constr-bc}.

\begin{lemma}
Let $\mathbf{X}$ and $\mathbf{Y}$ be disjoint node sets in a $\DAG$ $\g[D]$. Assume there is a set satisfying the generalized adjustment criterion $\mathbf{Z}$ relative to $(\mathbf{X,Y})$ in~$\g[D]$ such that $\mathbf{Z} \cap \De(\mathbf{X},\g[D])=\emptyset$. Then there is no back-door set relative to $(\mathbf{X,Y})$ in~$\g[D]$ if and only if there is a path $p$ from $\mathbf{X}$ to $\mathbf{Y}$ such that:
\begin{enumerate}[label=(\roman*)]
\item\label{l:gbc-bc1} $p$ is a back-door path, and
\item\label{l:gbc-bc2} a proper subpath of $p$ is a causal path, and
\item\label{l:gbc-bc3} there are no colliders on $p$, and
\item\label{l:gbc-bc4} all nodes on $p$ are in $\De(\mathbf{X},\g[D])$.
\end{enumerate}
\label{lemma:gbc bc}
\end{lemma}
\begin{corollary}(\textbf{Constructive back-door set})
Let $\mathbf{X}$ and $\mathbf{Y}$ be disjoint node sets in a $\DAG$ $\g[D]$.
The following statements are equivalent:
\begin{enumerate} [label = (\roman*)]
\item\label{constr-bc-1} There exists a set that satisfies Pearl's back-door criterion relative to $(\mathbf{X,Y})$ in~$\g[D]$.
\item\label{constr-bc-2} $\adjustb{\g[D]}\setminus \De(\mathbf{X},\g[D])$ satisfies Pearl's back-door criterion relative to $(\mathbf{X,Y})$ in~$\g[D]$.
\item\label{constr-bc-3} 
For all $X \in \mathbf{X}$, $Y \in \mathbf{Y}$, $\adjustb{\g[D]} \setminus \De(\mathbf{X},\g[D])$ satisfies condition~\ref{d:bc2} of Pearl's back-door criterion.
\end{enumerate}
\label{cor:constr-bc}
\end{corollary}

\citet{maathuis2013generalized} presented a constructive generalized back-door set for a $\DAG$, $\CPDAG$, $\MAG$ or $\PAG$ $\g$ when $|\mathbf{X}| = 1$. In Lemma~\ref{lemma:adjust xy gbc cond 2}, we show that any set $\mathbf{Z}$ that satisfies our generalized adjustment criterion relative to $(\mathbf{X,Y})$ in~$\g$ such that $\mathbf{Z} \cap \PossDe(\mathbf{X},\g) = \emptyset$ is a generalized back-door set relative to $(\mathbf{X,Y})$ in~$\g$. Using this result and  Theorem~\ref{theorem:general set proof}, we give a constructive set that satisfies the generalized back-door criterion, when such a set exists. This set is given in Corollary~\ref{cor:constr-set-gbc}.
If $|\mathbf{X}| = 1$, our constructive set for the generalized back-door criterion is a superset of the set presented in \citet{maathuis2013generalized} (Corollary~\ref{cor:dsep adjust} in Appendix~\ref{subsec:proofs relation}).

\begin{lemma}
Let $\mathbf{X,Y}$ and $\mathbf{Z}$ be pairwise disjoint node sets in a $\DAG$, $\CPDAG$, $\MAG$ or $\PAG$ $\g$. If $\g$ is amenable relative to $(\mathbf{X,Y})$ and $\mathbf{Z}$ satisfies \blck{} relative to $(\mathbf{X,Y})$ in~$\g$, then $\mathbf{Z}$ satisfies condition~\ref{gbc:cond2} of the generalized back-door criterion relative to $(\mathbf{X,Y})$ in~$\g$.
\label{lemma:adjust xy gbc cond 2}
\end{lemma}

\begin{corollary}(\textbf{Constructive generalized back-door set})
Let $\mathbf{X}$ and $\mathbf{Y}$ be disjoint node sets in a $\DAG$, $\CPDAG$, $\MAG$ or $\PAG$ $\g$.
The following statements are equivalent:
\begin{enumerate} [label = (\roman*)]
\item\label{constr-gbc-1} There exists a set that satisfies the generalized back-door criterion relative to $(\mathbf{X,Y})$ in~$\g$.
\item\label{constr-gbc-2} $\adjustb{\g}\setminus \PossDe(\mathbf{X},\g)$ satisfies the generalized back-door criterion relative to $(\mathbf{X,Y})$ in~$\g$.
\item\label{constr-gbc-3} $\g$ is amenable and $\adjustb{\g}\setminus \PossDe(\mathbf{X},\g)$ satisfies condition~\ref{gbc:cond2} of the  generalized back-door criterion relative to $(\mathbf{X,Y})$ in~$\g$.
\end{enumerate}
\label{cor:constr-set-gbc}
\end{corollary}

\subsection{Graphs for Which the Criteria Differ} \label{subsubsec:equivalence gac gbc}

We now define graphical conditions for the existence of a set satisfying one, two, or all three of the mentioned criteria. The main result of this section is given in Theorem~\ref{theorem:unified-nope}, which describes the $4$ graphical patterns that can appear when there is no set satisfying at least one of the criteria.

Previously, in Lemma~\ref{lemma:adj amen gen} and Lemma~\ref{lemma:general path proof}, we described such patterns for our generalized adjustment criterion.
In Section~\ref{subsec:GBC}, we showed that any set $\mathbf{Z}$ that satisfies our generalized adjustment criterion such that $\mathbf{Z} \cap \PossDe(\mathbf{X},\g) = \emptyset$ satisfies the generalized back-door criterion.
Thus, to describe an additional pattern that appears when there is no set that satisfies the generalized back-door criterion relative to $(\mathbf{X,Y})$ in~$\g$ we give Lemma~\ref{lemma:no gbc path}.
Lastly, to complete Theorem~\ref{theorem:unified-nope} we add the pattern described in Lemma~\ref{lemma:gbc bc} that additionally appears when there is no set that satisfies Pearl's back-door criterion relative to $(\mathbf{X,Y})$ in~$\g$.  
Thus, Theorem~\ref{theorem:unified-nope} summarizes and subsumes the results of Lemmas~\ref{lemma:adj amen gen}, \ref{lemma:general path proof}, \ref{lemma:gbc bc} and \ref{lemma:no gbc path}.

\begin{lemma} 
Let $\mathbf{X}$ and $\mathbf{Y}$ be disjoint node sets in a $\DAG$, $\CPDAG$, $\MAG$ or $\PAG$ $\g$ such that there exists an adjustment set relative to $(\mathbf{X,Y})$ in~$\g$. There is no set that satisfies the generalized back-door criterion relative to $(\mathbf{X,Y})$ in~$\g$ if and only if there is a path $p$ from $X \in \mathbf{X}$ to $Y \in \mathbf{Y}$ in~$\g$ and a node $V$ on $p$ such that:
\begin{enumerate}[label=(\roman*)]
\item\label{l:nogbc0} $p$ is a proper definite status non-causal path from $\mathbf{X}$ to $\mathbf{Y}$ in~$\g$, and
\item\label{l:nogbc1} $V \in \adjustb{\g} \cap \PossDe(\mathbf{X},\g)$ and is a definite non-collider on $p$, and
\item\label{l:nogbc2} any collider on $p$ is in $\adjustb{\g} \setminus \PossDe(\mathbf{X},\g)$ and any definite non-collider on $p(V,Y)$ is in $\fb{\g}$, and
\item\label{l:nogbc3} path $p$ is of the form $X \leftarrow \dots \leftarrow V \arrowbullet W \dots Y$, where $W=Y$ is possible and if $W \neq Y$ then $V \leftrightarrow W$ is on $p$.
\end{enumerate}
\label{lemma:no gbc path}
\end{lemma}
\begin{theorem}
Let $\mathbf{X}$ and $\mathbf{Y}$ be disjoint node sets in a $\DAG$, $\CPDAG$, $\MAG$ or $\PAG$ $\g$.
Consider the following criteria:
\begin{enumerate}[label=(\arabic*)]
\item\label{nogac1} $\g$ violates \amen{} relative to $(\mathbf{X,Y})$.
\item\label{nogac2} There is a proper definite status non-causal path $p$ from $\mathbf{X}$ to $\mathbf{Y}$ in~$\g$ such that every collider on $p$ is in $\adjustb{\g}$ and every definite non-collider on $p$ is in $\fb{\g}$.
\item\label{nogbc3} There is a path $p$ from $\mathbf{X}$ to $\mathbf{Y}$ in~$\g$ that satisfies~\ref{l:nogbc0}$-$\ref{l:nogbc3} in Lemma~\ref{lemma:no gbc path}.
\item\label{nobc3} There is a back-door path $p$ from $\mathbf{X}$ to $\mathbf{Y}$ in~$\g$ that satisfies~\ref{l:gbc-bc1}$-$\ref{l:gbc-bc4} in Lemma~\ref{lemma:gbc bc}.
\end{enumerate}
The following hold:
\begin{enumerate}[label=(\roman*)]
\item\label{nogac} There is no set that satisfies the generalized adjustment criterion relative to $(\mathbf{X,Y})$ in~$\g$ if and only if~\ref{nogac1}~or~\ref{nogac2} are satisfied.
\item\label{nogbc} There is no generalized back-door set relative to $(\mathbf{X,Y})$ in~$\g$ if and only if~\ref{nogac1}, \ref{nogac2} or \ref{nogbc3} are satisfied.
\item\label{nobc} If $\g$ is a $\DAG$, then there is no back-door set relative to $(\mathbf{X,Y})$ in~$\g$ if and only if~\ref{nogac1},~\ref{nogac2},~\ref{nogbc3}~or~\ref{nobc3} are satisfied.
\end{enumerate}
\label{theorem:unified-nope}
\end{theorem}

We now further explore condition~\ref{nogac2} in Theorem~\ref{theorem:unified-nope} under the assumption that the DAG, CPDAG, MAG or PAG $\g$ is amenable relative to disjoint node sets $(\mathbf{X,Y})$ (that is, \ref{nogac1} in Theorem~\ref{theorem:unified-nope} is violated). Condition~\ref{nogac2} in Theorem~\ref{theorem:unified-nope} is satisfied relative to $\mathbf{(X,Y)}$ in $\g$ if and only if there is no adjustment set relative to $(\mathbf{X,Y})$ in $\g$ (Theorem~\ref{theorem:gac}). Corollary~\ref{cor:noforbx} provides a simple sufficient condition for condition \ref{nogac2} in Theorem~\ref{theorem:unified-nope} to be satisfied in DAGs, CPDAGs, MAGs and PAGs, as well as a necessary and sufficient condition for condition \ref{nogac2} in Theorem~\ref{theorem:unified-nope} in certain DAGs and CPDAGs. 

\begin{corollary}
Let $\mathbf{X}$ and $\mathbf{Y}$ be disjoint node sets in a DAG, CPDAG, MAG or PAG $\g$ such that $\g$ is amenable relative to $(\mathbf{X,Y})$. The following statements hold:
\begin{enumerate}[label=(\roman*)]
\item\label{noforbx1}  If $\mathbf{X} \cap \fb{\g} \neq \emptyset$, then there is no adjustment set relative to $(\mathbf{X,Y})$ in $\g$.
\item\label{noforbx2} Let $\g$ be a DAG or CPDAG and $\mathbf{Y} \subseteq \PossDe(\mathbf{X},\g)$. Then $\mathbf{X} \cap \fb{\g} \neq \emptyset$ if and only if there is no adjustment set relative to $(\mathbf{X,Y})$ in $\g$.
\end{enumerate}
\label{cor:noforbx}
\end{corollary}

\noindent{}A necessary condition for both \ref{nogbc3}~and~\ref{nobc3}~in~Theorem~\ref{theorem:unified-nope} is that $\g$ contains a (possibly) directed path from one node in $\mathbf{X}$ to another node in $\mathbf{X}$. Thus, if $|\mathbf{X}|=1$, both \ref{nogbc3}~and~\ref{nobc3}~in~Theorem~\ref{theorem:unified-nope} are violated. Hence, if $|\mathbf{X}|=1$ there is a (generalized) back-door set relative to $(\mathbf{X,Y})$ in a $\DAG$ (or a $\CPDAG$, $\MAG$ or $\PAG$) $\g$, if and only if there is a set satisfying the generalized adjustment criterion relative to $(\mathbf{X,Y})$ in~$\g$.

We finish this section by giving two corollaries that describe some additional simple conditions under which there exists a set satisfying two or all three of the discussed adjustment criteria. The intuition behind these results is as follows. A necessary condition for the existence of a path $p$ from $\mathbf{X}$ to $\mathbf{Y}$ in a $\DAG$ $\g[D]$ that satisfies Lemma~\ref{lemma:gbc bc} is the existence of a causal path from one node in $\mathbf{X}$ to another node in $\mathbf{X}$ in $\g[D]$. This gives us Corollary~\ref{cor:equivalence gac bc}. Similarly, a necessary condition for the existence of a path $p$ from $\mathbf{X}$ to $\mathbf{Y}$ in~a $\DAG$, $\CPDAG$, $\MAG$ or $\PAG$ $\g$ that satisfies Lemma~\ref{lemma:no gbc path} is the existence of a causal path from one node in $\mathbf{X}$ to another node in $\mathbf{X}$ that contains at least one node not in $\mathbf{X}$ in $\g$. Or, in the case when $\g$ is a $\DAG$ or $\CPDAG$, another necessary condition for the existence of a path $p$ from $\mathbf{X}$ to $\mathbf{Y}$ that satisfies Lemma~\ref{lemma:no gbc path} in $\g$ is that $\mathbf{Y} \not\subseteq \PossDe(\mathbf{X},\g)$. This gives us Corollary~\ref{cor:equivalence gac gbc}.

\begin{corollary}
Let $\mathbf{X}$ and $\mathbf{Y}$ be disjoint node sets in a $\DAG$ $\g[D]$. If there is no directed path from one node in $\mathbf{X}$ to another node in $\mathbf{X}$ in~$\g[D]$, then the following statements are equivalent:
\begin{enumerate}[label=(\roman*)]
\item\label{l2:1} There exists a set that satisfies the generalized adjustment criterion relative to $(\mathbf{X,Y})$ in~$\g[D]$.
\item\label{l2:2} There exists a back-door set relative to $(\mathbf{X,Y})$ in~$\g[D]$.
\end{enumerate}
\label{cor:equivalence gac bc}
\end{corollary}
\begin{corollary}
Let $\mathbf{X}$ and $\mathbf{Y}$ be disjoint node sets in a $\DAG$, $\CPDAG$, $\MAG$ or $\PAG$ $\g$. If $\g$ contains no possibly directed path $p = \langle V_1, \dots ,V_k \rangle$ with $k \ge 3$ such that $ \{V_{1},V_{k}\} \subseteq \mathbf{X}$ and $\{V_{2}, \dots ,V_{k-1}\} \cap \mathbf{X} = \emptyset$, or if $\g$ is a $\DAG$ or $\CPDAG$ and $\mathbf{Y} \subseteq \PossDe(\mathbf{X},\g)$, then the following statements are equivalent:
\begin{enumerate}[label=(\roman*)]
\item\label{l:1} There exists a set that satisfies the generalized adjustment criterion relative to $(\mathbf{X,Y})$ in~$\g$.
\item\label{l:2} There exists a generalized back-door set relative to $(\mathbf{X,Y})$ in~$\g$.
\end{enumerate}
\label{cor:equivalence gac gbc}
\end{corollary}

\subsection{Examples} \label{subsec:ex3}

Figure~\ref{figex2} (see Example~\ref{ex2}) in Section~\ref{subsec:examples1} shows a non-amenable graph (that is, condition \ref{nogac1} in Theorem~\ref{theorem:unified-nope} is satisfied).
Figure~\ref{fig:amenable-noset}(b) (see Example~\ref{ex:amenable-noset}) in Section~\ref{subsec:examples1} shows an amenable graph for which there is no set that satisfies the generalized adjustment criterion (that is, condition \ref{nogac1}~in~Theorem~\ref{theorem:unified-nope} is violated, but condition \ref{nogac2} is satisfied).

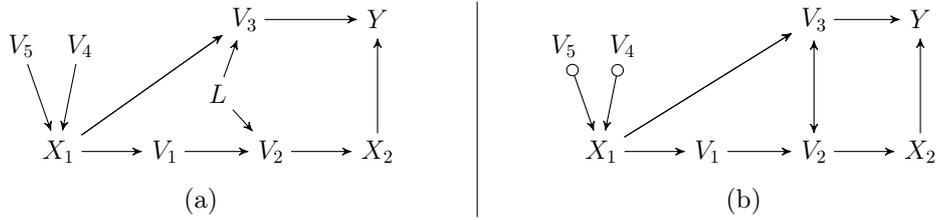
\begin{figure}
   \centering
   \begin{subfigure}{.45\textwidth}
     \centering
     \begin{tikzpicture}[>=stealth',shorten >=1pt,auto,node distance=2cm,main node/.style={minimum size=0.5cm,font=\sffamily\Large\bfseries},scale=0.7,transform shape]
   \node[main node]         (X1)                       {$X_1$};
   \node[main node]         (A)  at (.4,2)        {$V_4$};
   \node[main node]         (B)  at (-.7,2)        {$V_5$};
   \node[main node]         (V1) [right of= X1]  		{$V_{1}$};
   \node[main node]         (V2) [right of = V1]  		{$V_2$};
   \node[main node]         (L) at(3,1.1)				{$L$};
   \node[main node]         (V3)  at(3.5,2.5)			{$V_3$};
   \node[main node]         (X2) [right of= V2]  		{$X_{2}$};
   \node[main node]       	 (Y)  at (6, 2.5)      		{$Y$};
   
   \draw[->] (A) edge    (X1);
   \draw[->] (B) edge    (X1);
   \draw[->] (X1) edge    (V1);
   \draw[->] (X1) edge    (V3);
   \draw[->] (V1) edge    (V2);
   \draw[->] (V2) edge    (X2);
   \draw[->] (L) edge    (V3);
   \draw[->] (L) edge    (V2);
   \draw[->] (X1) edge    (V3);
   \draw[->] (V3) edge    (Y);
   \draw[->] (X2) edge    (Y);
   \end{tikzpicture}
    \caption{}
   \end{subfigure}
   \vrule
   \begin{subfigure}{.45\textwidth}
     \centering
     \begin{tikzpicture}[>=stealth',shorten >=1pt,auto,node distance=2cm,main node/.style={minimum size=0.5cm,font=\sffamily\Large\bfseries},scale=.7,transform shape]
   \node[main node]         (X1)                       {$X_1$};
   \node[main node]         (A)  at (.4,2)        {$V_4$};
   \node[main node]         (B)  at (-.7,2)        {$V_5$};
   \node[main node]         (V1) [right of= X1]  		{$V_{1}$};
   \node[main node]         (V2) [right of = V1]  		{$V_2$};
   \node[main node]         (V3)  at(4,2.5)			{$V_3$};
   \node[main node]         (X2) [right of= V2]  		{$X_{2}$};
   \node[main node]       	 (Y)  at (6, 2.5)      		{$Y$};
   
   \draw[o->] (A) edge    (X1);
   \draw[o->] (B) edge    (X1);
   \draw[->] (X1) edge    (V1);
   \draw[->] (X1) edge    (V3);
   \draw[->] (V1) edge    (V2);
   \draw[->] (V2) edge    (X2);
   \draw[<->] (V2) edge    (V3);
   \draw[->] (X1) edge    (V3);
   \draw[->] (V3) edge    (Y);
   \draw[->] (X2) edge    (Y);
   \end{tikzpicture}
    \caption{}
   \end{subfigure}
   \caption{\small (a) $\DAG$ $\g[D]$, (b) $\PAG$ $\g[P]$ used in Example~\ref{ex:nogbc}.}
   \label{fig:nogbc}
\end{figure}

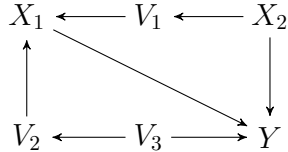
\begin{figure}
     \centering
     \begin{tikzpicture}[>=stealth',shorten >=1pt,auto,node distance=2cm,main node/.style={minimum size=0.6cm,font=\sffamily\Large\bfseries},scale=0.8,transform shape]
   \node[main node]         (X1)                       		{$X_1$};
   \node[main node]       	(V1) [right of= X1]	    		{$V_1$};
   \node[main node]         (X2) [right of= V1]				{$X_{2}$};
   \node[main node]         (Y)  [below of= X2]  			{$Y$};
   \node[main node]       	(V2) [below of= X1]	 	    	{$V_2$};
   \node[main node]       	(V3) [right of= V2]	   			{$V_3$};

   \draw[->] 	(V1)	edge    (X1);
   \draw[->] 	(X2)	edge    (V1);
   \draw[->] 	(X2)	edge    (Y);
   \draw[->] 	(V2) 	edge    (X1);
   \draw[->] 	(V3)	edge    (V2);
   \draw[->] 	(V3) 	edge    (Y);
   \draw[->] 	(X1) 	edge    (Y);
   \end{tikzpicture}
   \caption{ $\DAG$ $\g[D]$ used in Example~\ref{ex:constr-gbc-nodsep}.}
   \label{fig:constr-gbc-nodsep}
\end{figure}

We now give three additional examples. Figure~\ref{fig:nogbc1}(b) (see Example~\ref{ex:nogbc1}) shows an amenable graph for which there is no set that satisfies the generalized adjustment criterion (that is, \ref{nogac2} in Theorem~\ref{theorem:unified-nope} is satisfied). This example illustrates the result of Corollary~\ref{cor:noforbx}. Figure~\ref{fig:nogbc1}(a) (see Example~\ref{ex:nogbc1}) and Figure~\ref{fig:nogbc} (see Example~\ref{ex:nogbc}) show cases where there is a set that satisfies the generalized adjustment criterion, but there is no generalized back-door set (that is, conditions \ref{nogac1}~and~\ref{nogac2} in Theorem~\ref{theorem:unified-nope} are violated, but condition~\ref{nogbc3} is satisfied). Figure~\ref{fig:constr-gbc-nodsep} (see Example~\ref{ex:constr-gbc-nodsep}) shows an example of a $\DAG$ in which there are sets that satisfy the generalized adjustment criterion and the generalized back-door criterion, but no set satisfies Pearl's back-door criterion (that is, conditions~\ref{nogac1},~\ref{nogac2}~and~\ref{nogbc3} in Theorem~\ref{theorem:unified-nope} are violated, but condition~\ref{nobc3} is satisfied).

\begin{figure}
   \centering
   \begin{subfigure}{.5\textwidth}
   \vspace{1cm}
     \centering
     \begin{tikzpicture}[>=stealth',shorten >=1pt,auto,node distance=2cm,main node/.style={minimum size=0.5cm,font=\sffamily\Large\bfseries},scale=0.7,transform shape]
   \node[main node]         (X1)                       {$X_1$};
   \node[main node]         (V1) [right of= X1]  		{$V_{1}$};
   \node[main node]         (V2) [right of = V1]  		{$V_2$};
   \node[main node]         (X2) [below left of= V2]  		{$X_{2}$};
   \node[main node]       	 (Y1)  [right of= V2]        {$Y_1$};
    \node[main node]       	 (Y2)  [below left of= X2]        {$Y_2$};   
   \draw[->] (X1) edge    (V1);
   \draw[->] (V1) edge    (V2);
   \draw[->] (Y1) edge    (V2);
   \draw[->] (V2) edge    (X2);
   \draw[->] (X1) edge    (Y2);
   \draw[->] (X2) edge    (Y2);
   \draw[->] (Y1) ..  controls (4.5,-1.7) ..    (Y2);
   \end{tikzpicture}
    \caption{}
   \end{subfigure}
   \hspace{-1cm}
   \vrule
   \begin{subfigure}{.5\textwidth}
     \centering
     \begin{tikzpicture}[>=stealth',shorten >=1pt,auto,node distance=2cm,main node/.style={minimum size=0.5cm,font=\sffamily\Large\bfseries},scale=.7,transform shape]
   \node[main node]         (X1)                       {$X_1$};
   \node[main node]         (V1) [right of= X1]  		{$V_{1}$};
   \node[main node]         (V2) [right of = V1]  		{$V_2$};
   \node[main node]         (X2) [below left of= V2]  		{$X_{2}$};
   \node[main node]       	 (Y1)  [right of= V2]        {$Y_1$};
    \node[main node]       	 (Y2)  [below left of= X2]        {$Y_2$};   
   \draw[->] (X1) edge    (V1);
   \draw[->] (V1) edge    (V2);
   \draw[->] (Y1) edge    (V2);
   \draw[->] (V2) edge    (X2);
   \draw[->] (X1) edge    (Y2);
   \draw[->] (X2) edge    (Y2);
   \draw[->] (Y1) ..  controls (4.5,-1.7) ..    (Y2);
   \draw[->] (X1) ..  controls (2.5,1.7) and (3.5,1.7) ..   (Y1);
   \end{tikzpicture}
    \caption{}
   \end{subfigure}
   \caption{\small (a) $\DAG$ $\g[D]_1$, (b) $\DAG$ $\g[D]_2$ used in Example~\ref{ex:nogbc1}.}
   \label{fig:nogbc1}
\end{figure}
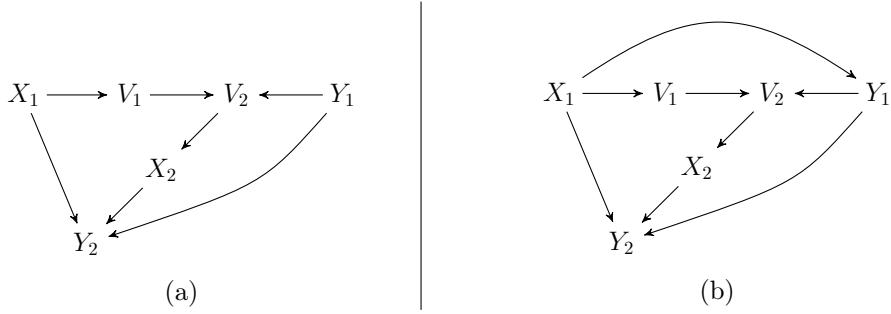
\begin{example} 
 Let $\mathbf{X} = \{X_1,X_2\}$ and $\mathbf{Y} = \{Y_1,Y_2\}$ and consider the $\DAG$s $\g[D]_1$ and $\g[D]_2$ in Figure~\ref{fig:nogbc1}(a) and~\ref{fig:nogbc1}(b) respectively.
   We first consider $\DAG$ $\g[D]_1$.
   The proper non-causal path $X_2 \leftarrow V_2 \leftarrow Y_1$ satisfies~\ref{l:nogbc0}$-$\ref{l:nogbc3} in Lemma~\ref{lemma:no gbc path}. Hence, there is no generalized back-door set relative to $(\mathbf{X},\mathbf{Y})$ in~$\g[D]_1$. However, $\{ V_1,V_2\}$ satisfies the generalized adjustment criterion relative to $(\mathbf{X},\mathbf{Y})$ in~$\g[D]_1$.

We now consider $\DAG$ $\g[D]_2$. Note that the only difference between $\g[D]_1$ and $\g[D]_2$ is the additional edge $X_1 \rightarrow Y_1$ in $\g[D]_2$. This edge implies that $Y_1 \in \fb{\g[D]_2}$. Hence, the proper non-causal path $X_2 \leftarrow V_2 \leftarrow Y_1$ satisfies~\ref{nogac2} in Theorem~\ref{theorem:unified-nope} and thus, there is no set that satisfies the generalized adjustment criterion relative to $(\mathbf{X,Y})$ in $\g[D]_2$. Then by Theorem~\ref{theorem:gac} there is no adjustment set relative to $(\mathbf{X,Y})$ in $\g[D]_2$. Since $\g[D]_2$ is a DAG such that $\mathbf{Y} \subseteq \De(\mathbf{X},\g[D]_2)$, \ref{noforbx2} in Corollary~\ref{cor:noforbx} implies that $\mathbf{X} \cap \fb{\g[D]_2} \neq \emptyset$. This is indeed true, since  $\fb{\g[D]_2} = \{ X_2,V_2,Y_1,Y_2\}$. 
   \label{ex:nogbc1}
\end{example}

\begin{example}

 Let $\mathbf{X} = \{X_1,X_2\}$ and $\mathbf{Y} = \{Y\}$ and consider $\DAG$ $\g[D]$ and $\PAG$ $\g[P]$ in Figures~\ref{fig:nogbc}(a) and~\ref{fig:nogbc}(b).
   We first consider $\DAG$ $\g[D]$. Any generalized back-door set relative to $(\mathbf{X},\mathbf{Y})$ in~$\g[D]$ must contain $L$. However, the same is not true for the generalized adjustment criterion. For example, $\{ V_1,V_2\}$ satisfies the generalized adjustment criterion relative to $(\mathbf{X},\mathbf{Y})$ in~$\g[D]$.

We now consider $\g[P]$. Note that $\g[P]$ is the PAG of $\g[D]$ when $L$ is unobserved. 
The proper non-causal path $X_2 \leftarrow V_2 \leftrightarrow V_3 \rightarrow Y$ satisfies~\ref{l:nogbc0}$-$\ref{l:nogbc3} in Lemma~\ref{lemma:no gbc path}. Hence, there is no generalized back-door set relative to $(\mathbf{X},\mathbf{Y})$ in~$\g[P]$.  However, the sets $\{ V_1,V_2\}$, $\{V_1,V_2,V_4\}$, $\{V_1,V_2,V_5\}$, $\{V_1,V_2,V_4,V_5\}$ all satisfy the generalized adjustment criterion relative to $(\mathbf{X},\mathbf{Y})$ in~$\g[P]$.
   \label{ex:nogbc}
\end{example}

\begin{example}
  Let $\mathbf{X} = \{X_1,X_2\}$ and $\mathbf{Y} = \{ Y \}$ and consider $\DAG$ $\g[D]$ in Figure~\ref{fig:constr-gbc-nodsep}.
The non-causal path $X_1 \leftarrow V_1 \rightarrow X_2 \rightarrow Y$ satisfies~\ref{l:gbc-bc1}$-$\ref{l:gbc-bc4} in Lemma~\ref{lemma:gbc bc}. Hence, no set can satisfy the back-door criterion relative to $(\mathbf{X},\mathbf{Y})$ in~$\g[D]$. However, $\{V_2\}$, $\{V_3\}$, $\{V_2,V_3\}$, $\{V_1,V_2\}$, $\{V_1,V_3\}$, and $\{V_1, V_2, V_3\}$ all satisfy the generalized adjustment criterion and $\{V_2\}$, $\{V_3\}$ and $\{V_2,V_3\}$ all satisfy the generalized back-door criterion relative to $(\mathbf{X},\mathbf{Y})$ in~$\g[D]$.
\label{ex:constr-gbc-nodsep}
\end{example}
\section{Discussion} \label{sec:discussion}

We have derived a generalized adjustment criterion that is sound and complete for adjustment in $\DAG$s,
$\MAG$s, $\CPDAG$s and $\PAG$s (see Definition~\ref{def:gac}, Theorem~\ref{theorem:gac}). 
This is relevant in practice,
in particular in combination with algorithms that can learn $\CPDAG$s or $\PAG$s from observational data.

In addition to the criterion itself, we have also given all necessary ingredients for implementing efficient algorithms to test the criterion for a given set and to construct all sets that fulfill it, or to learn that no set fulfilling the criterion exists. Thus, we obtain a complete generalization of the algorithmic framework for $\DAG$s and $\MAG$s by \citet{vanconstructing} to $\CPDAG$s or $\PAG$s.
In this sense, our work presented in this paper is a theoretical contribution that closes the chapter on covariate adjustment for $\DAG$s, $\CPDAG$s, $\MAG$s and $\PAG$s without selection variables.

 \cite{correa2017causal} define necessary and sufficient graphical conditions for covariate adjustment in latent projection graphs in the presence of selection variables. Future work may explore how to extend their results to MAGs and PAGs with selection variables. 
Other future work may study the estimation accuracy of estimators based on different adjustment sets.
Any valid adjustment set can be used to produce an unbiased estimator of the total causal effect. However, the efficiency of the estimators induced by distinct adjustment sets varies \citep{greenland1999causal,kuroki2003covariate,kuroki04selection,hahn2004functional,hahn1998role,guo2010sufficient,DeLunaEtAl11}. Moreover, \cite{kuroki2003covariate} and \cite{kuroki04selection} indicate that a minimal adjustment set does not necessarily lead to the most efficient estimator.  
Defining a best adjustment set in terms of efficiency is still an open question.

\paragraph*{Acknowledgements.} This work was supported in part by Swiss NSF Grants 200021\_149760 and 200021\_172603.

\section*{Appendix}

\appendix

\section{Preliminaries} \label{subsec:additional}

We first state various existing results and definitions.

\textbf{Adjacencies, discriminating paths and d-separations.} The set of all nodes adjacent to $X$ in a graph $\g$ is denoted by $\Adj(X,\g)$. A path $p = \langle X,\dots,Z,V,Y \rangle$ is a \textit{discriminating path} from $X$ to $Y$ for $V$ in graph $\g$, if it consists of at least four nodes, $X$ is not adjacent to $Y$ in $\g$ and every non-endpoint node on $p(X,V)$ is a collider on $p$ and a parent of $Y$.
If $\mathbf{X}$ and  $\mathbf{Y}$ are d-separated given $\mathbf{Z}$ in a $\DAG$ $\g[D]$, we write $\mathbf{X} \dsepp \mathbf{Y} \given \mathbf{Z}$.

\begin{definition} (\textbf{Distance-from-$\mathbf{Z}$})
Let $\mathbf{X,Y}$ and $\mathbf{Z}$ be pairwise disjoint node sets in a $\DAG$, $\CPDAG$, $\MAG$ or $\PAG$ $\g$. Let $p$ be a path between $\mathbf{X}$ and $\mathbf{Y}$ in~$\g$ such that every collider $C$ on $p$ has a possibly directed path (possibly of length $0$) to $\mathbf{Z}$. Define the $\distancefrom{\mathbf{Z}}$ of $C$ to be the length of a shortest possibly directed path (possibly of length $0$) from $C$ to $\mathbf{Z}$, and define the $\distancefrom{\mathbf{Z}}$ of $p$ to be the sum of the distances from $\mathbf{Z}$ of the colliders on $p$.
\label{def:distance from Z}
\end{definition}

\noindent{}If $\g$ is a $\MAG$ and $p$ is a path from $\mathbf{X}$ to $\mathbf{Y}$ that is m-connecting given $\mathbf{Z}$ in~$\g$, then Definition~\ref{def:distance from Z} reduces to the notion of \emph{distance-from-$\mathbf{Z}$} in \citet[][p213]{zhang2006causal}.

\begin{theorem} (\textbf{Wright's rule} \citealp[cf.][]{wright1921correlation}) Let $\mathbf{X} = \mathbf{A}\mathbf{X} + \mathbf{\epsilon}$, where $\mathbf{A} \in \mathbb{R}^{k \times  k}$, $\mathbf{X}$ $= (X_1,\dots, X_k)^T$  and $\mathbf{\epsilon} = (\epsilon_1,\dots, \epsilon_k)^T$ is a vector of mutually independent errors with means zero. Moreover, let $\Var(\mathbf{X}) = \mathbf{I}$.
Let $\g[D] = (\mathbf{X},\mathbf{E})$, be the corresponding $\DAG$ such that $X_i \rightarrow X_j$ in $\g[D]$ if and only if $A_{ji} \neq 0$. A nonzero entry $A_{ji}$ is called the edge coefficient of $X_i \rightarrow X_j$.
For two distinct nodes $X_i,X_j \in \mathbf{X}$, let $p_1, \dots, p_r$ be all paths between $X_i$ and $X_j$ in $\g[D]$ that do not contain a collider. Then $\Cov(X_i,X_j) = \sum_{s=1}^{r}\pi_s$, where $\pi_s$ is the product of all edge coefficients along path $p_s$, $s \in \{1,\dots, r\}$.
\label{theorem:wright}
\end{theorem}
\begin{theorem} \citep[cf.\ Theorem 3.2.4][p63]{mardia1980multivariate}
Let $\mathbf{X} = (\mathbf{X_1}^T,\mathbf{X_2}^T)^T$ be a $p$-dimensional multivariate Gaussian random vector with mean vector $\mathbf{\mu} = (\mathbf{\mu_1}^T,\mathbf{\mu_2}^T)^T$ and covariance matrix $\mathbf{\Sigma} = \begin{bmatrix}
\mathbf{\Sigma_{11}} & \mathbf{\Sigma}_{12} \\
\mathbf{\Sigma_{21}} & \mathbf{\Sigma}_{22}
\end{bmatrix}$, so that $\mathbf{X_1}$ is a $q$-dimensional multivariate Gaussian random vector with mean vector $\mathbf{\mu_1}$ and covariance matrix $\mathbf{\Sigma}_{11}$ and  $\mathbf{X_2}$ is a $(p-q)$-dimensional multivariate Gaussian random vector with mean vector $\mathbf{\mu_2}$ and covariance matrix $\mathbf{\Sigma}_{22}$.
Then $E[\mathbf{X_2} \given \mathbf{X_1 = x_1}] = \mathbf{\mu_2} + \mathbf{\Sigma}_{21}\mathbf{\Sigma}_{11}^{-1} (\mathbf{x_1} - \mathbf{\mu_1})$.
\label{theorem:mardia-condexp}
\end{theorem}

\begin{definition}{(\textbf{Moralization}; \citealp{lauritzen1988local})} 
Let $\g[D]$ be a $\DAG$. The \textit{moral} graph $\g[D]^m$ is formed by adding the edge $A-B$ to any structure of the form $A \rightarrow C \leftarrow B$ with $A \notin \Adj(B,\g[D])$ (marrying unmarried parents) and subsequently making all edges in the resulting graph undirected.
\label{def:moralization}
\end{definition}

\begin{definition}{(\textbf{Induced subgraph})} 
Let $\mathbf{X} \subseteq \mathbf{V}$ be a node set in a $\DAG$ $\g[D] = (\mathbf{V},\mathbf{E})$. Then $\g[D]_\mathbf{X} = (\mathbf{X},\mathbf{E}_\mathbf{X})$, where $\mathbf{E}_\mathbf{X}$ consists of all edges in $\mathbf{E}$ for which both endpoints are in $\mathbf{X}$, is the induced subgraph of $\g[D]$ on $\mathbf{X}$.
\label{def:induced-subgraph}
\end{definition}

\begin{theorem} (\textbf{Reduction of d-separation to node cuts}; \citealp[cf.\ Proposition 3 in][]{lauritzen1990independence}, \citealp[cf.\ Corollary 2 in][]{richardson2003markov}) 
Let $\g[D]$ be a $\DAG$ and let $\mathbf{X},\mathbf{Y}$ and $\mathbf{Z}$ be pairwise disjoint node sets in~$\g[D]$. Then $\mathbf{Z}$ d-separates $\mathbf{X}$ and $\mathbf{Y}$ in~$\g[D]$ if and only if all paths between $\mathbf{X}$ and $\mathbf{Y}$ in $(\g[D]_{\An(\mathbf{X} \cup \mathbf{Y} \cup \mathbf{Z},\g[D])})^m$ contain at least one node in $\mathbf{Z}$.
\label{theorem:moralization}
\end{theorem}

\begin{lemma}
   {(\textbf{Basic property of $\CPDAG$s and $\PAG$s}; \citealp[cf.\ Lemma~1 in][]{meek1995causal}, \citealp[cf.\ Lemma~3.3.1 in][]{zhang2006causal})}
   Let $X,Y$ and $Z$ be distinct nodes in a $\CPDAG$ or $\PAG$ $\g$. If $X \bulletarrow Y \circbullet Z$, then there is an edge between $X$ and $Z$ with an arrowhead at $Z$. Furthermore, if the edge between $X$ and $Y$ is $X \rightarrow Y$, then the edge between $X$ and $Z$ is either $X \circarrow Z$ or $X \rightarrow Z$ (that is, not $X \leftrightarrow Z$).
   \label{lemma:basic property of cpdags and pags}
\end{lemma}

\begin{lemma}{\citep[cf.\ Lemma~1 in][]{richardson2003markov}} Let $\mathbf{X,Y}$ and $\mathbf{Z}$ be pairwise disjoint node sets in a $\DAG$ or $\MAG$ $\g$. If there is a path $p$ from $X \in \mathbf{X}$ to $Y \in \mathbf{Y}$, on which no non-collider is in $\mathbf{Z}$ and every collider on $p$ is in $\An(\mathbf{X} \cup \mathbf{Y} \cup \mathbf{Z}, \g)$, then there exists a path $q$ from $X' \in \mathbf{X}$ to $Y' \in \mathbf{Y}$ that is m-connecting given $\mathbf{Z}$ in~$\g$.
\label{lemma:rich}
\end{lemma}

\begin{lemma}
   \citep[Lemma~0 in][p208]{zhang2006causal}
  Let $X$ and $Y$ be distinct nodes in a $\MAG$ $\g[M]$. If $p = \langle X, \dots ,Z,V,Y \rangle$ is a discriminating path from $X$ to $Y$ for $V$ in a $\MAG$ $\g[M]$, and the corresponding path to $p(X,V)$ in the $\PAG$ $\g[P]$ of $\g[M]$ is (also) a collider path, then the corresponding path to $p$ in $\g[P]$ is also a discriminating path from $X$ to $Y$ for $V$.
   \label{lemma0}
\end{lemma}

\begin{lemma} \citep[cf.\ Lemma~1 in][p208]{zhang2006causal}
    Let $X$ and $Y$ be distinct nodes and let $\mathbf{Z}$ be a node set that does not contain $X$ and $Y$ in a $\MAG$ $\g[M]$ ($\DAG$ $\g[D]$). Let $p$ be a shortest path from $X$ to $Y$ that is m-connecting given $\mathbf{Z}$ in $\g[M]$ $(\g[D])$. Let $\g[P]$ be the $\PAG$ of $\g[M]$ ($\CPDAG$ of $\g[D]$) and let $\pstar$ in $\g[P]$ be the corresponding path to $p$ in $\g[M]$ $(\g[D])$. Then $\pstar$ is a definite status path in~$\g[P]$.
   \label{lemma1}
\end{lemma}

\begin{lemma} \citep[cf.\ Lemma~2 in][p213]{zhang2006causal}
    Let $X$ and $Y$ be distinct nodes and let $\mathbf{Z}$ be a node set that does not contain $X$ and $Y$ in a $\MAG$ $\g[M]$ ($\DAG$ $\g[D]$). Let $p$ be a shortest path from $X$ to $Y$ that is m-connecting given $\mathbf{Z}$ in $\g[M]$ $(\g[D])$ such that no equally short m-connecting path has a shorter distance-from-$\mathbf{Z}$ (see Definition~\ref{def:distance from Z}) than $p$ does. Let $\g[P]$ be the $\PAG$ of $\g[M]$ ($\CPDAG$ of $\g[D]$) and let $\pstar$ in $\g[P]$ be the corresponding path to $p$ in $\g[M]$ $(\g[D])$. Then $\pstar$ is a definite status path from $X$ to $Y$ that is m-connecting given $\mathbf{Z}$ in $\g[P]$.
   \label{lemma2}
\end{lemma}

\begin{lemma}{\citep[cf.\ Lemma~B.1 in][]{zhang2008completeness}}
  Let $X$ and $Y$ be distinct nodes in a $\CPDAG$ or $\PAG$ $\g$. If $p = \langle X, \dots ,Y \rangle$ is a possibly directed path from $X$ to $Y$ in~$\g$, then some subsequence of $p$ forms an unshielded possibly directed path from $X$ to $Y$ in~$\g$.
   \label{lemma:unshielded}
\end{lemma}

\begin{lemma}
   {(\citealp[cf.\ Lemma~B.2 in][]{zhang2008completeness},\citealp[ Lemma 7.2 in][]{maathuis2013generalized})}
   Let $X$ and $Y$ be distinct nodes in a $\CPDAG$ or $\PAG$ $\g$. If $p = \langle X=V_0, \dots , V_k=Y \rangle$, $k \geq 2$ is an unshielded possibly directed path from $X$ to $Y$ in~$\g$, and $V_{i-1} \bulletarrow V_{i}$ for some $i \in \{1,\dots,k\}$, then  $V_{j-1} \rightarrow V_{j}$ for all $j \in \{i+1,\dots,k\}$.
   \label{lemma:unshielded edges}
\end{lemma}

\begin{lemma}
   {(\citealp[cf.\ Theorem~2 in][]{zhang2008completeness}, \citealp[Lemma~7.6 in][]{maathuis2013generalized})} Let $\g$ be a $\PAG$ $(\CPDAG)$. Let $\g[M]$ $(\g[D])$ be the graph resulting from the following procedure applied to a $\g$:
 \begin{enumerate}[label=(\arabic*)]
\item replace all partially directed edges $\circarrow$ in~$\g$ with directed edges $\rightarrow$, and
\item\label{orient1} orient the subgraph of $\g$ consisting of all non-directed edges $\circcirc$ into a $\DAG$ with no unshielded colliders.
\end{enumerate}
Then $\g[M]$ $(\g[D])$ is in $[\g]$. Moreover, if $X$ is a node in $\g$, then one can always find an orientation of \ref{orient1} that does not create any new edges into $X$.
   \label{lemmamarlrx}
\end{lemma}

\begin{lemma}
   {\citep[Lemma~7.5 in][]{maathuis2013generalized}} Let $X$ and $Y$ be two distinct nodes in a $\DAG$, $\CPDAG$, $\MAG$ or $\PAG$ $\g$. Then $\g$ cannot have both a possibly directed path from $X$ to $Y$ and an edge of the form $Y \bulletarrow X$.
   \label{lemmamarlcycle}
\end{lemma}

\begin{definition}{(\textbf{$\dsep{\g}$}; \citealp[][p136]{spirtes2000causation})} Let $X$ and $Y$ be two distinct nodes in a $\DAG$ or $\MAG$ $\g$. We say that $V \in \dsep{\g}$ if $V \neq X$ and there is a collider path between $X$ and $V$ in~$\g$ such that every node on this path is an ancestor of $X$ or $Y$ in~$\g$.
\label{def:dsep}
\end{definition}

\begin{definition}($\g[R]$ and $\g[R]_{\underline{X}}$; \citealp{maathuis2013generalized}) 
Let $X$ be a node in a $\DAG$, $\CPDAG$, $\MAG$ or $\PAG$ $\g$.
Let $\g[R]$ be a $\DAG$ or $\MAG$ represented by $\g$, in the following sense. If $\g$ is a $\DAG$ or $\MAG$, let $\g[R] = \g$. If $\g$ is a $\CPDAG$ $(\PAG)$, let $\g[R]$ be a $\DAG$ $(\MAG)$ in $[\g]$ as defined in Lemma~\ref{lemmamarlrx}, so that $\g[R]$ has the same number of edges into $X$ as $\g$.
Let $\g[R]_{\underline{X}}$ be the graph obtained from $\g[R]$ by removing all directed edges out of $X$ that are visible in~$\g$.
\label{def:rx}
\end{definition}

\begin{theorem} \citep[Theorem~4.1 in][]{maathuis2013generalized}
Let $X$ and $Y$ be distinct nodes in a $\DAG$, $\CPDAG$, $\MAG$ or $\PAG$ $\g$. Let $\g[R]$ and $\g[R]_{\underline{X}}$ be defined as in Definition
\ref{def:rx}. Then there exists a generalized back-door set relative to $(X,Y)$ in~$\g$ if and only if $Y \notin \Adj(X,\g[R]_{\underline{X}})$ and $
\dsep{\g[R]_{\underline{X}}} \cap \PossDe(X,\g) = \emptyset$. Moreover, if such a set exists, then $\dsep{\g[R]_{\underline{X}}}$ is a generalized back-door set
relative to $(X,Y)$ in~$\g$.
\label{theorem:marloes constructive gbc}
\end{theorem}

\subsection{Rules of the Do-calculus \citep[][Chapter~3.4]{Pearl2009}}

Let $\mathbf{X',Y',Z',W'}$ be pairwise disjoint (possibly empty) sets of nodes in a causal $\DAG$ $\g[D]$.
Let $\g[D]_{\overline{\mathbf{X'}}}$ denote the graph obtained by deleting all edges into $\mathbf{X'}$ from $\g[D]$. Similarly, let $\g[D]_{\underline{\mathbf{X'}}}$ denote the graph obtained by deleting all edges out of $\mathbf{X'}$ in~$\g[D]$ and let $\g[D]_{\overline{\mathbf{X'}}\underline{\mathbf{Z'}}}$ denote the graph obtained by deleting all edges into $\mathbf{X'}$ and all edges out of $\mathbf{Z'}$ in~$\g[D]$.
Then the following three rules are valid for every density function consistent with $\g[D]$.

\textbf{Rule 1} (Insertion/deletion of observations) If $\mathbf{Y'} \dsepp \mathbf{Z'} \given \mathbf{X'} \cup \mathbf{W'}$ in~$\g[D]_{\overline{\mathbf{X'}}}$, then
\begin{align}
f(\mathbf{y'} \given \ddo(\mathbf{x'}),\mathbf{w'}) = f(\mathbf{y'} \given \ddo(\mathbf{x'}),\mathbf{z',w'}). \label{rule1}
\end{align}

\textbf{Rule 2} (Action/observation exchange) If $\mathbf{Y'} \dsepp \mathbf{Z'} \given \mathbf{X'} \cup \mathbf{W'}$ in~$\g[D]_{\overline{\mathbf{X'}}\underline{\mathbf{Z'}}}$, then
\begin{align} \label{rule2}
\begin{split}
f(\mathbf{y'} \given \ddo(\mathbf{x'}),\ddo(\mathbf{z'}),\mathbf{w'}) & = f(\mathbf{y'} \given \ddo(\mathbf{x'}),\mathbf{z',w'}).
\end{split}
\end{align}

\textbf{Rule 3} (Insertion/deletion of actions)If $\mathbf{Y'} \dsepp \mathbf{Z'} \given \mathbf{X'} \cup \mathbf{W'}$ in~$\g[D]_{\overline{\mathbf{X' \cup Z'(W')}}}$, then
\begin{align} \label{rule3}
\begin{split}
 f(\mathbf{y'} \given \ddo(\mathbf{x'}),\mathbf{w'}) & = f(\mathbf{y'} \given \ddo(\mathbf{x'}),\ddo(\mathbf{z'}),\mathbf{w'}),
\end{split}
\end{align}
where $\mathbf{Z'(W')}$ denotes the set of $\mathbf{Z'}$-nodes that are not ancestors of any $\mathbf{W'}$ node in~$\g[D]_{\overline{\mathbf{X'}}}$.
If $\mathbf{W'} = \emptyset$, then $\mathbf{Z'}(\mathbf{W'}) = \mathbf{Z'}$.

\subsection{FCI Orientation Rules \citep[][p183]{spirtes2000causation}}

Let $A,B,C$ and $D$ be distinct nodes in a $\PAG$ $\g[P]$.
Below, we give the first $4$ orientation rules of the FCI algorithm defined in \cite{spirtes2000causation}.

\noindent\textbf{R1} If $A \bulletarrow B \circbullet C$, and $A$ and $C$ are not adjacent, then orient the triple $\left\langle A,B,C \right\rangle$ as $A \bulletarrow B \rightarrow C$.

\noindent\textbf{R2} If $A \rightarrow B \bulletarrow C$ or $A \bulletarrow B \rightarrow C$ and $A \bulletcirc C$, then orient $A \bulletcirc C$ as $A \bulletarrow C$.

\noindent\textbf{R3} If $A \bulletarrow B  \arrowbullet C$, $A \bulletcirc D \circbullet C$, $A$ and $C$ are not adjacent, and $D \bulletcirc B$, then orient $D \bulletcirc B$ as $D \bulletarrow B$.

\noindent\textbf{R4} If $\left\langle D, \dots ,A,B, C \right\rangle$ is a discriminating path from $D$ to $C$ for $B$ and $B \circbullet C$, then orient $B \circbullet C$ as $B \rightarrow C$ if $B$ is in the separation set of $D$ and $C$; otherwise orient the triple $\left\langle A,B,C \right\rangle$ as $A \leftrightarrow B \leftrightarrow C$.

\noindent{}These four rules were proven to be sound in \cite{spirtes2000causation}, meaning that edge marks oriented by these rules correspond to invariant edge marks in the maximally informative $\PAG$ for the true causal $\MAG$.
Six additional orientation rules for the FCI algorithm were defined in \cite{zhang2008completeness}.
The augmented FCI algorithm, including all ten orientation rules was proven to be sound and complete in \cite{zhang2008completeness}.

\section{Proofs for Section \ref{sec:main}} \label{subsec:proofs GAC}

\begin{figure}
    \centering
      \begin{tikzpicture}[>=stealth',shorten >=1pt,node distance=3cm, main node/.style={minimum size=0.4cm}]
      [>=stealth',shorten >=1pt,node distance=3cm,initial/.style    ={}]
    \node[main node]         (T34) {\textbf{Theorem~\ref{theorem:gac}}};
    \node[main node]         (L53)  [left of= T34]          {Lemma~\ref{lemmaprovingeqofacconditiona}};
    \node[main node,yshift=1cm]  (L52) at (L53) {Lemma~\ref{lemma:adj amen gen}};
    \node[main node]         (L31)  [left of= L52,align=center]          {Lemma~\ref{lemmadirectedpath}};
    \node[main node]         (L46)  [left of= L31,align=center]          {Lemma~\ref{lemma:stays_invisible}};
    \node[main node,yshift=-1cm] (L54) at (L53) {Lemma~\ref{lemmaequivalenceofcondb}};
    \node[main node]            (L56)  [left of= L54]   {Lemma~\ref{lemma:2}};
    \node[main node]            (L57)  [left of= L56,align=center]   {Lemma~\ref{lemma:1}};
    \node[main node]            (L58)  [left of= L57,align=center]   {Lemma~\ref{lemma:non-causal path}};
    \draw[->] (-2.1,1) to    (-1.2,.3);
    \draw[->] (-2,-1) to    (-1.2,-.3);
    \draw[->] (L46) edge    (L31);
    \draw[->] (L53) edge    (T34);
    \draw[->] (L31) edge    (L52);
    \draw[->] (L56) edge    (L54);
    \draw[->] (L57) edge    (L56);
    \draw[->] (L58) edge    (L57);
    \end{tikzpicture}
   \caption{Proof structure of Theorem~\ref{theorem:gac}.}
    \label{figproof}
\end{figure}
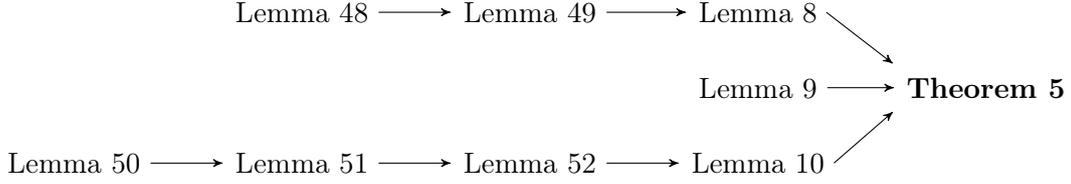

Figure~\ref{figproof} shows how all lemmas in this section fit together to prove Theorem~\ref{theorem:gac}.

\vspace{.3cm}

\begin{lemma}
Let $X$ be a node in a $\PAG$ $\g[P]$.
Let $\g[M]$ be a $\MAG$ $\g[M]$ in $[\g[P]]$ that satisfies Lemma~\ref{lemmamarlrx}. Then any edge that is either $X \circcirc Y$, $X \circarrow Y$ or invisible $X \rightarrow Y$ in~$\g[P]$ is invisible $X \rightarrow Y$ in~$\g[M]$.
\label{lemma:stays_invisible}
\end{lemma}

\begin{proofof}[Lemma~\ref{lemma:stays_invisible}]
Let $\g[M]$ be a $\MAG$ in $[\g[P]]$ that satisfies Lemma~\ref{lemmamarlrx}.
Then the edge $X \circcirc Y$, $X \circarrow Y$, or invisible $X \rightarrow Y$ in $\g[P]$ corresponds to edge $X \rightarrow Y$ in~$\g[M]$. It is left to prove that $X \rightarrow Y$ is invisible in $\g[M]$ in all these cases.

Suppose for a contradiction that $X \rightarrow Y$ is visible in~$\g[M]$. Then there is a node $D \notin \Adj(Y,\g[M])$ such that (1) $D \bulletarrow X$ is in~$\g[M]$, or (2) there is a collider path $\langle D,D_1,\dots,D_k,X \rangle,$ $k \ge 1$, into $X$ such that every $D_i, 1 \le i \le k$ is a parent of $X$ in~$\g[M]$.  We consider these cases separately and show that we arrive at a contradiction.

(1) Since $D \bulletarrow X$ is in $\g[M]$, the choice of $\g[M]$ implies that $D \bulletarrow X$ is in $\g[P]$. Since $D \notin \Adj(Y,\g[P])$,  $X \rightarrow Y$ must be in $\g[P]$, since otherwise rule \textbf{R1} of the FCI algorithm from \citet{spirtes2000causation} (see Appendix~\ref{subsec:additional}) would have been applied in $\g[P]$. But then $X \rightarrow Y$ is visible in $\g[P]$, which is a contradiction.

(2) Path $p = \langle D, D_1,\dots D_k, X, Y \rangle$ is a discriminating path from $D$ to $Y$ for $X$, that is into $X$ in~$\g[M]$. Let $\pstar$ be the path in $\g[P]$ corresponding to $p$ in $\g[M]$. Then since $p(D_1,X)$ contains only bi-directed edges in $\g[M]$, the choice of $\g[M]$ implies that $\pstar(D_1,X)$ also contains only bi-directed edges in~$\g[P]$.
Since $D \bulletarrow D_1$ is in $\g[M]$, $D \circcirc D_1$ or $D \bulletarrow D_1$ is in $\g[P]$.

Suppose first that $D \bulletarrow D_1$ is in $\g[P]$, then by Lemma~0 from \cite{zhang2006causal} (see Lemma~\ref{lemma0}), $\pstar$ is a discriminating path from $D$ to $Y$ for $X$, that is into $X$ in~$\g[P]$. Hence, $X \rightarrow Y$ is in~$\g[P]$, since otherwise rule \textbf{R4} in Appendix~\ref{subsec:additional} would have been applied. But then $X \rightarrow Y$ is a visible edge in $\g[P]$, which is a contradiction.

Next, suppose that $D \circcirc D_1$ is in $\g[P]$.
Since $D_1 \leftrightarrow D_2$ is in~$\g[P]$, by Lemma~\ref{lemma:basic property of cpdags and pags}
$D \arrowbullet D_2$ is in~$\g[P]$. This edge cannot be $D \arrowcirc D_2$ or $D \leftarrow D_2$, otherwise $D \arrowcirc D_1$ would be in $\g[P]$ (Lemma~\ref{lemma:basic property of cpdags and pags}, or \textbf{R2} of the FCI orientation rules in Appendix~\ref{subsec:additional}), contrary to our assumption.
Hence, the edge $D \leftrightarrow D_2$ is in~$\g[P]$. Then $D \leftrightarrow D_2$ is also in $\g[M]$ and $p_1  = \langle D, D_2,\dots D_k, X, Y \rangle$ is a discriminating path from $D$ to $Y$ for $X$, that is into $X$ in~$\g[M]$. Additionally, since $D \leftrightarrow D_2 \leftrightarrow \dots \leftrightarrow D_k \leftrightarrow X$ is in $\g[P]$, by Lemma~\ref{lemma0} the path $\pstar[p_1] =  \langle D, D_2,\dots D_k, X, Y \rangle$ is a discriminating path from $D$ to $Y$ for $X$, that is into $X$. But then $X \rightarrow Y$ is a visible edge in $\g[P]$, which is a contradiction.
\end{proofof}

\begin{lemma}
Let $X$ and $Y$ be distinct nodes in a $\PAG$ $\g[P]$ such that there is a possibly directed path $\pstar$ from $X$ to $Y$ in~$\g[P]$ that does not start with a visible edge out of $X$. Then there is a $\MAG$ $\g[M]$ in $[\g[P]]$ such that the path $p$ in~$\g[M]$, consisting of the same sequence of nodes as $\pstar$ in~$\g[P]$, contains a subsequence $p'$ that is a directed path from $X$ to $Y$ starting with an invisible edge in~$\g[M]$. In other words, $\g[M]$ violates \amen{} relative to $({X,Y})$.
    \label{lemmadirectedpath}
\end{lemma}

\begin{proofof}[Lemma~\ref{lemmadirectedpath}]
If $\g[P]$ violates \amen{} relative to $(X,Y)$, then there is a possibly directed path $\pstar$ from $X$ to $Y$ in~$\g[P]$ that does not start with a visible edge out of $X$.
  Let $\g[M]$ be a $\MAG$ in $[\g[P]]$ that satisfies Lemma~\ref{lemmamarlrx}. 
  We will show that $\g[M]$ violates \amen{} relative to $(X,Y)$.

    Let $\pstar[p']$ be a shortest subsequence of $\pstar$ such that $\pstar[p']$ is also a possibly directed path from $X$ to $Y$ in~$\g[P]$ that does not start with a visible edge out of $X$. We write $\pstar[p'] = \left\langle X = V_0,V_{1}, \dots ,V_r= Y \right\rangle, r \ge 1$.
   Let $p'$ in~$\g[M]$ be the path corresponding to $\pstar[p']$ in~$\g[P]$.
  Then by Lemma~\ref{lemma:stays_invisible}, $X \rightarrow V_1$ is an invisible edge in $\g[M]$. It is left to show that $p'$ is a directed path from $X$ to $Y$ in $\g[M]$.

  Suppose first that $\pstar[p']$ is a definite status path in~$\g[P]$. Then all non-endpoint nodes on $\pstar[p']$ are definite non-colliders. Hence,  $X \rightarrow V_1$ in $\g[M]$ implies that all the remaining edges on $p'$ are oriented towards $Y$.

  Else, $\pstar[p']$ is not of definite status in~$\g[P]$ and $r \ge 2$.
   Since $X \to V_1$ is in $\g[M]$, it is sufficient to show that $p'(V_{1},Y)$ is a directed path from $V_{1}$ to $Y$.
   Note that by the choice of $\pstar[p']$, $\pstar[p'](V_{1},Y)$ is a shortest possibly directed path from $V_1$ to $Y$ in $\g[P]$. Hence, it is unshielded (Lemma~B.1 from \cite{zhang2008completeness}, see Lemma~\ref{lemma:unshielded} in Appendix~\ref{subsec:additional}).
   If $V_1 \circarrow V_2$ or $V_1 \rightarrow V_2$ is in $\g[P]$, then by the choice of $\g[M]$ (Lemma~\ref{lemmamarlrx}), $V_1 \to V_2$ is in~$\g[M]$. Additionally, since $p'(V_{1},Y)$ is a possibly directed definite status path, $V_1 \rightarrow V_2$ in $\g[M]$ implies that all the remaining edges on $p'(V_1,Y)$ are oriented towards $Y$.

  Otherwise, $V_1 \circcirc V_2$ is in $\g[P]$.
   Path $\pstar[p']$ is not of definite status, whereas $\pstar[p'](V_{1},Y)$ is of definite status, as it is unshielded. Thus, node $V_1$ is not of definite status on $\pstar[p']$ and $X \in \Adj(V_2,\g[P])$. The edge $X \arrowbullet V_2$ is not in~$\g[P]$ since $\pstar[p'](X,V_2)$ is a possibly directed path from $X$ to $V_2$  in~$\g[P]$ (Lemma~7.5 from \citet{maathuis2013generalized}, see Lemma~\ref{lemmamarlcycle} in Appendix~\ref{subsec:additional}). Since $\pstar[p']$ is a shortest possibly directed path from $X$ to $Y$ in $\g[P]$ that does not start with a visible edge out of $X$, and $X \arrowbullet V_2$ is not in $\g[P]$, it follows that $X \to V_2$ is visible in~$\g[P]$.
   Since $X \to V_2$ is visible, there is a node $D \notin \Adj(V_2,\g[P])$ such that (1) $D \bulletarrow X$ is in $\g[P]$, or (2) there is a collider path $\langle D,D_1,\dots,D_k,X \rangle, k \ge 1$, that is into $X$ in $\g[P]$ such that every $D_i, 1 \le i \le k$ is a parent of $V_{2}$ in~$\g[P]$. We consider these cases separately and show that we arrive at a contradiction, implying that $\pstar[p'](V_1,Y)$ cannot start with $V_1 \circcirc V_2$.

(1) A node $D \notin \Adj(V_2,\g[P])$ such that $D \bulletarrow X$ is in $\g[P]$. Since $D \bulletarrow X$ and $X \circcirc V_1$, $X \circarrow V_1$ or $X \rightarrow V_1$ is invisible in $\g[P]$, by Lemma~\ref{lemma:basic property of cpdags and pags} and the definition of visibility, an edge between $D$ and $V_1$ is in $\g[P]$. This edge is of type $D \bulletarrow V_1$, since otherwise both a possibly directed path $\langle X,V_1,D \rangle$ and $D \bulletarrow X$ are in~$\g[P]$ (contrary to Lemma~\ref{lemmamarlcycle}). Then $D \bulletarrow V_1 \circcirc V_2$ is in $\g[P]$ and Lemma~\ref{lemma:basic property of cpdags and pags} implies that $D \in \Adj(V_2,\g[P])$, a contradiction.

(2) There is a node $D \notin \Adj(V_2,\g[P])$ and a collider path $\langle D,D_1,\dots,D_k,X \rangle, k \ge 1$, into $X$ such that every $D_i, 1 \le i \le k$ is a parent of $V_{2}$ in~$\g[P]$. Paths $D_{i} \to V_{2} \circcirc V_{1}$, $i=1,\dots,k$ are in~$\g[P]$, so by Lemma~\ref{lemma:basic property of cpdags and pags} either $D_{i} \circarrow V_{1}$ or $D_{i} \rightarrow V_{1}$ is in~$\g[P]$, for $i =1,\dots,k$. 
   If $D_{1} \circarrow V_{1}$ is in~$\g[P]$, then $ D \bulletarrow D_1 \circarrow V_1$ implies $D \bulletarrow V_1$ is also in~$\g[P]$ (Lemma~\ref{lemma:basic property of cpdags and pags}). But, then $D \bulletarrow V_1 \circcirc V_2$ implies $D \in \Adj(V_2,\g[P])$ (Lemma~\ref{lemma:basic property of cpdags and pags}), a contradiction.
  Hence, $D_{1} \to V_{1}$ is in~$\g[P]$.

  This allows us to deduce that $D \notin \Adj(V_1,\g[P])$, otherwise $D \bulletarrow D_1 \to V_1$ would imply $D \bulletarrow V_1$ (Lemma~\ref{lemmamarlcycle}) and we arrive at the contradiction $D \in \Adj(V_2,\g[P])$, as above.
   Hence, $\langle D, D_1, D_{2}, V_1 \rangle$ is a discriminating path from $D$ to $V_1$ for $D_2$, implying that $D_2$ is of definite status on this path (\textbf{R4} of the FCI orientation rules in Appendix~\ref{subsec:additional}). Thus, $D_2 \circarrow V_1$ is not possible, and since $D_2 \leftrightarrow V_{1}$ is already ruled out by Lemma~\ref{lemma:basic property of cpdags and pags}, $D_2 \to V_{1}$ is in~$\g[P]$. By the same reasoning, $D_i \to V_1$ is in~$\g[P]$, for $i = 3, \dots k$.
   It then follows that $\langle D, D_1, \dots ,D_k, X, V_1 \rangle $ is a discriminating path from $D$ to $V_1$ for $X$ in~$\g[P]$, so $X \to V_1$ is in~$\g[P]$ (\textbf{R4} in Appendix~\ref{subsec:additional} and the fact that $X \circcirc V_1$, $X \circarrow V_1$ or invisible $X \rightarrow V_1$ is in $\g[P]$) and $X \to V_1$ is visible. This contradicts the fact $X \circcirc V_1$, $X \circarrow V_1$ or invisible $X \rightarrow V_1$ is in $\g[P]$.
\end{proofof}

\begin{proofof}[Lemma~\ref{lemma:adj amen gen}]
   First suppose that $\g$ is amenable relative to $(\mathbf{X},\mathbf{Y})$, meaning that every proper possibly directed path from $\mathbf{X}$ to $\mathbf{Y}$ in~$\g$ starts with a visible edge out of $\mathbf{X}$. Any visible edge in~$\g$ is visible in all $\DAG$s ($\MAG$s) in $[\g]$, and any proper directed path in a $\DAG$ $(\MAG)$ in $[\g]$ corresponds to a proper possibly directed path in~$\g$. Hence, any proper directed path from $\mathbf{X}$ to $\mathbf{Y}$ in any $\DAG$ $(\MAG)$ in $[\g]$ starts with a visible edge out of $\mathbf{X}$.

   Next, suppose that $\g$ violates \amen{} relative to $(\mathbf{X},\mathbf{Y})$. We will show that this implies that there is no adjustment set relative to $(\mathbf{X},\mathbf{Y})$ in~$\g$. Since every directed edge in a $\CPDAG$ is visible and since the same does not hold true in a $\PAG$, we give separate proofs for $\CPDAG$s and $\PAG$s.

Suppose that $\g$ is a $\PAG$. Since $\g$ violates \amen{} relative to $(\mathbf{X,Y})$, there exists a proper possibly directed path $\pstar$ from some $X \in \mathbf{X}$ to some $Y \in \mathbf{Y}$ in $\g$ that does not start with visible edge out of $X$ in $\g$. Then by Lemma~\ref{lemmadirectedpath} there is a $\MAG$ $\g[M]$ in $[\g]$ such that the path $p$ in~$\g[M]$, consisting of the same sequence of nodes as $\pstar$ in~$\g[P]$, contains a subsequence $p'$ that is a proper directed path from $X$ to $Y$ starting with an invisible edge in~$\g[M]$. Hence, $\g[M]$ violates \amen{} relative to $(\mathbf{X,Y})$. Since the generalized adjustment criterion is complete for $\MAG$s (Theorem~5.8~in~\citealp{vanconstructing}) this means that there is no adjustment set relative to $(\mathbf{X},\mathbf{Y})$ in~$\g[M]$. Hence, there is no adjustment set relative to $(\mathbf{X},\mathbf{Y})$ in~$\g$.

Next, suppose $\g$ is a $\CPDAG$. We now show how to find $\DAG$s $\g[D]_1$ and $\g[D]_2$ in $[\g]$, such that a proper causal path $q'_1$ from $\mathbf{X}$ to $\mathbf{Y}$ in~$\g[D]_1$ corresponds to a proper non-causal path $q'_2$ from $\mathbf{X}$ to $\mathbf{Y}$ in  $\g[D]_2$ that does not contain colliders.
 Since $\g$ is not amenable relative to $(\mathbf{X,Y})$, there is a proper possibly directed path $\pstar[q]$ from a node $X \in \mathbf{X}$ to a node $Y \in \mathbf{Y}$ that starts with a non-directed edge ($\circcirc$).

   Let $\pstar[q'] = \langle X= V_0,V_1,\dots,V_k=Y\rangle, k \ge 1,$ be a shortest subsequence of $\pstar[q]$ such that $\pstar[q']$ is also a proper possibly directed path from $X$ to $Y$ starting with a non-directed edge in~$\g$. Since we chose $\pstar[q']$ using the additional constraint that it must start with a non-directed edge in~$\g$, we cannot use Lemma~\ref{lemma:unshielded} to guarantee that $\pstar[q']$ is of definite status. Hence, we first show that $\pstar[q']$ is a definite status path, by contradiction. Thus, suppose that $\pstar[q']$ is not a definite status path. Then $k \ge 2$. Since the subpath $\pstar[q'](V_1,Y)$ is a definite status path (otherwise, by Lemma~\ref{lemma:unshielded} we can choose a shorter path), this means that $V_1$ is not of definite status on $\pstar[q']$. This implies $X \in \Adj(V_2,\g)$. Moreover, we must have $X \to V_2$, since $X \circcirc V_2$ contradicts the choice of $\pstar[q']$, and $X \leftarrow V_2$ together with the possibly directed path $\pstar[q'](X,V_2)$ contradicts Lemma~\ref{lemmamarlcycle}. Then $X \to V_2$ implies  $V_1 \to V_2$, otherwise $X \rightarrow V_2$ and a possibly directed path $- \pstar[q'](V_2,X)$ are in $\g$, which contradicts Lemma~\ref{lemmamarlcycle}. But then $V_1$ is a definite non-collider on $\pstar[q']$, which contradicts that $V_1$ is not of definite status.

   Hence, $\pstar[q']$ is a proper definite status possibly directed path from $X$ to $Y$. By Lemma~\ref{lemmamarlrx}, there is a $\DAG$ $\g[D]_1$ in $[\g]$ such that there are no additional arrowheads into $X$, as well as a $\DAG$ $\g[D]_2$ in $[\g]$ such that there are no additional arrowheads into $V_1$. Let $q'_1$ in $\g[D]_1$ ($q'_2$ in $\g[D]_2$) be the path corresponding to $\pstar[q']$ in $\g$. Then $q'_1$ is of the form $= X \to V_1 \to \dots \to Y$ and $q'_2 $ is of the form $X \leftarrow V_1 \to \dots \to Y$. Since $q'_1$ is a proper causal path from $X$ to $Y$ and $q'_2$ is a proper non-causal path from $X$ to $Y$, $f(\mathbf{y}\mid do(\mathbf{x}))$ generally differs in $\g[D]_1$ and $\g[D]_2$. Since $\g[D]_1$ and $\g[D]_2$ are both represented by $\g[C]$, this implies that $f(\mathbf{y}\mid do(\mathbf{x}))$ is not identifiable in $\g[C]$.
\end{proofof}

\begin{proofof}[Lemma~\ref{lemmaprovingeqofacconditiona}]
   We first prove~\ref{l:eqa1}$\Rightarrow$\ref{l:eqa2}. Suppose that $\mathbf{Z} \cap \fb{\g} = \emptyset$.
   Since $\fb{\g[D]} \subseteq \fb{\g}$ ($\fb{\g[M]} \subseteq \fb{\g}$) for any $\DAG$ $\g[D]$ ($\MAG$ $\g[M]$) in $[\g]$, it follows directly that $\mathbf{Z}$ satisfies \forb relative to $(\mathbf{X},\mathbf{Y})$ in all $\DAG$s ($\MAG$s) in $[\g]$.

   Next, we prove $\neg$~\ref{l:eqa1}$\Rightarrow\neg$~\ref{l:eqa2}. Suppose that $\g$ is amenable relative to $(\mathbf{X},\mathbf{Y})$, but there is a node $V \in (\mathbf{Z} \cap \fb{\g})$. Then $V \in \PossDe(W,\g)$ for some $W= V_i, 1 \le i \le k$ on a proper possibly directed path $p = \langle X=V_0, V_1,\dots, V_k=Y \rangle, k\ge 1$. Then $q=p(X,W)$ is proper and $r=p(W,Y)$, where $r$ is allowed to be of zero length (if $W=Y$), does not contain a node in $\mathbf{X}$. Moreover, there is a possibly directed path $s$ from $W$ to $V$, where this path is allowed to be of zero length. We will show that this implies that there is a $\DAG$ $\g[D]$ ($\MAG$ $\g[M]$) in $[\g]$ such that $\mathbf{Z} \cap \fb{\g[D]} \neq \emptyset$ $(\mathbf{Z} \cap \fb{\g[M]} \neq \emptyset)$.

   By Lemma~\ref{lemma:unshielded}, there are subsequences $q'$, $r'$ and $s'$ of $q$, $r$ and $s$ that are unshielded possibly directed paths (again $r'$ and $s'$ are allowed to be paths of zero length). Moreover, $q'$ is proper and must start with a directed (visible) edge, since otherwise the concatenated path $q' \oplus r'$, which is a proper possibly directed path from $X$ to $Y$, would violate \amen.
 Lemma~\ref{lemma:unshielded edges} then implies that $q'$ is a directed path from $X$ to $W$ in~$\g$.

   By Lemma~\ref{lemmamarlrx}, there is at least one $\DAG$ $\g[D]$ ($\MAG$ $\g[M]$) in  $[\g]$ that has no additional arrowheads into $W$. In this graph $\g[D]$ ($\g[M]$), the first edge on the path corresponding to $r'$ is oriented out of $W$ and since $r'$ is an unshielded possibly directed path in~$\g[P]$, by Lemma~\ref{lemma:unshielded edges} the path in~$\g[D]$ ($\g[M]$) corresponding to $r'$ is a directed path from $W$ to $Y$. By the same reasoning, the path corresponding to $s'$ in~$\g[D]$ ($\g[M]$) is a directed path from $W$ to $V$.
   Hence,  $V \in \fb{\g[D]}$  ($V \in \fb{\g[M]}$), so that $\mathbf{Z}$ does not satisfy \forb relative to $(\mathbf{X},\mathbf{Y})$ in~$\g[D]$ ($\g[M]$).
\end{proofof}

We now start the path of proving Lemma~\ref{lemmaequivalenceofcondb}. The most involved part is proving the implication $\neg$~\ref{l:eqb2}$\Rightarrow\neg$~\ref{l:eqb1}, that is,
if there is a proper non-causal path $p$ from $\mathbf{X}$ to $\mathbf{Y}$ that is m-connecting given $\mathbf{Z}$ in a $\DAG$ $\g[D]$ ($\MAG$ $\g[M]$) in $[\g]$,
then there must be a proper non-causal definite status path $\pstar[p]$ from $\mathbf{X}$ to $\mathbf{Y}$ that is m-connecting given $\mathbf{Z}$ in~$\g$ as well. We proceed in three steps. First, we show that proper non-causal
paths from $\mathbf{X}$ to $\mathbf{Y}$ that are m-connecting given $\mathbf{Z}$ in~$\g[D]$ $(\g[M])$ correspond to proper non-causal paths in~$\g$ (Lemma~\ref{lemma:non-causal path}). Second, we show that a certain shortest proper non-causal path from $\mathbf{X}$ to $\mathbf{Y}$ that is m-connecting given $\mathbf{Z}$ in~$\g[D]$ $(\g[M])$ corresponds to a proper definite status non-causal path $\pstar$ from $\mathbf{X}$ to $\mathbf{Y}$ in~$\g$ (Lemma~\ref{lemma:1}). Lastly, we show that $\pstar$ is also m-connecting given $\mathbf{Z}$ in~$\g$ (Lemma~\ref{lemma:2}).

\begin{lemma}
Let $\mathbf{X,Y}$ and $\mathbf{Z}$ be pairwise disjoint node sets in a $\PAG$ $(\CPDAG)$ $\g[P]$. Let $\g[P]$ be amenable relative to $(\mathbf{X,Y})$ and let $\mathbf{Z}$ satisfy \forb{} relative to $(\mathbf{X},\mathbf{Y})$ in~$\g[P]$. Let $\g[M]$ be a $\MAG$ $(\DAG)$ in $[\g[P]]$ and let $p = \langle X = V_{0}, V_{1},\dots ,V_{n} = Y \rangle, n \ge 2$, be a proper non-causal path from $\mathbf{X}$ to $\mathbf{Y}$ that is m-connecting given $\mathbf{Z}$ in~$\g[M]$. Let $\pstar$ in~$\g[P]$ denote the path corresponding to $p$ in~$\g[M]$. Then:

\begin{enumerate}[label=(\roman*)]
\item\label{l:ncp:a} Let $i,j \in \mathbb{N}, 0 < i < j \leq n$ such that there is an edge $\langle V_i, V_j \rangle$ in~$\g[P]$. The path $ \pstar(X,V_i) \oplus \langle V_i, V_j \rangle \oplus \pstar(V_j,Y)$ ($\pstar(V_j, Y)$ is possibly of zero length) is a proper non-causal path in~$\g[P]$. For $j=i+1$, this implies that $\pstar$ is a proper non-causal path.
\item\label{l:ncp:b} If $V_1$ is not of definite status on $\pstar$, then $\langle X,V_2 \rangle \oplus \pstar(V_2,Y)$, ($\pstar(V_2, Y)$ is possibly of zero length) exists and is a proper non-causal path in~$\g[P]$.
\item\label{l:ncp:c} If $n \ge 3$ and $V_k, 2 \leq k < n$ is not of definite status on $\pstar$, and every non-endpoint node on $\pstar(X,V_k)$ is a collider on $\pstar$ and a parent of $V_{k+1}$ in~$\g[M]$, then $\langle X,V_{k+1} \rangle \oplus \pstar(V_{k+1},Y)$, ($\pstar(V_{k+1}, Y)$ is possibly of zero length) exists and is a proper non-causal path in~$\g[P]$.
\end{enumerate}
\label{lemma:non-causal path}
\end{lemma}.

\begin{proofof}[Lemma~\ref{lemma:non-causal path}]
All paths considered are proper as they are subsequences of $\pstar$, which consists of the same sequence of nodes as the proper path $p$.

\ref{l:ncp:a} Suppose for a contradiction that $\pstar[q] = \pstar(X,V_i) \oplus \langle V_i, V_j \rangle \oplus \pstar(V_j,Y) $ is possibly directed in~$\g[P]$. All nodes on $\pstar[q]$ except $X$ are in $\fb{\g[P]}$. Since $\g[P]$ is amenable relative to $(\mathbf{X,Y})$, $\pstar[q]$ starts with a visible edge $X \rightarrow V_1$ in~$\g[P]$.
Edge $X \rightarrow V_1$ is then also in~$\g[M]$ and since $p$ is a non-causal path in~$\g[M]$, there is at least one collider on $p$.
Let $V_r, r \ge 1$, be the collider closest to $X$ on $p$. Then $V_r \in \fb{\g[P]}$. Since $p$ is m-connecting given $\mathbf{Z}$, a descendant of $V_r$ is in $\mathbf{Z}$. This contradicts $\mathbf{Z} \cap \fb{\g[P]} = \emptyset$.

\ref{l:ncp:b} Since $V_1$ is not of definite status on $\pstar$ there is an edge between $X$ and $V_2$ in~$\g[P]$, so path $\pstar[q]= \langle X,V_2 \rangle \oplus \pstar(V_2,Y) $ exists in~$\g[P]$.
Suppose for a contradiction that $\pstar[q]$ is a possibly directed path from $X$ to $Y$ in~$\g[P]$. Then $X \to V_2$ is in~$\g[P]$, since $\g[P]$ is amenable relative to $(\mathbf{X,Y})$ and every node on $\pstar[q]$ except $X$ is in $\fb{\g[P]}$.
From~\ref{l:ncp:a} above we know that $\pstar$ is non-causal, so since $\pstar[q]$ is possibly directed, there is an arrowhead towards $X$ on $\pstar(X,V_2)$.

First, suppose $X \bulletbullet V_1 \arrowbullet V_2$. Then $X \rightarrow V_2 \bulletarrow V_1$ implies $X \bulletarrow V_1$ is in~$\g[P]$, since $\g[P]$ is ancestral. This contradicts that $V_1$ is not of definite status on $\pstar$.

Next, suppose $X \arrowbullet V_1 \bulletbullet V_2$ is in~$\g[P]$. If $\g[P]$ is a $\CPDAG$ then $V_1$ is a definite non-collider on $\pstar$, which contradicts that $V_1$ is not of definite status.
If $\g[P]$ is a $\PAG$, then since $X \rightarrow V_{2}$ is a visible edge in~$\g[P]$ there is a node $D \notin \Adj(V_2,\g[P])$ such that $D \bulletarrow X$ in~$\g[P]$ (there is a collider path $D \bulletarrow D_1 \leftrightarrow \dots \leftrightarrow D_s \leftrightarrow X, s \ge 1$, where every node $D_1,\dots , D_s$ is a parent of $V_{2}$ in~$\g[P]$).
   The path $\langle D,X,V_1,V_{2} \rangle$ ($\langle D,D_1, \dots ,D_s,X,V_1,V_{2} \rangle$) is a discriminating path from $D$ to $V_2$ for $V_1$ in~$\g[P]$ so $V_1$ is of definite status on $\pstar$, contrary to the assumption.

\ref{l:ncp:c}  If $\g[P]$ is a $\CPDAG$ and $p(X,V_k), k \ge 2$ is a collider path, then it must be that $k=2$ and $V_k$ is of definite status on $p$. Hence, let $\g[P]$ be a $\PAG$. There is an edge between $X$ and $V_{k+1}$ in~$\g[P]$, otherwise by Lemma~\ref{lemma0} the subpath $\pstar(X,V_{k+1})$ is a discriminating path for $V_k$ in~$\g[P]$, so $V_k$ would be of definite status on $\pstar$. Then path $\pstar[q] = \langle X,V_{k+1} \rangle \oplus \pstar(V_{k+1},Y)$ exists in~$\g[P]$. 
Suppose for a contradiction that $\pstar[q]$ is a possibly directed path from $X$ to $Y$ in~$\g[P]$. Because $\g[P]$ is amenable relative to $(\mathbf{X,Y})$ the edge $X \rightarrow V_{k+1}$ is visible in~$\g[P]$. Also, since $V_1 \rightarrow V_{k+1}$ is in~$\g[M]$, edge $\langle V_1, V_{k+1} \rangle$ is possibly directed towards $V_{k+1}$ in~$\g[P]$.

Consider the edge $X \bulletarrow V_1$ in~$\g[P]$.
If $X \bulletarrow V_1$ is not into $X$ in~$\g[P]$ then $\pstar(X,V_1) \oplus \langle V_1,V_{k+1} \rangle \oplus \pstar(V_{k+1},Y)$ is a proper possibly directed path from $X$ to $Y$ in~$\g[P]$ so $V_1 \in \fb{\g[P]}$. By assumption $V_1$ is a collider on $\pstar$ in~$\g[P]$ and hence also on $p$ in~$\g[M]$. Since $p$ is m-connecting given $\mathbf{Z}$, there is a node $Z \in \mathbf{Z}$ such that $Z \in \De(V_1,\g[M])$. Since $ \De(V_1,\g[M]) \subseteq \PossDe(V_1,\g[P]) \subseteq \fb{\g[P]}$, node $Z \in \fb{\g[P]}$. This contradicts \forb.

So $X \leftrightarrow V_1$ must be in~$\g[P]$.
Since $X \rightarrow V_{k+1}$ is a visible edge in~$\g[P]$ there is a node $D \notin \Adj(V_{k+1},\g[P])$ such that the edge $D \bulletarrow X$ is in~$\g[P]$ (there is a collider path $D \bulletarrow D_1 \leftrightarrow \dots \leftrightarrow D_s \leftrightarrow X, s \ge 1$, and every node $D_1,\dots , D_s$ is a parent of $V_{k+1}$ in~$\g[P]$).
   By Lemma~\ref{lemma0}, path $\langle D,X, V_1, \dots , V_k,V_{k+1} \rangle$ ($\langle D, D_1, \dots ,D_s,X,V_1, \dots , V_k,V_{k+1} \rangle$) is then a discriminating path from $D$ to $V_{k+1}$ for $V_k$ in~$\g[P]$. Hence, $V_k$ is of definite status on $\pstar$, contrary to the original assumption.
\end{proofof}

\begin{lemma}
Let $\mathbf{X,Y}$ and $\mathbf{Z}$ be pairwise disjoint node sets in a $\PAG$ $(\CPDAG)$ $\g[P]$.  Let $\g[P]$ be amenable relative to $(\mathbf{X,Y})$ and let $\mathbf{Z}$ satisfy \forb{} relative to $(\mathbf{X},\mathbf{Y})$ in~$\g[P]$. Let $\g[M]$ be a $\MAG$ $(\DAG)$ in $[\g[P]]$ and let $p$ be a shortest proper non-causal path from $\mathbf{X}$ to $\mathbf{Y}$ that is m-connecting given $\mathbf{Z}$ in~$\g[M]$. Let $\pstar$ in~$\g[P]$ be the corresponding path to $p$ in $\g[M]$. Then $\pstar$ is a proper definite status non-causal path in~$\g[P]$.
   \label{lemma:1}
\end{lemma}

\begin{proofof}[Lemma~\ref{lemma:1}]
Let $p = \langle X = V_{0}, V_{1},\dots ,V_{k} = Y \rangle, k \ge 1$, such that $X \in \mathbf{X}, Y \in \mathbf{Y}$. It follows directly that $\pstar$ is proper and by \ref{l:ncp:a}~in~Lemma~\ref{lemma:non-causal path}, it is also non-causal in~$\g[P]$.

It is left to prove that $\pstar$ is of definite status in~$\g[P]$. For this we rely on the proof of Lemma~1 from \cite{zhang2006causal} (see Lemma~\ref{lemma1} in Appendix~\ref{subsec:additional}) and prove the following claims for $\pstar$.

\begin{claim}
If $V_r, 1 \le r \le k-1$ is not of definite status on $\pstar$, then $V_{r+1}$ is a parent of $V_{r-1}$ in~$\g[M]$.
\label{claim zhang 1}
\end{claim}

\begin{claim}
If $V_r, 1 \le r \le k-1$ is not of definite status on $\pstar$, then $V_{r-1}$ is a parent of $V_{r+1}$ in~$\g[M]$.
 \label{claim zhang 2}
\end{claim}

\noindent{}These claims contradict each other, so every node on $\pstar$ must be of definite status.
\cite{zhang2006causal} proved these claims for a path $\pstar[q]$ in~$\g[P]$, which is the path corresponding to a shortest path $q$ from $\mathbf{X}$ to $\mathbf{Y}$ that is m-connecting given $\mathbf{Z}$ in~$\g[M]$.
These claims are proven using the following argument:

If a node on $\pstar[q]$ is not of definite status, then a subsequence $\pstar[q']$ that is formed by ``jumping'' over one, or a sequence of nodes on $\pstar[q]$ one of which is not of definite status, constitutes a path in~$\g[P]$.
Let $q'$ in~$\g[M]$ be the path corresponding to $\pstar[q']$ in~$\g[P]$. Then $q'$ is a path from $\mathbf{X}$ to $\mathbf{Y}$ that is shorter than $q$, so it is blocked by $\mathbf{Z}$ in~$\g[M]$. By the choice of $q'$, the collider/definite non-collider status of all nodes on $q'$, except two, is the same as on $q$. Therefore, one of these two nodes must block $q'$ in~$\g[M]$. All possible cases for the status of these two nodes are considered in \cite{zhang2006causal} and a contradiction is reached in every case that does not support the claim being proven.

Almost exactly the same argument as in \cite{zhang2006causal} can be carried out to prove that Claim~\ref{claim zhang 1} and~\ref{claim zhang 2} hold for $\pstar$. The only difference is in considering the paths that are subsequences of $p$. Since these paths are shorter than $p$ and proper they are either blocked by $\mathbf{Z}$ or causal in~$\g[M]$. However, the subsequences of $p$ considered in the proof of Lemma~\ref{lemma1} are either immediately non-causal or they are constructed as in \ref{l:ncp:a}-\ref{l:ncp:c} in Lemma~\ref{lemma:non-causal path} and thus non-causal by Lemma~\ref{lemma:non-causal path}. The argument from \cite{zhang2006causal} then still holds for Claims~\ref{claim zhang 1} and~\ref{claim zhang 2} for $\pstar$.
\end{proofof}

\begin{lemma}
  Let $\mathbf{X,Y}$ and $\mathbf{Z}$ be pairwise disjoint node sets in a $\PAG$ $(\CPDAG)$ $\g[P]$.  Let $\g[P]$ be amenable relative to $(\mathbf{X,Y})$ and let $\mathbf{Z}$ satisfy \forb{} relative to $(\mathbf{X},\mathbf{Y})$ in~$\g[P]$. Let $\g[M]$ be a $\MAG$ $(\DAG)$ in $[\g[P]]$ and let $p$ be a path with minimal $\distancefrom{\mathbf{Z}}$ among the shortest proper non-causal paths from $\mathbf{X}$ to $\mathbf{Y}$ that are m-connecting given $\mathbf{Z}$ in~$\g[M]$. Let $\pstar$ in $\g[P]$ be the corresponding path to $p$ in $\g[M]$. Then $\pstar$ is a proper definite status non-causal path from $\mathbf{X}$ to $\mathbf{Y}$ that is m-connecting given $\mathbf{Z}$ in~$\g[P]$.
   \label{lemma:2}
\end{lemma}

 \begin{proofof}[Lemma~\ref{lemma:2}]
By Lemma~\ref{lemma:1}, $\pstar$ is a proper definite status non-causal path in~$\g[P]$. It is only left to prove that $\pstar$ is m-connecting given $\mathbf{Z}$ in~$\g[P]$.

Every definite non-collider on $\pstar$ in~$\g[P]$ corresponds to a non-collider on $p$ in~$\g[M]$, and every collider on $\pstar$ is also a collider on $p$. Since $p$ is m-connecting given $\mathbf{Z}$, no non-collider is in $\mathbf{Z}$ and every collider has a descendant in $\mathbf{Z}$. Let $Q$ be an arbitrary collider (if there is one) on $p$. Then there is a directed path (possibly of zero length) from $Q$ to a node in $\mathbf{Z}$ in~$\g[M]$. Let $d$ be a shortest such path from $Q$ to a node in $\mathbf{Z}$. Let $\pstar[d]$ in~$\g[P]$ denote the corresponding path to $p$ in $\g[M]$. Then $\pstar[d]$ is a possibly directed path from $Q$ to $Z$ in~$\g[P]$.  It is only left to prove that $\pstar[d]$ is a directed path. If $\pstar[d]$ is of zero length, this is trivially true. Otherwise, suppose for a contradiction that there is a circle mark on $\pstar[d]$. Then $\pstar[d]$ must start with a circle mark at $Q$ (Lemma 7.2 from \citet{maathuis2013generalized} and Lemma~\ref{lemma:unshielded edges}, see Appendix~\ref{subsec:additional}).

We first prove that $\pstar[d]$ is unshielded in~$\g[P]$. Suppose for a contradiction that $\pstar[d]$ is shielded. Then there exists a subpath $\langle A,B,C \rangle$ of $\pstar[d]$ such that the edge $\langle A, C \rangle $ is in~$\g[P]$. The path corresponding to $\pstar[d](Q,A) \oplus \langle A, C \rangle \oplus \pstar[d](C,Z)$ must be a non-causal path from $Q$ to $Z$ in~$\g[M]$, otherwise we could have chosen a shorter path $d$. Hence, the edge $A \arrowbullet C$ is in~$\g[M]$. But path $d$ is directed from $Q$ to $Z$ in~$\g[M]$ so $A \rightarrow B \rightarrow C$ is also in~$\g[M]$. This contradicts that $\g[M]$ is ancestral.

 Let $S$ be the first node on $d$ after $Q$.
If $S$ is not a node on $p$, then following the proof of Lemma~2~from~\citet{zhang2006causal} (Lemma~\ref{lemma2}~in~Appendix~\ref{sec:prelim}) there exist nodes $W$ on $p(X,Q)$ and $V$ on $p(Q,Y)$, distinct from $Q$, such that the path $W \bulletarrow S \arrowbullet V$ is in~$\g[M]$ and both $W$ and $V$ have the same colliders/non-collider status of both $p$ and $p' = p(X,W) \oplus \langle W,S,V \rangle \oplus p(V,Y)$. Then $p'$ is m-connecting given $\mathbf{Z}$.
 Since $p'$ is non-causal and shorter than $p$, or as long as $p$ but with a shorter $\distancefrom{\mathbf{Z}}$ than $p$, $p'$ must be non-proper, that is, $S \in \mathbf{X}$.  But then $\langle S,V \rangle \oplus p(V,Y)$ is a proper non-causal m-connecting path from $\mathbf{X}$ to $\mathbf{Y}$ given $\mathbf{Z}$ that is shorter than $p$ in~$\g[M]$. This contradicts our assumption about $p$.

 If $S$ is a node on $p$, then it lies either on $p(X,Q)$ or $p(Q,Y)$. Assume without loss of generality that $S$ is on $p(Q,Y)$. Following the proof of Lemma 2 from \citet{zhang2006causal}, there exists a node $W, W \neq Q$, on $p(X,Q)$ such that $W \bulletarrow S$ is in~$\g[M]$ and $W$ has the same collider/non-collider status on both $p$ and $p' =p(X,W) \oplus  \langle W, S \rangle  \oplus p(S,Y)$. Then $p'$ is m-connecting given $\mathbf{Z}$.
Since $p'$  is proper, and shorter than $p$, or as long as $p$ but with a shorter $\distancefrom{\mathbf{Z}}$ than $p$, $p'$ must be causal in~$\g[M]$. Let $\pstar[p']$ in $\g[P]$ denote the corresponding path to $p$ in $\g[M]$.  Then $\pstar[p']$ is a possibly directed path from $\mathbf{X}$ to $\mathbf{Y}$, $S$ is on $\pstar[p']$ and $Z \in \PossDe(S,\g[P])$, so $Z \in \fb{\g[P]} \cap \mathbf{Z}$. This is a contradiction with $\mathbf{Z} \cap \fb{\g[P]} = \emptyset$.
\end{proofof}

\begin{proofof}[Lemma~\ref{lemmaequivalenceofcondb}]
   Let the $\CPDAG$ $(\PAG)$ $\g$ be amenable relative to $(\mathbf{X},\mathbf{Y})$, and let $\mathbf{Z} \cap \fb{\g}=\emptyset$.

   We first prove $\neg$~\ref{l:eqb1}$\Rightarrow\neg$~\ref{l:eqb3}. Assume $\mathbf{Z}$ does not satisfy \blck relative to $(\mathbf{X},\mathbf{Y})$ in~$\g$.
   Thus, there is a proper definite status non-causal path $p$ from $X \in \mathbf{X}$ to $Y \in \mathbf{Y}$ that is m-connecting given $\mathbf{Z}$ in~$\g$. Consider any $\DAG$ $\g[D]$ ($\MAG$ $\g[M]$) in $[\g]$. The path corresponding to $p$ in~$\g[D]$ ($\g[M]$) is a proper non-causal path from $\mathbf{X}$ to $\mathbf{Y}$ that is m-connecting given $\mathbf{Z}$.
   Hence, $\mathbf{Z}$ does not satisfy \blck relative to $(\mathbf{X,Y})$ in~$\g[D]$ $(\g[M])$ for all $\g[D]$ $(\g[M])$ in $[\g]$.

The implication $\neg$~\ref{l:eqb3}$\Rightarrow\neg$~\ref{l:eqb2} trivially holds, so it is only left to prove that $\neg$~\ref{l:eqb2}$\Rightarrow\neg$~\ref{l:eqb1}.
   Thus, assume that there is a $\DAG$ $\g[D]$ ($\MAG$ $\g[M]$) in $[\g]$ such that there is a proper non-causal path from $\mathbf{X}$ to $\mathbf{Y}$ in~$\g[D]$ ($\g[M]$) that is m-connecting given $\mathbf{Z}$.
   We choose a path $p$ with minimal $\distancefrom{\mathbf{Z}}$ among the shortest proper non-causal paths from $\mathbf{X}$ to $\mathbf{Y}$ that are m-connecting given $\mathbf{Z}$ in~$\g[D]$ $(\g[M])$. By Lemma~\ref{lemma:2}, the corresponding path $\pstar$ in~$\g$ is a proper definite status non-causal path from $\mathbf{X}$ to $\mathbf{Y}$ that is m-connecting given $\mathbf{Z}$.
\end{proofof}

\begin{proofof}[Theorem~\ref{theorem:gac-alt}]
Assume that $\g$ is amenable  relative to $(\mathbf{X,Y})$ and $\mathbf{Z}$ satisfies \blck{} relative to $(\mathbf{X,Y})$ in~$\g$. Then it is left to prove that $\mathbf{Z}$ satisfies \blck relative to $(\mathbf{X,Y})$ in~$\g$ if and only if it satisfies \condtwoprime relative to $(\mathbf{X,Y})$ in~$\gpbd{\mathbf{XY}}$.

We first prove that if $\mathbf{Z}$ satisfies \blck{}, then $\mathbf{Z}$ m-separates $\mathbf{X}$ and $\mathbf{Y}$ in~$\gpbd{\mathbf{XY}}$, that is, $\mathbf{Z}$ blocks all definite status paths from $\mathbf{X}$ to $\mathbf{Y}$ in~$\gpbd{\mathbf{XY}}$. Since every definite status path from $\mathbf{X}$ to $\mathbf{Y}$ has a proper definite status path as a subpath, it is enough to show that $\mathbf{Z}$ blocks all proper definite status paths from $\mathbf{X}$ to $\mathbf{Y}$ in~$\gpbd{\mathbf{XY}}$. Let $p$ be a proper definite status path from $\mathbf{X}$ to $\mathbf{Y}$ in~$\gpbd{\mathbf{XY}}$. Then $p$ must be non-causal, since otherwise $\g$ is not amenable. Let $\pstar$ in~$\g$ be the path corresponding to $p$ in~$\gpbd{\mathbf{XY}}$, consisting of the same sequence of nodes as $p$. Then $\pstar$ is a proper definite status non-causal path from $\mathbf{X}$ to $\mathbf{Y}$ in~$\g$. Thus, $\pstar$ is blocked by $\mathbf{Z}$. Since removing edges cannot m-connect a previously blocked path, $p$ is blocked by $\mathbf{Z}$ in~$\gpbd{\mathbf{XY}}$.

Next, we prove that if $\mathbf{Z}$ m-separates $\mathbf{X}$ and $\mathbf{Y}$ in~$\gpbd{\mathbf{XY}}$, then $\mathbf{Z}$ satisfies \blck{}.
Suppose for a contradiction that there exists a proper definite status non-causal path $\pstar$ from $\mathbf{X}$ to $\mathbf{Y}$ in~$\g$ that is m-connecting given $\mathbf{Z}$. Let $\tilde{p}^*$ in~$\gpbd{\mathbf{XY}}$ be the path corresponding to $\pstar$ in~$\g$, constituted by the same sequence of nodes as $\pstar$, if such a path exists in~$\gpbd{\mathbf{XY}}$.

First, suppose $\tilde{p}^*$ does not exist in~$\gpbd{\mathbf{XY}}$. Then since $\pstar$ is proper, it must start with a visible edge $X \rightarrow D$ in~$\g$ such that $D$ lies on a proper causal path from $X$ to $\mathbf{Y}$, that is, $D \in \fb{\g}$. Since $\pstar$ is non-causal and of definite status it must contain a collider $C \in \PossDe(D,\g)$. Since $\fb{\g}$ is descendral (see Definition~\ref{def:descendral}), $C \in \fb{\g}$ and similarly all descendants of $C$ are in $\fb{\g}$. Considering that $\pstar$ is m-connecting given $\mathbf{Z}$ and $C$ is a collider on $\pstar$, it follows that $C \in \An(\mathbf{Z},\g)$. This contradicts $\mathbf{Z} \cap \fb{\g} = \emptyset$.

Otherwise, $\tilde{p}^*$ is a path in~$\gpbd{\mathbf{XY}}$. Then $\tilde{p}^*$ is a proper definite status non-causal path from $\mathbf{X}$ to $\mathbf{Y}$ in~$\gpbd{\mathbf{XY}}$ that is blocked by $\mathbf{Z}$. Since $\pstar$ is m-connecting given $\mathbf{Z}$ in~$\g$, all colliders on $\pstar$ are in $\An(\mathbf{Z}, \g)$. Since no definite non-collider on $\pstar$ is in $\mathbf{Z}$, no definite non-collider on $\tilde{p}^*$ is in $\mathbf{Z}$. However, $\tilde{p}^*$ is blocked by $\mathbf{Z}$ in~$\gpbd{\mathbf{XY}}$, so at least one collider $C$ on $\tilde{p}^*$ (and therefore $\pstar$ as well) is not in $\An(\mathbf{Z},\gpbd{\mathbf{XY}})$. Thus, any directed path from $C$ to $\mathbf{Z}$ must contain a visible edge $X \rightarrow D$, where $D \in \fb{\g}$. This implies $D \in \An(\mathbf{Z},\g)$, which contradicts $\mathbf{Z} \cap \fb{\g} = \emptyset$.
\end{proofof}

\section{Proofs for Section~\ref{subsec:constructive set gac}} \label{subsec:proofs constructive}

\begin{proofof}[Theorem~\ref{thm:preprocessing}]
The first part of the theorem follows from from Rule 3 of the do-calculus  for DAGs from \cite{Pearl2009}, Rule 3 of the do-calculus for MAGs and PAGs from \cite{zhang2006causal} and the properties of the CPDAG.

Since the generalized adjustment criterion is sound and complete with respect to adjustment in DAGs, CPDAGs, MAGs and PAGs (Theorem~\ref{theorem:gac}), we will prove the second part of this theorem by proving that if $\mathbf{Z}$ satisfies the generalized adjustment criterion relative to $(\mathbf{X,Y})$ in $\g$, then $\mathbf{Z}$ satisfies the generalized adjustment criterion relative to $(\mathbf{X',Y})$ in $\g$. We prove this by showing that $\mathbf{Z}$ satisfies the three conditions of Definition~\ref{def:gac} relative to $(\mathbf{X',Y})$ in $\g$.

We first show that a possibly directed path from $\mathbf{X'}$ to $\mathbf{Y}$ is proper with respect to $\mathbf{X'}$ if and only if it is proper with respect to $\mathbf{X}$. Since $\mathbf{X'} \subseteq \mathbf{X}$ any path that is proper with respect to $\mathbf{X}$ will also be proper with respect to $\mathbf{X'}$. Hence, we only need to show that if $p$ is a possibly directed path from $X' \in \mathbf{X'}$ to $Y \in \mathbf{Y}$ that is proper with respect to $\mathbf{X'}$, then $p$ is also proper with respect to $\mathbf{X}$. Suppose for a contradiction that $p$ is not proper with respect to $\mathbf{X}$. Let $p(X,Y)$ be the subpath of $p$ that is proper with respect to $\mathbf{X}$. Since $X \notin \mathbf{X'}$, $X$ must be in $\mathbf{X} \setminus \mathbf{X'}$, which contradicts the assumption that there are no possibly directed paths from $\mathbf{X} \setminus \mathbf{X'}$ to $\mathbf{Y}$ that are proper with respect to $\mathbf{X}$.

Then $\fb[\mathbf{X',Y}]{\g} \subseteq \fb{\g}$, so $\mathbf{Z}$ must satisfy \forb{} relative to $\mathbf{(X',Y)}$ in $\g$.
Additionally, since $\mathbf{Z}$ satisfies the generalized adjustment criterion relative to $(\mathbf{X,Y})$ in $\g$, $\g$ must be amenable relative to $(\mathbf{X,Y})$. This means that every possibly directed path from $\mathbf{X}$ to $\mathbf{Y}$ that is proper with respect to $\mathbf{X}$ starts with a visible edge out of $\mathbf{X}$. Since $\mathbf{X'} \subseteq \mathbf{X}$ and since possibly directed paths from $\mathbf{X'}$ to $\mathbf{Y}$ that are proper with respect to $\mathbf{X'}$ are proper with respect to $\mathbf{X}$, it follows that \amen{} is satisfied relative to $(\mathbf{X',Y})$ and $\g$.

It is only left to show that $\mathbf{Z}$ satisfies \blck{} relative to $\mathbf{(X',Y)}$ in $\g$. Since $\mathbf{Z}$ satisfies \blck{} relative to $\mathbf{(X,Y)}$ in $\g$ and since $\mathbf{X'} \subseteq \mathbf{X}$, $\mathbf{Z}$ blocks every definite status non-causal path from $\mathbf{X'}$ to $\mathbf{Y}$ that is proper with respect to $\mathbf{X}$. Hence, we only need to show that $\mathbf{Z}$ also blocks every definite status non-causal path from $\mathbf{X'}$ to $\mathbf{Y}$ that is proper with respect to $\mathbf{X'}$, but not proper with respect to $\mathbf{X}$.
Let $p$ be one such path from $X' \in \mathbf{X'}$ to $Y \in \mathbf{Y}$. Since $p$ is proper with respect to $\mathbf{X'}$, but not with respect to $\mathbf{X}$, let $X \in \mathbf{X} \setminus \mathbf{X'}$ be the node on $p$ such that $p(X,Y)$ is proper with respect to $\mathbf{X}$. Since there is no possibly directed path from $\mathbf{X} \setminus \mathbf{X'}$ to $\mathbf{Y}$ that is proper with respect to $\mathbf{X}$, $p(X,Y)$ must be a non-causal path from $X$ to $Y$. Additionally, since $p(X,Y)$ is a subpath of $p$, $p(X,Y)$ is of definite status. Then $p(X,Y)$ is a non-causal definite status path from $\mathbf{X}$ to $\mathbf{Y}$ that is proper with respect to $\mathbf{X}$, so $p(X,Y)$ is blocked by $\mathbf{Z}$. Thus, $p$ is also blocked by $\mathbf{Z}$.
\end{proofof}

\begin{proofof}[Lemma~\ref{lemma:general path proof}]
There is a proper definite status non-causal path from $\mathbf{X}$ to $\mathbf{Y}$ that is m-connecting given $\adjustb{\g} \setminus \mathbf{I} $. Among all such paths consider the ones with minimal length and among those let $p = \langle X, \dots,Y \rangle, X \in \mathbf{X}, Y \in \mathbf{Y}$ be the path with a shortest distance-from-($\mathbf{\mathbf{X} \cup \mathbf{Y}}$) in~$\g$.
By the choice of $p$,~\ref{l:gpp0} holds. It is left to prove that~\ref{l:gpp1}$-$\ref{l:gpp3} also hold for $p$.

\ref{l:gpp1} Since $p$ is m-connecting given $\adjustb{\g} \setminus \mathbf{I}$ and since $\adjustb{\g} \setminus \mathbf{I}$ $=\PossAn(\mathbf{X} \cup \mathbf{Y}, \g) \setminus (\mathbf{X} \cup \mathbf{Y} \cup \mathbf{I})$, any collider on $p$ is in $\PossAn(\mathbf{X} \cup \mathbf{Y}, \g)$. Since $p$ is proper, no collider on $p$ is in $\mathbf{X}$. Additionally, no collider $C$ on $p$ is in $\mathbf{Y} \setminus \mathbf{I}$, otherwise $p(X,C)$ is a non-causal path and we could have chosen a shorter $p$.
It is only left to show that no collider on $p$ is in $\mathbf{I}$. Suppose for a contradiction that a collider on $p$ is in $\mathbf{I}$. Since $\mathbf{I}$ is a descendral set, all (possible) descendants of this collider are also in $\mathbf{I}$. But then $p$ is not m-connecting given $\adjustb{\g} \setminus \mathbf{I}$, a contradiction.

\ref{l:gpp2} Any definite non-collider on $p$ is a possible ancestor of an endpoint node of $p$ or of a collider on $p$. Then it follows from~\ref{l:gpp1} that any definite non-collider on $p$ is in $\PossAn(\mathbf{X} \cup \mathbf{Y}, \g)$. Furthermore, $p$ is m-connecting given $\adjustb{\g} \setminus \mathbf{I}$, so no definite non-collider on $p$ is in $\adjustb{\g} \setminus \mathbf{I}$. Since $\adjustb{\g} \setminus \mathbf{I}= \PossAn(\mathbf{X} \cup \mathbf{Y}, \g) \setminus (\mathbf{X} \cup \mathbf{Y} \cup \mathbf{I})$, it follows that all definite non-colliders on $p$ are in $\PossAn(\mathbf{X} \cup \mathbf{Y}, \g) \cap (\mathbf{X} \cup \mathbf{Y} \cup \mathbf{I})$. Since $p$ is proper, no definite non-collider on $p$ is in $\mathbf{X}$. Additionally, no definite non-collider $C$ on $p$ is in $\mathbf{Y} \setminus \mathbf{I}$, otherwise we could have chosen path $p(X,C)$ instead of $p$. Thus, all definite non-colliders on $p$ are in $\PossAn(\mathbf{X} \cup \mathbf{Y},\g) \cap \mathbf{I} \subseteq \mathbf{I}$.

\ref{l:gpp3} Let $C$ be a collider on $p$. From~\ref{l:gpp1}, it follows that $C \notin \mathbf{X} \cup \mathbf{Y}$ and that there is an unshielded possibly directed path from $C$ to a node $V \in \mathbf{X} \cup \mathbf{Y}$. Suppose that this path starts with an edge of type $C \circcirc Q$ (possibly $Q =V$). We derive a contradiction by constructing a proper definite status non-causal path from $\mathbf{X}$ to $\mathbf{Y}$ that is m-connecting given $\adjustb{\g} \setminus \mathbf{I}$ and shorter than $p$, or of the same length as $p$, but with a shorter distance-from-($\mathbf{X} \cup \mathbf{Y}$).

Let $A$ and $B$ be nodes on $p$ such that $A \bulletarrow C \arrowbullet B$ is a subpath of $p$ (possibly $A=X$, $B=Y$).
Then paths $A \bulletarrow C \circcirc Q$ and $B \bulletarrow C \circcirc Q$ together with Lemma~\ref{lemma:basic property of cpdags and pags} imply that $A \bulletarrow Q \arrowbullet B$ is in~$\g$.

Since all colliders on $p$ are not in $\mathbf{I}$, and all definite non-colliders on $p$ are in $\mathbf{I}$, and since $\mathbf{I}$ is a descendral set, it follows that no collider on $p$ is a possible descendant of a definite non-collider on $p$.
Thus, if $A \neq X$ ($B \neq Y$), then $A \leftrightarrow C$ ($C \leftrightarrow B$) is in~$\g$. Moreover, if $A \leftrightarrow C$ ($C \leftrightarrow B$) is in~$\g$,  then $A \leftrightarrow Q$ $(Q \leftrightarrow B)$ is in~$\g$, otherwise a possibly directed path $\langle A,Q,C \rangle$ ($\langle B,Q,C \rangle$) and $C \leftrightarrow A$ $(C \leftrightarrow B)$ are in $\g$, which contradicts Lemma~\ref{lemmamarlcycle}.
Hence, if $A \neq X$ $(B \neq Y)$, the collider/definite non-collider status of $A$ $(B)$ is the same on $p$ and on $p(X,A) \oplus \langle A,Q\rangle$ $(\langle Q, B \rangle \oplus p(B,Y))$.

Suppose first that $Q \in \mathbf{X} \cup \mathbf{Y}$.
If $Q \in \mathbf{X}$, then $\langle Q,B \rangle \oplus p(B,Y)$ is a proper definite status non-causal path that is m-connecting given $\adjustb{\g} \setminus \mathbf{I}$ in~$\g$ and shorter than $p$.
Otherwise, $Q \in \mathbf{Y}$. If $Q \in \mathbf{Y} \cap \mathbf{I}$, this would imply that $C \in \mathbf{I}$, which contradicts \ref{l:gpp1}. So $Q$ must be in $\mathbf{Y} \setminus \mathbf{I}$. Then $p(X,A) \oplus \langle A,Q \rangle$ is a non-causal path. Hence, we found a proper definite status non-causal path that is m-connecting given $\adjustb{\g} \setminus \mathbf{I}$ in~$\g$ and shorter than $p$, which is a contradiction.

Next, suppose that $Q \notin \mathbf{X} \cup \mathbf{Y}$. Then if $Q$ is not on $p$, $p(X,A) \oplus \langle A,Q,B \rangle \oplus p(B,Y)$ is a proper definite status non-causal path that is m-connecting given $\adjustb{\g} \setminus \mathbf{I}$ in~$\g$ and of the same length as $p$, but with a shorter distance-from-($\mathbf{X} \cup \mathbf{Y}$).
Otherwise, $Q$ is on $p$. Then $Q$ is a collider on $p$, otherwise $Q \in \mathbf{I}$ and $C \in \mathbf{I}$. Suppose first that $Q$ is on $p(C,Y)$. Then  $p(X,A) \oplus \langle A,Q \rangle \oplus p(Q,Y)$ is a non-causal path because $p(Q,Y)$ is into $Q$. Hence, there is a proper definite status non-causal path that is m-connecting given $\adjustb{\g} \setminus \mathbf{I}$ in~$\g$ and shorter than $p$, which is a contradiction. Next, suppose that $Q$ is on $p(X,C)$. Then $p(X,Q) \oplus \langle Q,B \rangle \oplus p(B,Y)$ is a non-causal path since it contains $Q \leftarrow B$. This path is a proper definite status non-causal path that is m-connecting given $\adjustb{\g} \setminus \mathbf{I}$ in~$\g$ and shorter than $p$, which is a contradiction.
\end{proofof}

The following proof is similar to the proof of Lemma~1 from \citet{richardson2003markov} (see Lemma~\ref{lemma:rich} in Appendix~\ref{subsec:additional}). The difference lies in the fact that Lemma~\ref{lemma:rich prime} additionally considers $\CPDAG$s and $\PAG$s $\g$ in which we define m-separation (m-connection) only for paths of definite status, as well as the fact that we require the resulting path $p$ to be proper and non-causal from $\mathbf{X}$ to $\mathbf{Y}$ in~$\g$.

\begin{proofof}[Lemma~\ref{lemma:rich prime}]
Let $p$ be a path in $\g$ satisfying \ref{l:rp0}$-$\ref{l:rp2}.
If there is no collider on $p$, or all colliders on $p$ are in $\An(\mathbf{Z}, \g)$, then $p$ is a proper definite status non-causal path that is m-connecting given $\mathbf{Z}$ in~$\g$ and the lemma holds.

Hence, assume there is at least one collider $C$ on $p$ that is not in $\An(\mathbf{Z},\g)$. By~\ref{l:rp2}, $ C \in \An(\mathbf{X} \cup \mathbf{Y},\g) \setminus \An( \mathbf{Z}, \g) $.
We now construct a path $q$ from $\mathbf{X}$ to $\mathbf{Y}$ in~$\g$ that is m-connecting given $\mathbf{Z}$ and prove it is proper, of definite status and non-causal.

Let $D$ be the node closest to $Y$ on $p$ such that $D \in \An(\mathbf{X},\g) \setminus \An(\mathbf{Z},\g)$ if such a node exists, otherwise let $D=X$.
Let $E$ be the node closest to $D$ on $p(D,Y)$ such that $E \in \An(\mathbf{Y}, \g) \setminus \An(\mathbf{Z},\g)$ if such a node exists, otherwise let $E=Y$.
Since at least one collider on $p$ is in $\An(\mathbf{X} \cup \mathbf{Y}) \setminus \An(\mathbf{Z},\g)$, either $D \neq X$ or $E \neq Y$ must hold. Moreover, if $D=X$, then since there is at least one collider $C$ on $p$ that is in $ \An(\mathbf{X} \cup \mathbf{Y},\g) \setminus \An( \mathbf{Z}, \g) $, it follows that $E \neq D$. However, if $D \neq X$, then $D=E$ is possible.

Since $D \in \An(\mathbf{X},\g) \setminus \An(\mathbf{Z},\g)$, let $v_D$ be a shortest directed path from $D$ to a node in $\mathbf{X}$ (possibly of length zero). Since $E \in \An(\mathbf{Y}, \g)\setminus \An(\mathbf{Z},\g)$, let $v_{E}$ be a shortest directed path from $E$ to a node in $\mathbf{Y}$ (possibly of length zero). Thus, all non-endpoint nodes on $v_D$ and $v_E$ are in $\An(\mathbf{X} \cup \mathbf{Y}, \g) \setminus \An(\mathbf{Z},\g)$. Also, by the choice of $D$ and $E$, no non-endpoint node on $p(D,E)$ is in $\An(\mathbf{X} \cup \mathbf{Y},\g) \setminus \An(\mathbf{Z},\g)$.
Hence, no non-endpoint node on either $v_D$ or $v_E$ is also on $p(D,E)$.

Let $q = (-v_{D}) \oplus p(D,E) \oplus v_{E}$. We prove that $q$ is a proper definite status non-causal path from $\mathbf{X}$ to $\mathbf{Y}$ in~$\g$ that is m-connecting given $\mathbf{Z}$. Path $q$ is of definite status by construction. To prove that $q$ is proper we must show that no non-endpoint node on $v_D$ or $v_E$ is in $\mathbf{X}$. No non-endpoint node on $v_D$ is in $\mathbf{X}$, otherwise we could have chosen a shorter $v_D$. Similarly, no non-endpoint node on $v_E$ is in $\mathbf{X}$, as this would contradict the choice of $D$. It is left to show that $q$ is non-causal from $\mathbf{X}$ to $\mathbf{Y}$ and m-connecting given $\mathbf{Z}$.

We first show that $q$ is m-connecting given $\mathbf{Z}$. By assumption, no definite non-collider on $p$ is in $\mathbf{Z}$. Additionally, if $D$ and $E$ are non-endpoints on $p$, then all nodes on $v_D$ and $v_E$ are in $\An(\mathbf{X} \cup \mathbf{Y}, \g) \setminus \An(\mathbf{Z},\g)$, that is, no node on either $v_D$ or $v_E$ is in $\mathbf{Z}$. Hence, no definite non-collider on $q$ is in $\mathbf{Z}$. For $q$ to be m-connecting given $\mathbf{Z}$ we still have to show that all colliders on $q$ are in $\An(\mathbf{Z},\g)$. Any collider on $q$ is a non-endpoint node on $p(D,E)$. Since no non-endpoint node on $p(D,E)$ is in $\An(\mathbf{X} \cup \mathbf{Y},\g) \setminus \An(\mathbf{Z},\g)$, by choice of $D$ and $E$, and all colliders on $p$ are in $\An(\mathbf{X} \cup \mathbf{Y} \cup \mathbf{Z}, \g)$, by assumption, it follows that any collider on $p(D,E)$ must be in $\An(\mathbf{Z},\g)$.

It is only left to show that $q$ is a non-causal path from $\mathbf{X}$ to $\mathbf{Y}$ in~$\g$.
If $v_{D}$ is not of zero length, this is obviously true.
If $v_D$ is of zero length, then $q = p(X,E) \oplus v_E$.  Since $v_E$ is a directed path from $E$ to $Y$, we need to show that $p(X,E)$ is non-causal from $X$ to $E$.  Suppose for a contradiction that $p(X,E)$ is possibly directed from $X$ to $E$. Then $q$ is a proper possibly directed path from $\mathbf{X}$ to $\mathbf{Y}$, so $E \in \fb{\g} \subseteq \mathbf{I}$. By assumption, there is at least one collider on $p$, and in particular on $p(E,Y)$. Let $C$ be the collider on $p(E,Y)$ closest to $E$. Then $C \in \De(E,\g)$. Since $\De(E,\g) \subseteq \fb{\g}$, it follows that $C \in \fb{\g}$. Since $\fb{\g} \subseteq \mathbf{I}$, this is in contradiction with~\ref{l:rp2}.
\end{proofof}

\section{Proofs for Section~\ref{sec:relation}} \label{subsec:proofs relation}

\begin{proofof}[Lemma~\ref{lemma:gbc bc}]
If there is a path $p$ from $\mathbf{X}$ to $\mathbf{Y}$ for which the listed statements hold, then no set $\mathbf{Z}$ such that $\mathbf{Z} \cap \De(\mathbf{X},\g[D])  = \emptyset$ can block $p$.

Conversely, since $\mathbf{Z} \cap \De(\mathbf{X},\g[D])  = \emptyset$ and $\mathbf{Z}$ satisfies the generalized adjustment criterion relative to $(\mathbf{X,Y})$ in~$\g[D]$, then $\mathbf{Z'} = \adjustb{\g[D]} \setminus \De(\mathbf{X}, \g[D])$ satisfies the generalized adjustment criterion relative to $(\mathbf{X,Y})$ in~$\g[D]$ (Theorem~\ref{theorem:general set proof}).
Suppose there is no back-door set relative to $(\mathbf{X,Y})$ in~$\g[D]$. 
Then $\mathbf{Z'}$ violates condition \ref{d:bc2} of Pearl's back-door criterion relative to at least one $X \in \mathbf{X}$ and $Y \in \mathbf{Y}$ in~$\g[D]$.
Hence, we can choose $p$ to be a shortest back-door path from a node $X \in \mathbf{X}$ to a node $Y \in \mathbf{Y}$ that is d-connecting given $\mathbf{Z'}$ in~$\g[D]$. Then the statement in~\ref{l:gbc-bc1} holds for $p$.
We prove that the statements~\ref{l:gbc-bc2}$-$\ref{l:gbc-bc4} also hold for $p$.

\ref{l:gbc-bc2} Since $\mathbf{Z'}$ satisfies the generalized adjustment criterion and $\mathbf{Z'}$ does not block the back-door path $p$, it follows that $p$ is not proper. Since $p$ is not proper, a subpath of $p$ forms a proper path $q$ from $\mathbf{X}$ to $\mathbf{Y}$. If $q$ is non-causal, then it is blocked by $\mathbf{Z'}$. But in this case $\mathbf{Z'}$ would block $p$ as well. Hence, $q$ is causal.

\ref{l:gbc-bc3} Any non-collider on $p$ is an ancestor of an endpoint or a collider on $p$. Since $p$ is d-connecting given $\mathbf{Z'}$, all colliders on $p$ are in $\An(\mathbf{Z'},\g[D])$. By our choice of $\mathbf{Z'}$, $\An(\mathbf{Z'},\g[D]) \subseteq \An(\mathbf{X} \cup \mathbf{Y}, \g[D])$, so all colliders on $p$ are in $\An(\mathbf{X} \cup \mathbf{Y}, \g[D])$. Since any non-collider on $p$ is an ancestor of an endpoint or a collider on $p$, it follows that all non-colliders on $p$ are also in $\An(\mathbf{X} \cup \mathbf{Y}, \g[D])$. Moreover, since $p$ is d-connecting given $\mathbf{Z'}$, no non-collider on $p$ is in $\mathbf{Z'}$. Thus, since $\mathbf{Z'} = \An(\mathbf{X} \cup \mathbf{Y}, \g[D]) \setminus (\De(\mathbf{X},\g[D]) \cup \mathbf{Y})$, any non-collider on $p$ must be in $\An(\mathbf{X} \cup \mathbf{Y},\g[D]) \cap (\De(\mathbf{X},\g[D]) \cup \mathbf{Y})$. Path $p$ is a shortest back-door path from $\mathbf{X}$ to $\mathbf{Y}$, so no non-collider on $p$ is in $\mathbf{Y}$. Hence, all non-colliders on $p$ are in $\De(\mathbf{X},\g[D])$.

Now, assume that there is a collider $C$ on $p$. This collider is a descendant of $X$ or a non-collider on $p$, so $C \in \De(\mathbf{X}, \g[D])$ as well.
However, since $\De(\mathbf{X},\g[D]) \cap \mathbf{Z} =  \emptyset$ and $\De(C,\g[D]) \subseteq \De(\mathbf{X},\g[D])$ it follows that $p$ is blocked by $\mathbf{Z'}$. This contradicts that $p$ is d-connecting given $\mathbf{Z'}$.

\ref{l:gbc-bc4} From \ref{l:gbc-bc2}, it follows that $Y \in \De(\mathbf{X},\g[D]$). Additionally, we've shown in \ref{l:gbc-bc3} that there is no collider on $p$ and that all non-collider on $p$ are in $\De(\mathbf{X},\g[D])$. Thus, all nodes on $p$ are in $\De(\mathbf{X},\g[D])$.
\end{proofof}

\begin{proofof}[Lemma~\ref{lemma:adjust xy gbc cond 2}]
Let $p$ be a definite status back-door path from a node $X \in \mathbf{X}$ to a node $Y \in \mathbf{Y}$ in~$\g$.
Then there exists a node $X' \in \mathbf{X}$ (possibly $X'=X$) on $p$ such that the subpath $p(X',Y)$ is a proper path from $\mathbf{X}$ to $\mathbf{Y}$ in~$\g$. Since $p(X',Y)$ is a subpath of $p$, $p(X',Y)$ is of definite status.

If $X' \neq X$ and $X'$ is a definite non-collider on $p$, then since $X' \in \mathbf{Z} \cup \mathbf{X} \setminus \{X\}$, $p$ is blocked by $\mathbf{Z} \cup \mathbf{X} \setminus \{X\}$.
Else, $X' \neq X$ and $X'$ is a definite collider on $p$, or $X'=X$ and $p$ is a proper back-door path in~$\g$. Since $\g$ is amenable, all proper back-door paths from $\mathbf{X}$ to $\mathbf{Y}$ are also proper non-causal paths from $\mathbf{X}$ to $\mathbf{Y}$ in~$\g$.
We prove that $\mathbf{Z} \cup \mathbf{X} \setminus \{X\}$ blocks $p$, by proving that it blocks $p(X',Y)$.

Suppose for a contradiction that $p(X',Y)$ is m-connecting given $\mathbf{Z} \cup \mathbf{X} \setminus \{X\}$ in~$\g$.
We show that it is then possible to construct a proper definite status non-causal path from $\mathbf{X}$ to $\mathbf{Y}$ in~$\g$, that is m-connecting given $\mathbf{Z}$, which contradicts that $\mathbf{Z}$ satisfies \blck{} relative to $(\mathbf{X,Y})$ in~$\g$.
Since $p(X',Y)$ is a proper back-door path and since $\g$ is amenable relative to $(\mathbf{X,Y})$, $p(X',Y)$ is a non-causal path from $\mathbf{X}$ to $\mathbf{Y}$ in~$\g$. Then since $p(X',Y)$ also of definite status, it must be blocked by $\mathbf{Z}$ in $\g$.
Since $p(X',Y)$ is m-connecting given $\mathbf{Z} \cup \mathbf{X} \setminus \{X\}$ and blocked by $\mathbf{Z}$, it follows that no definite non-collider on $p(X',Y)$ is in $\mathbf{Z}$ and there is at least one collider on $p(X',Y)$ that is in $\An(\mathbf{X} \setminus \{X\},\g) \setminus \An(\mathbf{Z},\g)$. Let $C$ be the collider closest to $Y$ on $p(X',Y)$ such that $C \in \An(\mathbf{X},\g) \setminus \An(\mathbf{Z},\g)$.
Let $q$ be of the form $C \rightarrow \dots \rightarrow X'',X'' \in \mathbf{X}$ be the shortest causal path from $C$ to $\mathbf{X}$.
Since $p(X',Y)$ is proper, it follows that $C \neq X''$.

Let $D$ be the node closest to $X''$ on $-q$ such that $D$ is also on $p(C,Y)$ (possibly $D=C$) and $r = -q(X'',D) \oplus p(D,Y)$. It is left to show that $r$ is a proper definite status non-causal path from $\mathbf{X}$ to $\mathbf{Y}$ that is m-connecting given $\mathbf{Z}$ in~$\g$. Since $p(C,Y)$ does not contain a node in $\mathbf{X}$, it follows that $-q(X'',D)$ is of non-zero length. So $r$ is a definite status non-causal path. Additionally, since $q$ was chosen as the shortest path from $C$ to $\mathbf{X}$, it follows that $r$ is proper. Lastly, since both $q$ and $p(C,Y)$ are m-connecting given $\mathbf{Z}$ and $D \notin \mathbf{Z}$ and $D$ is a definite non-collider on $r$, it follows that $r$ is m-connecting given $Z$.
\end{proofof}

\begin{corollary} Let $X$ and $Y$ be distinct nodes in a $\DAG$, $\CPDAG$, $\MAG$ or $\PAG$ $\g$ and let $\gout[R]{X}$ be a graph as defined in Definition~\ref{def:rx}. If there exists a generalized back-door set relative to $(X,Y)$ in~$\g$, then $\dsep{\gout[R]{X}} \subseteq \adjust{\g} \setminus \PossDe(X,\g)$.
\label{cor:dsep adjust}
\end{corollary}

\begin{proofof}[Corollary~\ref{cor:dsep adjust}]
Since there exists a generalized back-door set relative to $(X,Y)$ in~$\g$, by Theorem~\ref{theorem:marloes constructive gbc} $\dsep{\gout[R]{X}} \subseteq \PossAn(X \cup Y,\g) \setminus (\PossDe(X,\g) \cup Y)$ is a generalized back-door set relative to $(X,Y)$ in~$\g$.
Additionally, by Corollary~\ref{cor:constr-set-gbc} $\adjust{\g} \setminus \PossDe(X,\g)$ is also generalized back-door set relative to $(X,Y)$ in~$\g$ and $\adjust{\g} \setminus \PossDe(X,\g) = \PossAn(X \cup Y,\g) \setminus (\PossDe(X,\g) \cup Y)$.
\end{proofof}

\begin{proofof}[Lemma~\ref{lemma:no gbc path}]
Let there be a $p$ from $\mathbf{X}$ to $\mathbf{Y}$ in~$\g$ for which~\ref{l:nogbc0}$-$\ref{l:nogbc3} hold. Then $p$ is a proper definite status non-causal path that is m-connecting given $\adjustb{\g} \setminus \PossDe(\mathbf{X},\g)$. Thus, $\adjustb{\g} \setminus \PossDe(\mathbf{X},\g)$ violates \blck{} relative to $(\mathbf{X,Y})$ in~$\g$ and Theorem~\ref{theorem:general set proof} implies that there is no adjustment set $\mathbf{Z}$ relative to $(\mathbf{X,Y})$ in~$\g$ such that $\mathbf{Z} \cap \PossDe(\mathbf{X},\g) = \emptyset$.

Conversely, assume there is no adjustment set $\mathbf{Z}$ relative to $(\mathbf{X,Y})$ in~$\g$ such that $\mathbf{Z} \cap \PossDe(\mathbf{X},\g) = \emptyset$. Since there exists an adjustment set relative to $(\mathbf{X,Y})$ in~$\g$, $\g$ is amenable relative to $(\mathbf{X,Y})$. Then Theorem~\ref{theorem:general set proof} implies that $\adjustb{\g} \setminus \PossDe(\mathbf{X},\g)$ violates \blck{} relative to $(\mathbf{X,Y})$ in~$\g$.
Hence, there is a proper definite status non-causal path from $\mathbf{X}$ to $\mathbf{Y}$ in~$\g$ that is m-connecting given $\adjustb{\g} \setminus \PossDe(\mathbf{X},\g)$.
Then we can use Lemma~\ref{lemma:general path proof} with $\mathbf{I} = \PossDe(\mathbf{X},\g)$ to choose a shortest path $p$ from $X \in \mathbf{X}$ to $Y \in \mathbf{Y}$ for which~\ref{l:gpp0}$-$\ref{l:gpp3} in Lemma~\ref{lemma:general path proof} hold.

We now show that \ref{l:nogbc0}$-$\ref{l:nogbc3} in Lemma~\ref{lemma:no gbc path} hold for $p$.

\ref{l:nogbc0} Follows immediately from \ref{l:gpp0}~in~Lemma~\ref{lemma:general path proof}.

\ref{l:nogbc1} Any definite non-collider on $p$ is in $\PossDe(\mathbf{X},\g)$ (\ref{l:gpp2} in Lemma~\ref{lemma:general path proof}).
Since there is an adjustment set relative to $(\mathbf{X,Y})$ in~$\g$, it follows from Corollary~\ref{cor:adjust xy gac general} that $\adjustb{\g}$ satisfies the generalized adjustment criterion. Thus, $p$ is blocked by $\adjustb{\g}$ and m-connecting given $\adjustb{\g} \setminus\PossDe(\mathbf{X},\g)$. This implies that at least one definite non-collider on $p$ must be in $\adjustb{\g}\cap \PossDe(\mathbf{X},\g)$.

For the remainder of the proof let $V$ be the definite non-collider that is closest to $Y$ on $p$, among all definite non-colliders on $p$ in $ \adjustb{\g} \cap \PossDe(\mathbf{X},\g)$.

\ref{l:nogbc2} By \ref{l:gpp1}~in~Lemma~\ref{lemma:general path proof}, all colliders on $p$ are in $\adjustb{\g} \setminus \PossDe(\mathbf{X},\g)$.
It is left to show that all definite non-colliders on $p(V,Y)$ are in $\fb{\g}$.

All definite non-colliders on $p$ are possible ancestors of an endpoint node of $p$ or a collider on $p$. Hence, all definite non-colliders on $p$ are in $\PossAn(\mathbf{X} \cup \mathbf{Y},\g)$.
By \ref{l:gpp2}~in~Lemma~\ref{lemma:general path proof}, all definite non-colliders are also in $\PossDe(\mathbf{X},\g)$. Hence, all definite non-colliders on $p$ are in $\PossAn(\mathbf{X} \cup \mathbf{Y},\g) \cap \PossDe(\mathbf{X},\g)$. Additionally, by the choice of $V$, no definite non-collider on $p(V,Y)$ is in $\adjustb{\g} \cap \PossDe(\mathbf{X},\g)$. Since $\adjustb{\g} \cap \PossDe(\mathbf{X},\g) =  (\PossAn(\mathbf{X} \cup \mathbf{Y},\g) \cap \PossDe(\mathbf{X},\g)) \setminus (\mathbf{X} \cup \mathbf{Y} \cup \fb{\g})$, it follows that all definite non-colliders on $p(V,Y)$ are in $\mathbf{X} \cup \mathbf{Y} \cup \fb{\g}$.
It is only left to show that no definite non-collider on $p(V,Y)$ is in $\mathbf{X}$ or $\mathbf{Y} \setminus \fb{\g}$.
Since $p$ is proper, no definite non-collider on $p$ is in $\mathbf{X}$. Also, no non-endpoint node $C$ on $p$ is in $\mathbf{Y} \setminus \fb{\g}$, otherwise $p(X,C)$ is a non-causal path and we could have chosen a shorter $p$. Hence, any definite non-collider on $p(V,Y)$ is in $\fb{\g}$.

\ref{l:nogbc3} Let $V_2$ be the node closest to $Y$ on $p$ such that $p(V,V_2)$ is a possibly directed path from $V$ to $V_2$ (possibly of zero length). We will show that $V_2 =V$. Note that $V_2$ is either $Y$, $V$ or a collider on $p$.  Since $V_2 \in \PossDe(V, \g)$ and $\PossDe(V,\g) \subseteq  \PossDe(\mathbf{X},\g)$, by~\ref{l:nogbc2} $V_2$ cannot be a collider on $p$. Hence, $V_2$ is either $Y$ or $V$.

Let $V_1$ be the node closest to $X$ on $p$ such that $-p(V,V_1)$ is an possibly directed path from $V$ to $V_1$ (possibly of zero length).  We will show that $V_1 =X$. Note that $V_1$ is either $X$, $V$ or a collider on $p$. Using the same reasoning as for $V_2$, $V_1$ cannot be a collider on $p$. So $V_1$ is either $X$ or $V$.
As $V$ is a definite non-collider on $p$, either $V_1 \neq V$ or $V_2 \neq V$.

Suppose that $V_2 \neq V$. Then $V_2 = Y$ and $p(V,Y)$ is a possibly directed path from $V$ to $Y$. Since $V \in \PossDe(\mathbf{X},\g)$, let $q$ be a proper possibly directed path from $X' \in \mathbf{X}$ (possibly $X = X'$) to $V$ in~$\g$. Let $W'$ (possibly $W'=V$) be the node closest to $X'$ on $q$ that is also on $p(V,Y)$. Then $q(X',W') \oplus p(W',Y)$ is a proper possibly directed path from $\mathbf{X}$ to $\mathbf{Y}$ in~$\g$, so $W' \in \fb{\g}$. Since $V \in \PossDe(W',\g)$, it follows that $V \in \fb{\g}$. This contradicts~\ref{l:nogbc1}.

Hence, $V_2 = V$. Then $V_1 =X$ and $p$ is of the form $X \dots \leftarrow V \arrowbullet W \dots Y,$ possibly $W =Y$ is in $\g$. Since $p$ is of definite status, $p(X,V)$ must be of the form $X \leftarrow \dots \leftarrow V$. If $W \neq Y$, then if $W$ is a definite non-collider on $p(V,Y)$, \ref{l:nogbc2} implies that $W \in \fb{\g}$. Hence, $V \leftrightarrow W$ is in $\g$, otherwise $V \in \PossDe(W,\g) \subseteq \fb{\g}$, which contradicts \ref{l:nogbc1}. Else, $W$ is a collider on $p(V,Y)$, so $V \leftrightarrow W$ must be on $p$.
\end{proofof}

\begin{proofof}[Corollary~\ref{cor:noforbx}]
\ref{noforbx1}
Let $\mathbf{X'} = \mathbf{X} \cap \fb{\g}$. Since $\mathbf{X'} \neq \emptyset$, there exists a proper possibly directed path
$p = \langle X = V_1,\dots, V_k = Y\rangle, k >1$, from $X \in \mathbf{X}$ to $Y \in \mathbf{Y}$ in $\g$ such that for some $V_i$, $i \in \{2,\dots, k\}$, $\PossDe(V_i,\g) \cap \mathbf{X'} \neq \emptyset$. Let $V_j, j \in \{2,\dots, k\},$ be the node closest to $Y$ on $p$ such that $\PossDe(V_j,\g) \cap \mathbf{X'} \neq \emptyset$. 
Let $q$ be a shortest possibly directed path from $V_j$ to a node $X'$ in $\mathbf{X'}$. Since $p$ was chosen to be proper with respect to $\mathbf{X}$, we know that $V_j \neq X'$.  

By the choice of $V_j$ on $p$, no other node from $p(V_j,Y)$ is on $q$. By Lemma~\ref{lemma:unshielded}, let $\overline{p(V_j,Y)}$ be a subsequence of $p(V_j,Y)$ that forms a possibly directed definite status path from $V_j$ to $Y$.
Now, we can concatenate $-q$ and $\overline{p(V_j,Y)}$ to form the path $r= (-q) \oplus \overline{p(V_j,Y)}$. We will show that $r$ is a proper definite status non-causal path from $\mathbf{X}$ to $\mathbf{Y}$ that does not contain a collider and consists of nodes in $\fb{\g}$. This means that $r$ satisfies condition \ref{nogac2} in Theorem~\ref{theorem:unified-nope}, so there is no set that satisfies the generalized adjustment criterion relative to $(\mathbf{X,Y})$ in $\g$. By completeness of the generalized adjustment criterion (Theorem~\ref{theorem:gac}), there is no adjustment set relative to $(\mathbf{X,Y})$ in $\g$. 

We first show that $r$ is proper with respect to $\mathbf{X}$. Since $p$ is proper with respect to $\mathbf{X}$ and since $V_j \neq X$, $\overline{p(V_j,Y)}$ does not contain a node in $\mathbf{X}$. Additionally, by choice of $q$, $X'$ is the only node from $\mathbf{X}$ on $q$. Hence, $r$ is proper with respect to $\mathbf{X}$.

Since $q$ and $\overline{p(V_j,Y)}$ are possibly directed paths from $V_j$ to $X'$ and from $V_j$ to $Y$, there is no collider on $r$. Additionally, since $V_j \in \fb{\g}$ and since $\fb{\g}$ is a descendral set in $\g$, all nodes on $r$ are in $\fb{\g}$. 

Next, we show that $r$ is of definite status. First, $\overline{p(V_j,Y)}$ is of definite status.
Since $q$ is a shortest possibly directed path from $V_j$ to $X'$, it is of definite status by Lemma~\ref{lemma:unshielded}. Hence, if $V_j = Y$,  then $r = (-q)$ is of definite status. If $V_j \neq Y$, it is left to show that $V_j$ is of definite status on $r$.

We prove this by contradiction. Thus, suppose that $V_j$ is not of definite status on $r$. Let $\langle A,V_j,B \rangle$ be a subpath of $r$, so that $A$ is the node adjacent to $V_j$ on $q$ and $B$ is the node adjacent to $V_j$ on $\overline{p(V_j,Y)}$. Then $A \bulletarrow V_j \circbullet B$, $A \bulletcirc V_j \arrowbullet B$, or $A \bulletcirc V_j \circbullet B$ is in $\g$ and there is an edge $\langle A, B\rangle$ in $\g$. Since $q$ and $\overline{p(V_j,Y)}$ are possibly directed paths from $V_j$ to $X'$ and from $V_j$ to $Y$, $A \bulletarrow V_j \circbullet B$ and $A \bulletcirc V_j \arrowbullet B$ cannot be in $\g$. Hence, $A \bulletcirc V_j \circbullet B$ is a subpath of $r$.

Since $A \bulletcirc V_j \circbullet B$ is in $\g$, neither $A \bulletarrow B$ nor $A \arrowbullet B$ can be in $\g$ (Lemma~\ref{lemma:basic property of cpdags and pags}). This implies that $A \circcirc B$ is in $\g$ and so, $\langle B,A \rangle \oplus q(A,X')$ is a possibly directed path from $B$ to $X'$, which contradicts the choice of $V_{j}$ (since $B$ is closer to $Y$ on $p$). 

It is only left to show that $r$ is a non-causal path from $X'$ to $Y$. We again use a proof by contradiction. Thus, suppose that $r$ is a possibly directed path from $X'$ to $Y$. Then, since $r(X',V_j) = (-q)$ and $q$ are both possibly directed paths in $\g$, $r$ must start with a non-directed edge. Hence, $r$ is a proper possibly directed path from $\mathbf{X}$ to $\mathbf{Y}$ in $\g$ that starts with a non-directed edge, which contradicts that $\g$ is amenable relative to $(\mathbf{X,Y})$.

\ref{noforbx2} Let $\g$ be a $\DAG$ or $\CPDAG$ that is amenable relative to $(\mathbf{X,Y})$ and let $\mathbf{Y} \subseteq \PossDe(\mathbf{X},\g)$. By \ref{noforbx1}, it follows that if $\mathbf{X} \cap \fb{\g} \neq \emptyset$, then there is no adjustment set relative to $(\mathbf{X,Y})$ in $\g$. Hence, we only prove the converse statement.  

If there is no adjustment set relative to $(\mathbf{X,Y})$ in $\g$, then by the soundness of the generalized adjustment criterion (Theorem~\ref{theorem:gac}) there is no set that satisfies the generalized adjustment criterion relative to $(\mathbf{X,Y})$ in $\g$. Since $\g$ is amenable relative to $(\mathbf{X,Y})$, this means that there is a path $p$ from $X \in \mathbf{X}$ to $Y \in \mathbf{Y}$ in $\g$ that satisfies condition \ref{nogac2} in Theorem~\ref{theorem:unified-nope} i.e., $p$ is a proper definite status non-causal path from $\mathbf{X}$ to $\mathbf{Y}$ in~$\g$ such that every collider on $p$ is in $\adjustb{\g}$ and every definite non-collider on $p$ is in $\fb{\g}$.

We first note that, since $\mathbf{Y} \subseteq \PossDe(\mathbf{X},\g)$, $\mathbf{Y} \subseteq \fb{\g}$. We now show, by contradiction, that there is no collider on $p$. Thus, suppose that there is a collider $C$ on $p$. Then $C$ must be either a descendant of a non-collider on $p$ or a descendant of both $X$ and $Y$. Since $C \in \adjustb{\g}$, $C \notin \fb{\g}$. Moreover, every non-collider on $p$ is in $\fb{\g}$ and since $\fb{\g}$ is a descendral set, $C$ cannot be a descendant of a non-collider on $p$. Then $p$ must be of the form $X \rightarrow C \leftarrow Y$.  Since $Y \in \PossDe(\mathbf{X},\g)$ and $C \in \De(Y,\g)$, this contradicts that $C \notin \fb{\g}$. 

Hence, $p$ does not contain a collider. Additionally, $p$ is a non-causal path from $X$ to $Y$. This implies that there is a node $A$ on $p$, $A \neq X$, such that $-p(A,X)$ is a directed path from $A$ to $X$. Since $A$ is either a non-collider on $p$, or $A =Y$, $A \in \fb{\g}$. Hence, $X \in \mathbf{X} \cap \fb{\g}$.  
\end{proofof}

\begin{proofof}[Corollary~\ref{cor:equivalence gac bc}]
Since there is no directed path from one node in $\mathbf{X}$ to another node in $\mathbf{X}$ in $\g[D]$, it is easy to see that a path of the form $ X \leftarrow V \dots Y$, where $X \in \mathbf{X}, Y \in \mathbf{Y}$ and $V \in \De(\mathbf{X},\g[D])$ cannot occur in~$\g[D]$. Then there can be no path $p$ from $\mathbf{X}$ to $\mathbf{Y}$ that satisfies condition~\ref{l:gbc-bc1}~in~Lemma~\ref{lemma:gbc bc} (condition~\ref{l:nogbc3}~in~Lemma~\ref{lemma:no gbc path}) in $\g[D]$. Hence, conditions \ref{nogbc3}~and~\ref{nobc3}~in~Theorem~\ref{theorem:unified-nope} are violated relative to $(\mathbf{X,Y})$ in $\g[D]$. Then by \ref{nogac}~and~\ref{nobc}~in~Theorem~\ref{theorem:unified-nope} it follows that there exists a set that satisfies the generalized adjustment criterion relative to $(\mathbf{X,Y})$ in $\g[D]$ if and only if there exists a back-door set relative to $(\mathbf{X,Y})$ in $\g[D]$.
\end{proofof}

\begin{proofof}[Corollary~\ref{cor:equivalence gac gbc}]
It is enough to prove that if there exists an adjustment set relative to $(\mathbf{X,Y})$ in~$\g$ and the assumptions of the corollary hold, then there is no proper definite status non-causal path $p_1$ that satisfies \ref{l:nogbc0}$-$\ref{l:nogbc3} in Lemma~\ref{lemma:no gbc path}.

If $\g$ contains no possibly directed path $p = \langle V_1, \dots ,V_k \rangle$, with $k \ge 3, \{V_{1},V_{k}\} \subseteq \mathbf{X}$ and $\{V_{2}, \dots ,V_{k-1}\} \cap \mathbf{X} = \emptyset$, it follows that there cannot be a node $V \in \PossDe(\mathbf{X},\g)$ and the path $X \leftarrow \dots \leftarrow V$ in~$\g$. So there cannot be a proper definite status non-causal path $p_1$ that satisfies ~\ref{l:nogbc0}$-$\ref{l:nogbc3} in Lemma~\ref{lemma:no gbc path}.

Next, suppose that $\g$ is a $\DAG$ or $\CPDAG$ and that $\mathbf{Y} \subseteq \PossDe(\mathbf{X},\g)$ and that there exists an adjustment set relative to $(\mathbf{X,Y})$ in~$\g$. By Corollary~\ref{cor:adjust xy gac general}, it follows that $\adjustb{\g}$ then satisfies \blck{} relative to $(\mathbf{X,Y})$ in~$\g$.
Suppose for a contradiction that there is a proper definite status non-causal path $p_1$ that satisfies ~\ref{l:nogbc0}$-$\ref{l:nogbc3} in Lemma~\ref{lemma:no gbc path}. Then $p_1$ is of the form $ X \leftarrow \dots \leftarrow V \leftarrow Y$. By assumption $Y \in \PossDe(\mathbf{X},\g)$, so it follows that $Y \in \fb{\g}$. By the definition of the forbidden set, every other node on $p$ is also in $\fb{\g}$. But then $\adjustb{\g}$ does cannot block $p$. This is in contradiction with $\adjustb{\g}$ satisfying \blck{} relative to $(\mathbf{X,Y})$ in~$\g$.
\end{proofof}

\section{Adjustment Criterion for DAGs} \label{subsec:proofs shpitser}

In this section we provide the soundness and completeness proof for the adjustment criterion from \cite{shpitser2012avalidity,vanconstructing} (see Definition~\ref{def:johannes ac shpitser}). The main result is given in Theorem~\ref{theorem:sound and complete ac}.

This section can be read independently from the rest of the paper. Since we restrict our proof to $\DAG$s, we first define adjustment sets (see Definition~\ref{def:pearl adjustment}) and the adjustment criterion (see Definition~\ref{def:johannes ac shpitser}) in $\DAG$s. Thus, Definition~\ref{def:pearl adjustment} and Definition~\ref{def:johannes ac shpitser} are special cases of Definition~\ref{defadjustment} and Definition~\ref{def:gac} for $\DAG$s.

\begin{definition}{(\textbf{Adjustment set}; \citealp[Chapter~3.3.1]{Pearl2009})}
  Let $\mathbf{X,Y}$ and $\mathbf{Z}$ be pairwise disjoint node sets in a causal $\DAG$ $\g[D]$. Then $\mathbf{Z}$ is an adjustment set relative to $(\mathbf{X,Y})$ in~$\g[D]$ if for any density $f$ consistent with $\g[D]$:
   \begin{center}
      \resizebox{20pc}{!}{
   $
   f(\mathbf{y}\mid do(\mathbf{x}))=
   \begin{cases}
   f(\mathbf{y}\mid \mathbf{x}) & \text{if }\mathbf{Z} = \emptyset,\\
   \int_{\mathbf{z}}f(\mathbf{y}\mid \mathbf{x,z})f(\mathbf{z})d\mathbf{z} & \text{otherwise.}
   \end{cases}
   $
   }
  \end{center}
   \label{def:pearl adjustment}
\end{definition}

\begin{definition} (\textbf{Adjustment criterion}; cf.~\citealp{shpitser2012avalidity}, \citealp{vanconstructing})
Let $\mathbf{X,Y}$ and $\mathbf{Z}$ be pairwise disjoint node sets in a $\DAG$ $\g[D]$. Let $\fb{\g[D]}$ denote the set of all descendants in $\g[D]$ of any $W \notin \mathbf{X}$ which lies on a proper causal path from $\mathbf{X}$ to $\mathbf{Y}$ in $\g[D]$. Then $\mathbf{Z}$ satisfies the adjustment criterion relative to $(\mathbf{X,Y})$ in~$\g[D]$ if the following two conditions hold:
\begin{enumerate}[label = (\ccctext*), leftmargin=0.5cm,align=left]
\item\label{d:acoriginal1} $\mathbf{Z} \cap \fb{\g[D]} = \emptyset$, and
\item\label{d:acoriginal2} all proper non-causal paths from $\mathbf{X}$ to $\mathbf{Y}$ in~$\g[D]$ are blocked by $\mathbf{Z}$.
\end{enumerate}
\label{def:johannes ac shpitser}
\end{definition}

Definition~\ref{def:johannes ac shpitser} was introduced in \citet{vanconstructing} and differs from the definition of the adjustment criterion in \citet{shpitser2012avalidity} in that it uses $\g[D]$ in \forb{}, as opposed to $\g[D]_{\overline{\mathbf{X}}}$, where $\g[D]_{\overline{\mathbf{X}}}$ is the graph obtained by removing all edges into $\mathbf{X}$ from $\g[D]$. These two formulations of the adjustment criterion are equivalent (\citealp[Remark~4.3~in][]{vanconstructing}). We now give the main result in Theorem~\ref{theorem:sound and complete ac}, which follows directly from Theorem~\ref{theorem:completeness ac} and Theorem~\ref{theorem:soundness ac}. To prove Theorem~\ref{theorem:soundness ac} we rely on Lemma~\ref{lemma:zprime satisfies ac} and Lemma~\ref{lemma:shpitser ynyd}, which are given later in this section.

\begin{theorem}
Let $\mathbf{X,Y}$ and $\mathbf{Z}$ be pairwise disjoint node sets in a causal $\DAG$ $\g[D] = (\mathbf{V},\mathbf{E})$. Then $\mathbf{Z}$ satisfies the adjustment criterion (see Definition~\ref{def:johannes ac shpitser}) if and only if $\mathbf{Z}$ is an adjustment set (see Definition~\ref{def:pearl adjustment}).
\label{theorem:sound and complete ac}
\end{theorem}

\begin{theorem} {(\textbf{Completeness of the adjustment criterion for $\DAG$s})}
Let $\mathbf{X,Y}$ and $\mathbf{Z}$ be pairwise disjoint node sets in a causal $\DAG$ $\g[D]$. If $\mathbf{Z}$ does not satisfy the adjustment criterion relative to $(\mathbf{X,Y})$ in $\g[D]$, then there exists a density $f$ consistent with $\g[D]$ such that $f(\mathbf{y}\mid \text{do}(\mathbf{x})) \neq \int_\mathbf{z}  f(\mathbf{y} \mid \mathbf{x},\mathbf{z})f(\mathbf{z})d\mathbf{z}$.
\label{theorem:completeness ac}
\end{theorem}

\begin{proofof}[Theorem~\ref{theorem:completeness ac}]
Suppose that $\mathbf{Z}$ does not satisfy the adjustment criterion relative to $(\mathbf{X,Y})$ in $\g[D] = (\mathbf{V},\mathbf{E})$. It suffices to show that there is a density consistent with $\g[D]$ such that $E[ Y \mid do(\mathbf{X =1}) ] \neq \int_z E[ Y \mid \mathbf{X=1},\mathbf{z} ]f(\mathbf{z})d\mathbf{z}$ for at least one node $Y \in \mathbf{Y}$.

We consider multivariate Gaussian densities with mean vector zero, constructed using linear structural equation models (SEMs) with Gaussian noise. In particular, we let each random variable $A \in \mathbf{V}$ be a linear combination of its parents in $\g[D]$ and a designated Gaussian noise variable $\epsilon_{A}$ with zero mean and a fixed variance. We also assume that the Gaussian noise variables $\{\mathbf{\epsilon}_{A}: A \in \mathbf{V}\}$, are mutually independent. Thus, this model can be parameterized using one
number per node (the residual variance) and one number per edge
(the edge coefficient).

Since $\mathbf{Z}$ does not satisfy the adjustment criterion relative to $\mathbf{(X,Y)}$ in $\g[D]$, $\mathbf{Z}$ violates \forb{} or \blck{}.
\begin{enumerate}[label=\arabic*]
\item\label{I} If $\mathbf{Z}$ violates \forb{}, then there is a proper causal path $\langle X, V_1 , \dots ,$ $V_k=Y \rangle$, $k\ge 1$, from $X \in \mathbf{X}$ to $Y \in \mathbf{Y}$ in $\g[D]$ and a node $Z \in \mathbf{Z}$ such that
\begin{enumerate}[label=(\alph*)]
\item\label{noforb1} $Z =V_i$, for some $i \in \{1,\dots k-1\}$, or
\item\label{noforb2} $Z \in \De(\mathbf{Y},\g[D])$, or
\item\label{noforb3} $Z \in \De(V_i,\g[D]) \setminus \{V_1,\dots,V_{k-1}\}$ for some $i \in \{1,\dots k-1\}$.
\end{enumerate}
\item\label{II} If $\mathbf{Z}$ violates \blck{}, then there exists a proper non-causal path from $X \in \mathbf{X}$ to $Y \in \mathbf{Y}$ that is d-connecting given $\mathbf{Z}$ in $\g[D]$ such that:
\begin{enumerate}[label=(\alph*)]
\item\label{noblck1} the path does not contain any colliders, or
\item\label{noblck2} the path contains at least one collider.
\end{enumerate}
\end{enumerate}
We now discuss these cases systematically.

\begin{enumerate}[label=(\roman*)]

\item\label{case2} Suppose there is a path $p$ from $X \in \mathbf{X}$ to $Y \in \mathbf{Y}$ that satisfies \ref{II}\ref{noblck1} in $\g[D]$. Since $p$ is a proper non-causal path that does not contain colliders, $p$ starts with an edge into $X$, that is, $p$ is of the form $ X \leftarrow \dots  Y$. 

We define our SEM so that all edge coefficients  except for the ones on $p$ are $0$, and all edge coefficients  on $p$ are in $(0,1)$ and small enough so that we can choose the residual variances so that the variance of every random variable in $\mathbf{V}$ is 1. Then the density $f$ generated by this SEM is consistent with the causal $\DAG$ $\g[D]$. Moreover, $f$ is also consistent with the causal $\DAG$ $\g[D']$ that is obtained from $\g[D]$ by removing all edges except for the ones on $p$.  

Since $Y \dsepp \mathbf{X}$ in $\g[D]^{'}_{\overline{\mathbf{X}}}$, we use Rule 3 of the do-calculus (see Equation~\ref{rule3} in Appendix~\ref{subsec:additional}), with $\mathbf{X'} = \emptyset, \mathbf{W'} = \emptyset$, $\mathbf{Z'} = \mathbf{X}$ and $\mathbf{Y'} = \{Y\}$, so that $E[ Y \mid do(\mathbf{X}=\mathbf{1})]= E[Y]= 0$. 

Since $p$ is proper, no node in $\mathbf{X} \setminus \{X\}$ is on $p$. Additionally, since $p$ is d-connecting given $\mathbf{Z}$ and $p$ does not contain colliders, no node in $\mathbf{Z}$ is on $p$. This implies $Y \dsepp \mathbf{Z} \cup \mathbf{X} \setminus \{X\} $ in $\g[D']$. Furthermore, $Y \dsepp \mathbf{Z} \cup \mathbf{X}\given \mathbf{S}$ in $\g[D']$ for any subset $\mathbf{S}$ of the remaining nodes. In particular, we have $Y \dsepp \mathbf{Z} \cup \mathbf{X} \setminus \{X\} \given X$ in $\g[D']$, so that $\int_\mathbf{z} E[ Y \mid \mathbf{X}=\mathbf{1},\mathbf{z}]f(\mathbf{z})d\mathbf{z} = E[ Y \mid X=1]$. By Theorem~\ref{theorem:mardia-condexp}, $ E[ Y \mid X=1]= \Cov(X,Y)$. By Wright's rule (Theorem~\ref{theorem:wright}), $\Cov(X,Y) = a$, where $a$ is the product of all edge coefficients  on $p$. Since $\int_\mathbf{z} E[ Y \mid \mathbf{X}=\mathbf{1},\mathbf{z}]f(\mathbf{z})d\mathbf{z} = a \neq 0$, this case is completed.

\item\label{case1} Suppose no path satisfies \ref{II}\ref{noblck1} in $\g[D]$, but there is a path $p$ from $X \in \mathbf{X}$ to $Y \in \mathbf{Y}$ that satisfies \ref{I}\ref{noforb1} in $\g[D]$.
Let $\mathbf{\tilde{Z}}$ be the set of all nodes in $\mathbf{Z}$ that are on $p$ and let $Z \in \mathbf{\tilde{Z}}$. 

We define our SEM so that all edge coefficients except the ones on $p$ are $0$, and all edge coefficients  on $p$ are in $(0,1)$ and small enough so that we can choose the residual variances such that the variance of every random variable in $\mathbf{V}$ is 1. Then the density $f$ generated by this SEM is consistent with the causal $\DAG$ $\g[D]$, and also with the causal $\DAG$ $\g[D']$ that is obtained from $\g[D]$ by removing all edges except for the ones on $p$. 

Since no node from $\mathbf{X} \setminus \{X\}$ is on $p$, it follows that $Y \dsepp \mathbf{X}\setminus \{X\} $ in $\g[D']$. Furthermore,  $Y \dsepp \mathbf{X}\setminus \{X\} \given X $ in $\g[D]^{'}_{\overline{X}}$.
We use Rule 3 of the do-calculus, with $\mathbf{X'} = \{X\}, \mathbf{W'} = \emptyset$, $\mathbf{Z'} = \mathbf{X}\setminus \{X\}$ and $\mathbf{Y'} = \{Y\}$, so that $E[ Y \mid do(\mathbf{X}=\mathbf{1})]= E[ Y \mid do(X=1)]$. This means that we have reduced a joint intervention to a single intervention. 

Additionally, since $Y \dsepp X $ in $\g[D]^{'}_{\underline{X}}$, we use Rule 2 of the do-calculus, with $\mathbf{X'} = \emptyset, \mathbf{W'} = \emptyset$, $\mathbf{Z'} = X$ and $\mathbf{Y'} = \{Y\}$, so that $E[ Y \mid do(X=1)]= E[ Y \mid X=1]$. By Theorem~\ref{theorem:mardia-condexp}, $E[ Y \mid X=1] = \Cov(X,Y)$. Let $a$ be the product of all edge coefficients  on $p$. By Wright's rule (Theorem~\ref{theorem:wright}), we have that $\Cov(X,Y) = a \neq 0$.

We complete this case by showing that $\int_\mathbf{z} E[ Y \mid \mathbf{X}=\mathbf{1},\mathbf{z}]f(\mathbf{z})d\mathbf{z}  =0$. Since $p$ is proper, no node in $\mathbf{X} \setminus \{X\}$ is on $p$. Additionally, by our choice of $\mathbf{\tilde{Z}}$, no node in $\mathbf{Z} \setminus  \mathbf{\tilde{Z}}$ is on $p$. Thus, $Y \dsepp \mathbf{X} \cup (\mathbf{Z} \setminus \mathbf{\tilde{Z}}) \mid \mathbf{\tilde{Z}}$ in $\g[D']$. Then $\int_\mathbf{z} E[ Y \mid \mathbf{X}=\mathbf{1},\mathbf{z}]f(\mathbf{z})d\mathbf{z} = \int_{\mathbf{\tilde{z}}}E[ Y \mid \mathbf{\tilde{z}}]f(\mathbf{\tilde{z}})d\mathbf{\tilde{z}} = E[Y]=0$. 

\item\label{case3} Suppose no path satisfies \ref{II}\ref{noblck1} or \ref{I}\ref{noforb1}, but there is a path $p$ from $X \in \mathbf{X}$ to $Y \in \mathbf{Y}$ that satisfies \ref{I}\ref{noforb2} in $\g[D]$. Choose $Z \in \mathbf{Z}$ such that the causal path $q$ from $Y$ to $Z$ is a shortest causal path from $Y$ to a node in $\mathbf{Z}$. Then $Y$ is the only node that is on both $p$ and $q$, otherwise there is a cycle in $\g[D]$. Hence, $p \oplus q$ is a causal path from $X$ to $Z$ that contains $Y$ (see Figure~\ref{fig:case3}).

We define our SEM so that all edge coefficients  except the ones on $p \oplus q$ are $0$, and all edge coefficients  on $p \oplus q$ are in $(0,1)$ and small enough so that we can choose the residual variances such that the variance of every random variable in $\mathbf{V}$ is 1. Then the density $f$ generated by this SEM is consistent with the causal $\DAG$ $\g[D]$, as well as with the causal $\DAG$ $\g[D']$ that is obtained from $\g[D]$ by removing all edges except the ones on $p \oplus q$.
\begin{figure}
\centering
\begin{tikzpicture}[>=stealth',shorten >=1pt,auto,node distance=10cm,main node/.style={minimum size=0.4cm,font=\sffamily\Large\bfseries},scale=0.8,transform shape]
\node[main node] (xj) at (0,0) {$X$};
\node[main node] (m) at (2,0) {};
\node[main node] (n1) at (3.5,0) {};
\node[main node] (a) at (2.7,0.4) {$a$};
\node[main node] (y) at (5.5,0) {$Y$};
\node[main node] (n2) at (7.5,0) {};
\node[main node] (b) at (8.2,0.4) {$b$};
\node[main node] (n3) at (9,0) {};
\node[main node] (z) at (11,0) {$Z$};

\draw [->] (xj) edge  (m);
\draw [dotted] (m) edge (n1);
\draw [->] (n1) edge (y);
\draw [->] (y) edge   (n2);
\draw [dotted] (n2) edge (n3);
\draw [->] (n3) edge  (z);
\end{tikzpicture}
\caption{An example of path $p \oplus q$ in $\g[D]$ corresponding to \ref{case3}, where $Z \in \mathbf{Z}$.}
\label{fig:case3}
\end{figure}
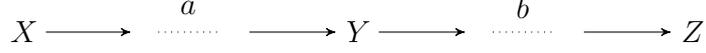

 Since there is no path that satisfies \ref{II}\ref{noblck1} in $\g[D]$, no node from $\mathbf{X}$ is on $q$. Additionally, since $p$ is proper, no node in $\mathbf{X} \setminus \{X\}$ is on $p \oplus q$. Thus, $Y \dsepp \mathbf{X}\setminus \{X\} \given X$ in $\g[D']$. Furthermore, $Y \dsepp \mathbf{X}\setminus \{X\} \given X$ in $\g[D]^{'}_{\overline{X}}$. Hence, we use Rule 3 of the do-calculus, with $\mathbf{X'} = \{X\}, \mathbf{W'} = \emptyset$, $\mathbf{Z'} = \mathbf{X}\setminus \{X\}$ and $ \mathbf{Y'} = \{Y\}$, so that $E[ Y \mid do(\mathbf{X}=\mathbf{1})]= E[ Y \mid do(X=1)]$.

Moreover, $Y \dsepp X $ in $\g[D]^{'}_{\underline{X}}$. Hence, we use Rule 2 of the do-calculus, with $\mathbf{X'} = \emptyset, \mathbf{W'} = \emptyset$, $\mathbf{Z'} = X$ and $ \mathbf{Y'} = \{Y\}$, so that $E[ Y \mid do(X=1)]= E[ Y \mid X=1]$. Lastly, using Theorem~\ref{theorem:mardia-condexp} and Wright's rule (Theorem~\ref{theorem:wright}), we have that $E[ Y \mid X=1] = \Cov(X,Y) = a$, where $a$ is the product of all edge coefficients on $p$ .

Next, we show that $\int_\mathbf{z}E[Y \mid \mathbf{X=1},\mathbf{z}]f(\mathbf{z})d\mathbf{z} \neq a$. Since no path  satisfies \ref{I}\ref{noforb1}, no node from $\mathbf{Z}$ is on $p$. Furthermore, by the choice of $q$, no node from $\mathbf{Z} \setminus \{Z\}$ is on $q$. Hence, $Z$ is the only node from $\mathbf{Z}$ that is on $p \oplus q$. From the above, we also know that $X$ is the only node from $\mathbf{X}$ that is on $p \oplus q$. Hence, $Y \dsepp (\mathbf{X} \cup \mathbf{Z}) \setminus \{X,Z\} \given \{X,Z\}$ in $\g[D']$ and we have that $\int_\mathbf{z}E[Y \mid \mathbf{X=1},\mathbf{z}]f(\mathbf{z})d\mathbf{z} =\int_z E[Y \mid X=1,z]f(z)dz$.

Let $b$ be the product of all edge coefficients on $q$. By Wright's rule (Theorem~\ref{theorem:wright}), we have that $\Cov(X,Y) = a, \Cov(Y,Z) =b$ and $\Cov(X,Z) =ab$. Now, we use Theorem~\ref{theorem:mardia-condexp} to calculate $E[Y \mid X=1,z]$:
\begin{align*}
 E[ Y \mid X=1,z] &= \begin{bmatrix}
       a & b
     \end{bmatrix}
     \begin{bmatrix}
       1 & ab \\
       ab & 1
     \end{bmatrix}^{-1}
     \begin{bmatrix}
       1  \\
       z
     \end{bmatrix}
=\frac{a(1-b^2)}{1-(ab)^2} + \frac{b(1-a^2)}{1-(ab)^2}z.
\end{align*}
\begin{align}
\int_z E[Y \mid X=1,z]f(z)dz
= a\frac{1-b^2}{1-(ab)^2} +  \frac{b(1-a^2)}{1-(ab)^2} E[Z] =a\frac{1-b^2}{1-(ab)^2}. \label{case3eq}
\end{align}
Since $0<a <1$ and $0< b <1$, right-hand side of Equation \eqref{case3eq} is strictly smaller than $a = E[Y \mid do(\mathbf{X} =\mathbf{1})]$.
\item\label{case4}  Suppose no path satisfies \ref{I}\ref{noforb1},  \ref{I}\ref{noforb2}, or \ref{II}\ref{noblck1}, but there is a path $p$ from $X \in \mathbf{X}$ to $Y \in \mathbf{Y}$ that satisfies \ref{I}\ref{noforb3} in $\g[D]$. Let $V_i$, $i \in \{1,\dots,k-1\},$ be a node on $p$ that has a shortest causal path to a node in $\mathbf{Z}$. Let $q_i$ be such a shortest causal path from $V_i$ to $\mathbf{Z}$. Then no node except $V_i$ is on both $p$ and $q_i$, otherwise we would have chosen a different $V_i$ (see Figure~\ref{fig:case4}).

We define our SEM so that all edge coefficients which are not on $p$ or $q_i$ are $0$, and all edge coefficients on $p$ and $q_i$ are in $(0,1)$ and small enough so that we can choose the residual variances so that the variance of every random variable in $\mathbf{V}$ is 1. Then the density $f$ generated by this SEM is consistent with the causal $\DAG$ $\g[D]$, as well as with the causal $\DAG$ $\g[D']$ that is obtained from $\g[D]$ by removing all edges except for the ones on $p$ and $q_i$.
\begin{figure}
\centering
\begin{tikzpicture}[>=stealth',shorten >=1pt,auto,node distance=10cm,main node/.style={minimum size=0.4cm,font=\sffamily\Large\bfseries},scale=0.8,transform shape]
\node[main node] (xj) at (0,0) {$X$};
\node[main node] (m) at (2,0) {};
\node[main node] (n1) at (3.5,0) {};
\node[main node] (a) at (2.7,0.4) {$a$};
\node[main node] (v) at (5.5,0) {$V_i$};
\node[main node] (n2) at (7.5,0) {};
\node[main node] (n3) at (9,0) {};
\node[main node] (y) at (11,0) {$Y$};
\node[main node] (b) at (8.2,0.4) {$b$};
\node[main node] (n4) at (5.5,-2) {};
\node[main node] (n5) at (7,-2) {};
\node[main node] (r) at (6.3,-1.6) {$c$};
\node[main node] (z) at (9,-2) {$Z$};

\draw [->] (xj) edge  (m);
\draw [dotted] (m) edge (n1);
\draw [->] (n1) edge (v);
\draw [->] (v) edge  (n2);
\draw [dotted] (n2) edge (n3);
\draw [->] (n3) edge  (y);
\draw [->] (v) edge (n4);
\draw [dotted] (n4) edge (n5);
\draw [->] (n5) edge  (z);

\end{tikzpicture}
\caption{An example of paths $p$ and $q_i$ in $\g[D]$ corresponding to \ref{case4}, where $Z \in \mathbf{Z}$.}
\label{fig:case4}
\end{figure}
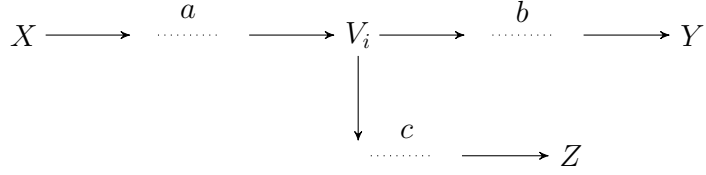

Since there is no path that satisfies \ref{II}\ref{noblck1}, no node from $\mathbf{X}$ is on $q_i$. Additionally, since $p$ is proper, no node from $\mathbf{X} \setminus \{X\}$ is on $p$.
Thus, $Y \dsepp \mathbf{X}\setminus \{X\} \given X$ in $\g[D]^{'}_{\overline{X}}$ and we use Rule 3 of the do-calculus, with $\mathbf{X'} = \{X\}, \mathbf{W'} = \emptyset$, $\mathbf{Z'} = \mathbf{X}\setminus \{X\}$ and $ \mathbf{Y'} = \{Y\}$, so that $E[ Y \mid do(\mathbf{X}=\mathbf{1})]= E[ Y \mid do(X=1)]$.

Additionally, $Y \dsepp X $ in $\g[D]^{'}_{\underline{X}}$. Hence, we use Rule 2 of the do-calculus with $\mathbf{X'} = \emptyset, \mathbf{W'} = \emptyset$, $\mathbf{Z'} = \{X\}$ and $ \mathbf{Y'} = \{Y\}$, so that $E[ Y \mid do(X=1)]= E[ Y \mid X=1]$. Let $a$ be the product of all edge coefficients on $p(X,V_i)$ and let $b$ be the product of all edge coefficients on $p(V_i, Y)$. Then, using Theorem~\ref{theorem:mardia-condexp} and Wright's rule (Theorem~\ref{theorem:wright}), we have that $E[ Y \mid X=1] = \Cov(X,Y) = ab$.

Since there is no path that satisfies \ref{I}\ref{noforb1}, no node from $\mathbf{Z}$ is on $p$. By the choice of $q_i$, no node from $\mathbf{Z} \setminus \{Z\}$ is on $q_i$. Hence, no node from $(\mathbf{X} \cup \mathbf{Z}) \setminus \{X,Z\}$ is on $p$ nor $q_i$. Then  $Y \dsepp (\mathbf{X} \cup \mathbf{Z}) \setminus \{X,Z\} \given \{X,Z\}$ in $\g[D']$ and it follows that $\int_\mathbf{z}E[Y \mid \mathbf{X=1},\mathbf{z}]f(\mathbf{z})d\mathbf{z} =\int_z E[Y \mid X=1,z]f(z)dz$.

Let $c$ be the product of all edge coefficients  on $q_i$.
By Wright's rule (Theorem~\ref{theorem:wright}), we have that $\Cov(X,Y) = ab, \Cov(Y,Z) =bc$ and $\Cov(X,Z) =ac$. We can now use Theorem~\ref{theorem:mardia-condexp} to calculate $E[ Y \mid X=1,z]$:
\begin{align*}
 E[ Y \mid X=1,z] &= \begin{bmatrix}
       ab & bc
     \end{bmatrix}
     \begin{bmatrix}
       1 & ac \\
       ac & 1
     \end{bmatrix}^{-1}
     \begin{bmatrix}
       1  \\
       z
     \end{bmatrix}
     = \frac{ab(1-c^2)}{1-(ac)^2} + \frac{bc(1-a^2)}{1-(ac)^2}z.
\end{align*}
Hence, 
\begin{align}
\int_z E[Y \mid X=1,z]f(z)dz =\frac{ab(1-c^2)}{1-(ac)^2} + \frac{bc(1-a^2)}{1-(ac)^2}E[Z]
=ab \frac{1-c^2}{1-(ac)^2}. \label{case4eq}
\end{align}
Since $0<a<1$, $0<b<1$ and $0<c<1$, right-hand side of Equation \eqref{case4eq} is strictly smaller than $ab = E[Y \mid do(\mathbf{X} =\mathbf{1})]$.
\item\label{case5} {Suppose there is no path that satisfies \ref{I}\ref{noforb1}, \ref{I}\ref{noforb2}, \ref{I}\ref{noforb3}, or \ref{II}\ref{noblck1}, but there is a path that satisfies \ref{II}\ref{noblck2} in $\g[D]$.
Let $p$ be such a path from $X \in \mathbf{X}$ to $Y \in \mathbf{Y}$ in $\g[D]$ that contains the smallest number of colliders among all such paths.
Since no path satisfies \ref{II}\ref{noblck1}, there is at least one collider on $p$. Let $C_1, \dots , C_r$, $r \ge1$, be all colliders on $p$ ordered from the collider closest to $X$ on $p$, which is $C_1$, to the collider closest to $Y$ on $p$, which is $C_r$. Since $p$ is d-connecting given $\mathbf{Z}$, we have $C_i \in \An(\mathbf{Z},\g[D])$ for all $i=1,\dots,r$. .
For each $i=1,\dots,r$, let $q_i$ be a shortest path (possibly of length zero) from $C_i$ to $\mathbf{Z}$. 
{Let $\mathbf{\tilde{Z}}$ be the collection of all nodes in $\mathbf{Z}$ that are enpoints of $ q_1, \dots, q_r$.}
 
We define our SEM so that all edge coefficients which are not on $p,q_1,\dots , q_r$ are $0$, and all edge coefficients which are on $p,q_1,\dots , q_r$ are in $(0,1)$ and are small enough so that we can choose the residual variances such that the variance of every random variable in $\mathbf{V}$ is 1. Then the density $f$ generated by this SEM is consistent with the causal $\DAG$ $\g[D]$, as well as with the causal $\DAG$ $\g[D']$ which is obtained from $\g[D]$ by removing all edges except the ones on $p$, $q_1, \dots , q_r$. {Moreover, $f$ is a non-degenerate multivariate Gaussian density on $\mathbf{V}$. }

Since there is no path that satisfies \ref{II}\ref{noblck1} in $\g[D]$, no node from $\mathbf{X}$ is on $q_r$. Additionally, no node from $\mathbf{X}$ is on $q_i$, for any $i \in \{1,\dots,r-1\},$ when $r >1$, since otherwise there is a proper non-causal path $p'$ from $\mathbf{X}$ to $\mathbf{Y}$ in $\g[D]$ that is d-connecting given $\mathbf{Z}$ and that contains a smaller number of colliders than $p$. Hence, $X$ is the only node from $\mathbf{X}$ that is on $p,q_1,\dots,q_r$.

Since $X$ is the only node from $\mathbf{X}$ that is on $p,q_1,\dots,q_r$, $Y \dsepp \mathbf{X} \setminus \{X\}$ in $\g[D']$. By assumption there is at least one collider on $p$. Since there is no path that satisfies \ref{I}\ref{noforb2} in $\g[D]$, $Y$ is not on $q_1$. Additionally, $Y$ is not on $q_i$, for any $i \in \{2,\dots,r\}$, when $r >1$, otherwise there is a proper non-causal path $p'$ from $\mathbf{X}$ to $\mathbf{Y}$ in $\g[D]$ that is d-connecting given $\mathbf{Z}$ and that contains a smaller number of colliders than $p$. Hence, $Y \dsepp \mathbf{X}$ in $\g[D']$. Furthermore, $Y \dsepp \mathbf{X}$ in $\g[D]^{'}_{\overline{X}}$, so we use Rule 3 of the do-calculus, with $\mathbf{X'} = \emptyset, \mathbf{W'} = \emptyset$, $\mathbf{Z'} = \mathbf{X}$ and $ \mathbf{Y'} = \{Y\}$, so that $E[ Y \mid do(\mathbf{X}=\mathbf{1})]= E[ Y]=0$.

By the choice of $p,q_1,\dots,q_r$, $\mathbf{\tilde{Z}}$ are the only nodes from $\mathbf{Z}$ that are on $p,q_1,\dots, q_r$. Then  $\{X\} \cup \mathbf{\tilde{Z}}$ are the only nodes from $\mathbf{X} \cup \mathbf{Z}$ that are on $p,q_1,\dots, q_r$. Hence, $Y \dsepp (\mathbf{X} \cup \mathbf{Z}) \setminus ( \{X\} \cup \mathbf{\tilde{Z}})$ in $\g[D']$. Furthermore, $Y \dsepp (\mathbf{X} \cup \mathbf{Z}) \setminus ( \{X\} \cup \mathbf{\tilde{Z}}) \given  \{X\} \cup \mathbf{\tilde{Z}}$ in $\g[D']$. Hence, $\int_{\mathbf{z}}E[Y\mid \mathbf{X=1},\mathbf{z}]f(\mathbf{z})d\mathbf{z} = \int_{\mathbf{\tilde{z}}}E[Y\mid X=1,\mathbf{\tilde{z}}]f(\mathbf{\tilde{z}})d\mathbf{\tilde{z}}$.

We now show $\int_{\mathbf{\tilde{z}}}E[Y\mid X=1,\mathbf{\tilde{z}}]f(\mathbf{\tilde{z}})d\mathbf{\tilde{z}} \neq 0$. For this we need the covariance matrix of $(X,\mathbf{\tilde{Z}},Y)^T$, which we will compute by applying Wright's rule to $\g[D']$ (see Figure~\ref{fig:case5}). In order to do this, we first need to show that no node on $q_i$ is on $q_j$, for all $i,j \in \{ 1,\dots, r\}$ with $i \neq j$ and that no node on $q_i$ except $C_i$ is on $p$, for all $i \in \{ 1,\dots, r\}$. From this it will follow that each path $q_i$ ends in a different node in $\mathbf{\tilde{Z}}$. We label these nodes as $\mathbf{\tilde{Z}} =(Z_1,\dots,Z_r)^T$ (see Figure~\ref{fig:case5}). 

\begin{figure}
\centering
\begin{tikzpicture}[>=stealth',shorten >=1pt,auto,node distance=10cm,main node/.style={minimum size=0.4cm,font=\sffamily\Large\bfseries},scale=0.8,transform shape]
\node[main node] (xj) at (0,0) {$X$};
\node[main node] (n1) at (1.5,0) {};
\node[main node] (a) at (1.7,0.4) {$a_1$};
\node[main node] (c1) at (3.5,0) {$C_1$};
\node[main node] (n2) at (5.5,0) {};
\node[main node] (n3) at (7,0) {};
\node[main node] (c2) at (9,0) {$C_r$};
\node[main node] (n4) at (11,0) {};
\node[main node] (a) at (11,0.4) {$a_{r+1}$};
\node[main node] (y) at (12.5,0) {$Y$};

\node[main node] (n5) at (3.5,-2) {};
\node[main node] (n6) at (5,-2) {};
\node[main node] (r) at (4.3,-1.5) {$b_1$};
\node[main node] (z1) at (7,-2) {$Z_1$};

\node[main node] (n7) at (9,-2) {};
\node[main node] (n8) at (10.5,-2) {};
\node[main node] (r) at (9.8,-1.5) {$b_r$};
\node[main node] (z2) at (12,-2) {$Z_r$};

\draw [dotted] (xj) edge (n1);
\draw [->] (n1) edge (c1);
\draw [<-] (c1) edge   (n2);
\draw [dotted] (n2) edge (n3);
\draw [->] (n3) edge  (c2);
\draw [<-] (c2) edge  (n4);
\draw [dotted] (n4) edge (y);

\draw [->] (c1) edge  (n5);
\draw [dotted] (n5) edge (n6);
\draw [->] (n6) edge   (z1);

\draw [->] (c2) edge  (n7);
\draw [dotted] (n7) edge (n8);
\draw [->] (n8) edge  (z2);
\end{tikzpicture}
\caption{An example of paths $p$ and $q_1,\dots,q_r$ in $\g[D]$ corresponding to \ref{case5}.}
\label{fig:case5}
\end{figure}
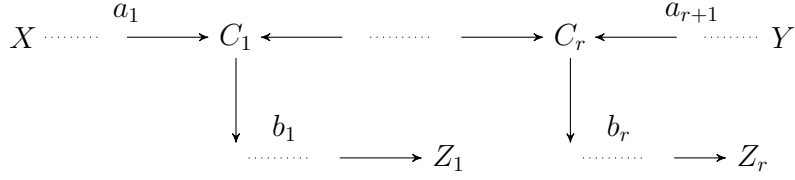

We start by showing that no node on $q_i$, except $C_i$, is on $p$, for any $ i \in \{ 1,\dots, r\}$.
Suppose for a contradiction that for some $i \in \{1,\dots,r\}$ a node on $q_i$, other than $C_i$, is on $p$. Then $q_i$ is at least of length $1$, that is, $C_i \notin \mathbf{Z}$. Let $D$ be the node closest to $C_i$ on $q_i$ that is also on $p$. Then $D$ is either on $p(X,C_i)$ or on $p(C_i,Y)$.

Suppose first that $D$ is on $p(X,C_i)$. Let $p' = p(X,D) \oplus (-q_i)(D,C_i) \oplus p(C_i,Y)$. From the above, we know that $p'$ is a proper path from $\mathbf{X}$ to $\mathbf{Y}$. Since $(-q_i)(D,C_i)$ is of the form $D \leftarrow \dots \leftarrow C_i$, $p'$ is a non-causal path from $X$ to $Y$. By construction $p'$ will have fewer colliders than $p$. Hence, $p'$ is a proper non-causal path from $X$ to $Y$ with fewer colliders than $p$, since $C_i$ is a non-collider on $p'$. Hence, in order to reach a contradiction, we only need to show that $p'$ is d-connecting given $\mathbf{Z}$.

Since $p,q_1,\dots,q_r$ are d-connecting given $\mathbf{Z}$, we only need to discuss the collider/non-collider status of $D$ and $C_i$ on $p'$. Since $C_i \notin \mathbf{Z}$ and since $C_i$ is a non-collider on $p'$, we have that $p'(D,Y) = (-q_i)(D,C_i) \oplus p(C_i,Y )$ is d-connecting given $\mathbf{Z}$. Since $D$ is on $q_i$, $D \in \An(\mathbf{Z},\g[D])$. So if $D$ is a collider on $p'$, $p'$ is d-connecting given $\mathbf{Z}$. If $D$ is a non-collider on $p'$, then $(-p)(D,X)$ is out of $D$, so $D$ is also a non-collider on $p$. Thus, in this case $D \notin \mathbf{Z}$ and hence, $p'$ is d-connecting given $\mathbf{Z}$.

Next, suppose that $D$ is on $p(C_i,Y)$. Let $p' = p(X,C_i) \oplus q_i(C_i,D) \oplus p(D,Y)$. From the above, we know that $p'$ is a proper path from $\mathbf{X}$ to $\mathbf{Y}$. Since $D \in \An(\mathbf{Z},\g[D])$ and since there is no path that satisfies \ref{I}\ref{noforb1}, \ref{I}\ref{noforb2}, or \ref{I}\ref{noforb3}, $p'$ cannot be a causal path. Thus, $p'$ is a proper non-causal path from $\mathbf{X}$ to $\mathbf{Y}$ that by construction has fewer colliders than $p$, since $C_i$ is a non-collider on $p'$. Hence, in order to reach a contradiction, we only need to show that $p'$ is d-connecting given $\mathbf{Z}$.

Since $p,q_1,\dots,q_r$ are d-connecting given $\mathbf{Z}$, we only need to discuss the collider/non-collider status of $D$ and $C_i$ on $p'$.
Since $C_i \notin \mathbf{Z}$ and since $C_i$ is a non-collider on $p'$, we have that $p'(X,D) = p(X,C_i) \oplus q_i(C_i,D)$ is d-connecting given $\mathbf{Z}$. Similarly as above, $D \in \An(\mathbf{Z},\g[D])$, so if $D$ is a collider on $p'$, $p'$ is d-connecting given $\mathbf{Z}$. If $D$ is a non-collider on $p'$, then $p(D,Y)$ starts with an edge out of $D$, so $D$ is also a non-collider on $p$. Thus, in this case $D \notin \mathbf{Z}$ and hence, $p'$ is d-connecting given $\mathbf{Z}$.

Thus, we have shown that no node on $q_i$, other than $C_i$, is on $p$, for all $i \in \{1,\dots,r\}$. Next, we consider the case $r>1$ and show, by contradiction, that no node on $q_i$ is on $q_j$, for any $i,j \in \{ 1,\dots, r\}$ with $i \neq j$.
Hence, suppose that there are $q_i$ and $q_j$ such that a node on $q_i$ is also on $q_j$, for some $i <j$. Let $D$ be the node closest to $C_i$ on $q_i$, that is also on $q_j$. Note that from the above, $D \neq C_j$ and $D \neq C_i$, so that $q_i$ and $q_j$ are at least of length $1$. Then let $p' = p(X,C_i) \oplus q_i(C_i,D) \oplus (-q_j)(D,C_j) \oplus p(C_j,Y)$. As discussed above, no node on $q_i$ (or $q_j$) is in $\mathbf{X}$. Hence, $p'$ is a proper path from $\mathbf{X}$ to $\mathbf{Y}$. Since $ (-q_j)(D,C_j) $ is of the form $D \leftarrow \dots \leftarrow C_j$, $p'$ is also a non-causal path from $X$ to $Y$. Thus $p'$ is a proper non-causal path from $\mathbf{X}$ to $\mathbf{Y}$ that by construction has fewer colliders than $p$, since $C_i$ and $C_j$ are non-colliders on $p'$. Hence, in order to reach a contradiction we only need to show that it is d-connecting given $\mathbf{Z}$.

Since $p,q_1,\dots,q_r$ are d-connecting given $\mathbf{Z}$, we only need to discuss the collider/non-collider status of $D$, $C_i$ and $C_j$ on $p'$. Since $C_i$ and $C_j$ are non-colliders on $p'$, we have that $p'(X,D)$ and $p'(D,Y)$ are both d-connecting given $\mathbf{Z}$. Since $D$ is on $q_i$, $D \in \An(\mathbf{Z},\g[D])$. Since $D$ is a collider on $p'$, $p'$ is d-connecting given $\mathbf{Z}$.

We have now established that $\g[D']$ looks like Figure~\ref{fig:case5}, where none of the paths intersect, and we can compute the covariance  matrix of $(X,\mathbf{\tilde{Z}}^T,Y)^T$ using Wright's rule (Theorem~\ref{theorem:wright}) on $\g[D']$. For this purpose, let $a_1$ and $a_{r+1}$ be the products of all edge coefficients on $p(X,C_1)$ and  $p(C_r,Y)$, respectively. Let $b_i$, $i \in \{1,\dots,r\}$ be product of all edge coefficients on $q_i$. If $r>1$, let $a_{j}$, $j \in \{2,\dots,r\}$ be the product of all edge coefficients  on $p(C_{j-1},C_{j})$.
Let $\mathbf{\Sigma}$ be the covariance matrix of $(X,\mathbf{\tilde{Z}}^T,Y)^T$. Then using Wright's rule (Theorem~\ref{theorem:wright}) on $\g[D']$ yields:
\[ \mathbf{\Sigma} =
\begin{blockarray}{ccccccc}
     	&1	   		&2 				& 3 		& \cdots 		& r+1 			&  r+2  		\\
\begin{block}{c[cccccc]}
      	& 1   		&  a_1b_1  		&  	  		&      			&  	   			&      			\\
      	& a_1b_1	&  1   			& b_1a_2b_2	&      			&    			&  				\\
      	&      		&  b_1a_2b_2 	& 1  		&  \ddots   	&       		& \bigzero  	\\
		&  			&      			& \ddots   	&  \ddots 		&  \ddots   	&  				\\
   		&       	&  \bigzero  	&        	& 1				& b_{r-1}a_rb_r	& 				\\
   		&       	&  			  	&        	& b_{r-1}a_rb_r	&  1        	& b_ra_{r+1}	\\
   		&       	&      			&        	&        		& b_ra_{r+1}	& 1				\\
\end{block}
\end{blockarray}
=
 \begin{bmatrix}
        \mathbf{\Sigma}_{11} & \mathbf{\Sigma}_{12}  	\\
      	\mathbf{\Sigma}_{21} & \mathbf{\Sigma}_{22}
   \end{bmatrix},
 \]
where $\mathbf{\Sigma}_{11}$ is the covariance matrix of $(X,\mathbf{\tilde{Z}}^T)^T$, $\mathbf{\Sigma}_{22} =1$, $\mathbf{\Sigma}_{21}$ $= \mathbf{\Sigma}^{T}_{12}$ and $\mathbf{\Sigma}_{21}= \begin{bmatrix}
0 & \cdots & 0 & b_ra_{r+1} 
\end{bmatrix}$, $\mathbf{\Sigma}_{21} \in \mathbb{R}^{1 \times (r+1)}$.
Then using Theorem~\ref{theorem:mardia-condexp} we have that
\begin{align*}
 E[ Y \mid X=1,\mathbf{\tilde{z}}] &= \begin{bmatrix}
       0 &  \dots & 0 & b_r a_{r+1}
     \end{bmatrix}
      \mathbf{\Sigma}_{11}^{-1}
     \begin{bmatrix}
       1  \\
       \mathbf{\tilde{z}}
     \end{bmatrix}.  
\end{align*}
{ Since our multivariate Gaussian distribution is non-degenerate, $\Det(\mathbf{\Sigma}_{11}) > 0$.}
Let $t_{i,j}$ denote $(i,j)^{th}$ element of $\mathbf{\Sigma}_{11}^{-1} $, $\mathbf{\Sigma}_{11}^{-1} \in \mathbb{R}^{(r+1)\times(r+1)}$.
Putting everything together, we have that
\begin{align}
 \int_{\mathbf{z}}E[Y |\mathbf{X =1}, \mathbf{z}]f(\mathbf{z})d\mathbf{z} &= \int_{\mathbf{\tilde{z}}}E[Y |X=1, \mathbf{\tilde{z}}]f(\mathbf{\tilde{z}})d\mathbf{\tilde{z}} \nonumber \\ 
 &=  \int_{\mathbf{\tilde{z}}}\begin{bmatrix}
       0 & \smash{\dots} & 0 & b_{r}a_{r+1}
     \end{bmatrix}
      \begin{bmatrix}
        t_{1,1}  	&  \smash{\dots}    &  t_{1,r+1}   	\\
      	\vdots   	&  \smash{\ddots}   &  \vdots    	\\
      	t_{r+1,1}  	& \smash{\dots}	&  t_{r+1,r+1}
   \end{bmatrix}
     \begin{bmatrix}
       1  \\
       \mathbf{\tilde{z}}
     \end{bmatrix}f(\mathbf{\tilde{z}})d\mathbf{\tilde{z}}  \nonumber\\
       &= b_{r}a_{r+1}t_{r+1,1} + \sum_{i=1}^r b_{r}a_{r+1}t_{r+1,i+1} E[Z_i]  = b_{r}a_{r+1}t_{r+1,1}. \nonumber
\end{align}
Since $b_r \neq 0$ and $a_{r+1} \neq 0$, it is only left to show that $t_{r+1,1} \neq 0$. Using standard linear algebra, we find
\begin{align*}
t_{r+1,1} &= \frac{(-1)^{r}}{\Det(\mathbf{\Sigma}_{11})} \begin{blockarray}{ccccccc}
     	&	   	& 				&  				 	&          &  		\\
\begin{block}{c|cccccc|}
		&  a_1b_1  		&  1	  	&  b_1a_2b_2    &  	   					&      					&  							\\
		&     			& b_1a_2b_2	&   1  			&   b_2a_3b_3			& 			 			&  \bigzero	 				\\
		&  			 	&   		&   		   	&   					&  						&  							\\
      	&  			 	&   		&  \ddots   	&    \ddots 			&  	\ddots				&  							\\
      	&   			&  \bigzero	&   		   	&   					&  						&  							\\
 		&   			& 		  	&  		 		& b_{r-3}a_{r-2}b_{r-2}	& 1						& b_{r-2}a_{r-1}b_{r-1} 	\\
		&    		 	& 		  	&  		 		&  			  			& b_{r-2}a_{r-1}b_{r-1}	& 1							\\
   		&   			&        	& 				&  		      			& 						& b_{r-1}a_{r}b_r 					\\
\end{block}
\end{blockarray} \\
&= \frac{(-1)^{r}a_{r}b_{r}}{\Det(\mathbf{\Sigma}_{11})}\prod_{i=1}^{r-1}a_{i}b_{i}^2.
\end{align*}
Since $a_i \neq 0$ and $b_i \neq 0$ for all $i \in \{1,\dots,r\}$, the proof is completed.}
\end{enumerate}

\end{proofof}

\begin{theorem} {(\textbf{Soundness of the adjustment criterion for $\DAG$s})}
Let $\mathbf{X,Y}$ and $\mathbf{Z_0}$ be pairwise disjoint node sets in a causal $\DAG$ $\g[D]$. If $\mathbf{Z_0}$ satisfies the adjustment criterion (see Definition~\ref{def:johannes ac shpitser}), then $\mathbf{Z_0}$ is an adjustment set (see Definition~\ref{def:pearl adjustment}).
\label{theorem:soundness ac}
\end{theorem}

\begin{proofof}[Theorem~\ref{theorem:soundness ac}]
Assume that $\mathbf{Z_0}$ satisfies the adjustment criterion (see Definition~\ref{def:johannes ac shpitser}) relative to $(\mathbf{X,Y})$ in~$\g[D]$ and let $f$ be a density consistent with $\g[D]$.
We need to show that

\begin{equation}
f(\mathbf{y} \given \ddo(\mathbf{x})) = \int_{\mathbf{z_0}} f( \mathbf{y} \given \mathbf{x,}\mathbf{z_0}) f(\mathbf{z_0})d\mathbf{z_0}.
\label{eq02}
\end{equation}

Let $\mathbf{Y_D} = \mathbf{Y} \cap \De(\mathbf{X},\g[D])$ and $\mathbf{Y_N} = \mathbf{Y} \setminus \De(\mathbf{X}, \g[D])$.
Then $\mathbf{Y_N} \dsepp \mathbf{X}$ in~$\g[D]_{\overline{\mathbf{X}}}$, since $\g[D]_{\overline{\mathbf{X}}}$ does not contain paths into $\mathbf{X}$ and all paths from $\mathbf{X}$ to $\mathbf{Y_N}$ that are out of $\mathbf{X}$ in $\g[D]_{\overline{\mathbf{X}}}$ must contain a collider by definition of $\mathbf{Y_N}$. 
Hence, using Rule 3 of the do-calculus (see Equation \ref{rule3} in Appendix~\ref{subsec:additional}), with $\mathbf{X'} = \emptyset, \mathbf{W'} = \emptyset$, $\mathbf{Z'} = \mathbf{X}$ and $\mathbf{Y'} = \mathbf{Y_N}$ we have
\begin{align}
\begin{split}
&f(\mathbf{y_N}) = f(\mathbf{y_N} \given \ddo(\mathbf{x})). \label{eq0}
\end{split}
\end{align}

Now, assume that $\mathbf{Y_D} = \emptyset$, so that $\mathbf{Y} = \mathbf{Y_N}$. Then all paths from $\mathbf{X}$ to $\mathbf{Y}$ are non-causal in~$\g[D]$.
Since $\mathbf{Z_0}$ satisfies the adjustment criterion relative to ($\mathbf{X,Y}$) in~$\g[D]$ it blocks all proper non-causal paths from $\mathbf{X}$ to $\mathbf{Y}$. Thus, $\mathbf{Z_0}$ also blocks all non-causal paths from $\mathbf{X}$ to $\mathbf{Y}$ in~$\g[D]$, so $\mathbf{X} \dsepp \mathbf{Y} \given \mathbf{Z_0}$ in~$\g[D]$.
Combining the probabilistic implications of d-separation and Equation \eqref{eq0} with the right-hand side of Equation \eqref{eq02} we obtain
\begin{align}
\int_{\mathbf{z_0}} f(\mathbf{y} \given \mathbf{x,z_0}) f(\mathbf{z_0})d\mathbf{z_0} &= \int_{\mathbf{z_0}} f(\mathbf{y} \given \mathbf{z_0}) f(\mathbf{z_0})d\mathbf{z_0}  \\ \nonumber
&= f(\mathbf{y}) = f(\mathbf{y} \given \ddo(\mathbf{x})). \nonumber
\end{align}

For the remainder of the proof we assume $\mathbf{Y_D} \neq \emptyset$. We enlarge the set $\mathbf{Z_0}$ to $\mathbf{Z} = \mathbf{Z_0} \cup \An(\mathbf{X} \cup \mathbf{Y},\g[D]) \setminus ( \De(\mathbf{X}, \g[D]) \cup \mathbf{Y})$.
Then applying \ref{l:zp-satisf-ac3}~in~Lemma~\ref{lemma:zprime satisfies ac} to the right-hand side of Equation \eqref{eq02} we get
\begin{align}
\int_{\mathbf{z_0}} f( \mathbf{y} \given \mathbf{x,}\mathbf{z_0}) f(\mathbf{z_0})d\mathbf{z_0} = \int_{\mathbf{z}} f( \mathbf{y} \given \mathbf{x,z}) f(\mathbf{z})d\mathbf{z} .
\label{eq03}
\end{align}

Let $\mathbf{Z_D} = \mathbf{Z} \cap \De(\mathbf{X}, \g[D])$ and $\mathbf{Z_N} = \mathbf{Z} \setminus \De(\mathbf{X}, \g[D])$.
Now suppose that $\mathbf{Y} = \mathbf{Y_D}$, so that $\mathbf{Y_N} = \emptyset$. From~\ref{l:s1} in Lemma~\ref{lemma:shpitser ynyd}, we have $\mathbf{Y_D} \dsepp \mathbf{Z_D} \given \mathbf{X} \cup \mathbf{Z_N}$.
Using the probabilistic implications of d-separation, the right-hand side of Equation \eqref{eq03} equals
\begin{align}
\int_{\mathbf{z_D,z_N}} f( \mathbf{y} \given \mathbf{x,z_D,z_N})  f(\mathbf{z_D,z_N})d\mathbf{z_D}d\mathbf{z_N} &= \int_{\mathbf{z_D,z_N}} f( \mathbf{y} \given \mathbf{x,z_N})  f(\mathbf{z_D,z_N})d\mathbf{z_D}d\mathbf{z_N} \nonumber \\
&= \int_{\mathbf{z_N}} f( \mathbf{y} \given \mathbf{x,z_N}) f(\mathbf{z_N})d\mathbf{z_N}. \label{eq3}
\end{align}

Since $\mathbf{Z_N}$ satisfies the generalized back-door criterion relative to $(\mathbf{X,Y})$ in~$\g[D]$ (see~\ref{l:s2} in Lemma~\ref{lemma:shpitser ynyd}) and since the generalized back-door criterion is sound by Theorem 3.1 in \cite{maathuis2013generalized}, we have
\begin{align}
\int_{\mathbf{z_N}} f( \mathbf{y} \given \mathbf{x,z_N}) f(\mathbf{z_N}) d\mathbf{z_N}
& = f(\mathbf{y} \given \ddo(\mathbf{x})). \label{eq33}
\end{align}

\noindent{}Combining Equations \eqref{eq03}, \eqref{eq3} and \eqref{eq33} completes the proof when $\mathbf{Y_N} = \emptyset$.
In the remainder of the proof we assume $\mathbf{Y_D} \neq \emptyset$ and $\mathbf{Y_N} \neq \emptyset$.
From~\ref{l:s1} in Lemma~\ref{lemma:shpitser ynyd}, $\mathbf{Y_D} \dsepp \mathbf{Z_D} \given \mathbf{Y_N} \cup \mathbf{X} \cup \mathbf{Z_N}$. Using the probabilistic implications of d-separation, the right-hand side of Equation \eqref{eq03} equals
\begin{align}
& \int_{\mathbf{z_D,z_N}} f( \mathbf{y_D,y_N} \given \mathbf{x,z_D,z_N})  f(\mathbf{z_D,z_N}) d\mathbf{z_D}d\mathbf{z_N} \nonumber \\
& =\int_{\mathbf{z_D,z_N}} f( \mathbf{y_D} \given \mathbf{y_N,x,z_D,z_N}) f(\mathbf{y_N} \given \mathbf{x,z_D,z_N}) f(\mathbf{z_D,z_N}) d\mathbf{z_D}d\mathbf{z_N} \nonumber \\
&= \int_{\mathbf{z_N}} f( \mathbf{y_D} \given \mathbf{y_N,x,z_N}) \int_{\mathbf{z_D}} f(\mathbf{y_N} \given \mathbf{x,z_D,z_N}) f(\mathbf{z_D,z_N}) d\mathbf{z_D}d\mathbf{z_N}. \label{eq333}
\end{align}

Since $\mathbf{Z}$ satisfies the adjustment criterion relative to $(\mathbf{X,Y})$ in~$\g[D]$ (\ref{l:zp-satisf-ac2}~in~Lemma~\ref{lemma:zprime satisfies ac}), $\mathbf{Z}$ blocks all proper non-causal paths from $\mathbf{X}$ to $\mathbf{Y}$. Hence, $\mathbf{Z}$ blocks all proper paths from $\mathbf{X}$ to $\mathbf{Y_N}$ in~$\g[D]$, so $\mathbf{X} \dsepp \mathbf{Y_N} \given \mathbf{Z}$ in~$\g[D]$. Additionally, the empty set satisfies the generalized back-door criterion relative to ($(\mathbf{Y_N} \cup \mathbf{X} \cup \mathbf{Z_N}),\mathbf{Y_D}$) (\ref{l:s3} in Lemma~\ref{lemma:shpitser ynyd}). 
Since the generalized back-door criterion is sound by Theorem 3.1 in \cite{maathuis2013generalized}, we apply this to the right-hand side of Equation \eqref{eq333}
\begin{align}
  & \int_{\mathbf{z_N}} f( \mathbf{y_D} \given \mathbf{y_N,x,z_N}) \int_{\mathbf{z_D}} f(\mathbf{y_N} \given \mathbf{x,z_D,z_N}) f(\mathbf{z_D,z_N})d\mathbf{z_D}d\mathbf{z_N}  \nonumber  \\
  &= \int_{\mathbf{z_N}} f( \mathbf{y_D} \given \ddo(\mathbf{x,y_N,z_N})) \int_{\mathbf{z_D}}  f(\mathbf{y_N} \given \mathbf{z_N,z_D}) f(\mathbf{z_D,z_N})d\mathbf{z_D}d\mathbf{z_N} \nonumber \\
  &= \int_{\mathbf{z_N}}  f( \mathbf{y_D} \given \ddo(\mathbf{x,y_N,z_N})) f(\mathbf{z_N} \given \mathbf{y_N}) f(\mathbf{y_N})d\mathbf{z_N}. \label{eq4}
\end{align}

To finish the proof we rely on the $\ddo$-calculus rules. By~\ref{l:s4}~in~Lemma~\ref{lemma:shpitser ynyd}, it follows that
$\mathbf{Y_N} \cup\mathbf{Z_N} \dsepp \mathbf{Y_D} \given \mathbf{X}$ in~$\g[D]_{\overline{\mathbf{X}}\underline{\mathbf{Y_N} \cup  \mathbf{Z_N}}}$.
Using Rule 2 of the $\ddo$-calculus (see Equation \ref{rule2} in Appendix~\ref{subsec:additional}) with $\mathbf{X'} = \mathbf{X}, \mathbf{W'} = \emptyset$, $\mathbf{Z'} = \mathbf{Y_N} \cup \mathbf{Z_N}$ and $\mathbf{Y'} = \mathbf{Y_D}$:
\begin{align}
f( \mathbf{y_D} \given \ddo(\mathbf{x,y_N,z_N})) = f( \mathbf{y_D} \given \ddo(\mathbf{x}),\mathbf{z_N,y_N}).
\label{eqdorules1}
\end{align}

\noindent{}Additionally, by~\ref{l:s5}~in~Lemma~\ref{lemma:shpitser ynyd}, $\mathbf{Z_N} \dsepp \mathbf{X} \given \mathbf{Y_N}$ in~$\g[D]_{\overline{\mathbf{X}}}$. Using Rule 3 of the $\ddo$-calculus (see Equation \ref{rule3} in Appendix~\ref{subsec:additional}) with $\mathbf{X'} = \emptyset, \mathbf{W'} = \mathbf{Y_N}$, $\mathbf{Z'} = \mathbf{X}$ and $\mathbf{Y'} = \mathbf{Z_N}$:
\begin{align}
f(\mathbf{z_N} \given \mathbf{y_N}) =  f(\mathbf{z_N} \given  \ddo(\mathbf{x}),\mathbf{y_N}).
\label{eqdorules2}
\end{align}

\noindent{}Finally, we combine Equations \eqref{eqdorules1},\eqref{eqdorules2} and \eqref{eq0} with the right-hand side of Equation \eqref{eq4}:
\begin{align}
& \int_{\mathbf{z_N}}  f( \mathbf{y_D} \given \ddo(\mathbf{y_N,x,z_N})) f(\mathbf{z_N} \given \mathbf{y_N}) f(\mathbf{y_N}) d\mathbf{z_N} \nonumber  \\
&= \int_{\mathbf{z_N}} f( \mathbf{y_D} \given \mathbf{z_N,y_N}, \ddo(\mathbf{x}))  f(\mathbf{z_N} \given \mathbf{y_N}, \ddo(\mathbf{x})) f(\mathbf{y_N} \given \ddo(\mathbf{x})) d\mathbf{z_N} \nonumber \\
&= \int_{\mathbf{z_N}} f( \mathbf{y_D,z_N} \given \mathbf{y_N}, \ddo(\mathbf{x}))   f(\mathbf{y_N} \given \ddo(\mathbf{x})) d\mathbf{z_N} \nonumber \\
&= f( \mathbf{y_D} \given \mathbf{y_N}, \ddo(\mathbf{x}))   f(\mathbf{y_N} \given \ddo(\mathbf{x})) =f(\mathbf{y_D,y_N} \given \ddo(\mathbf{x})) = f(\mathbf{y} \given \ddo(\mathbf{x})).  \nonumber
\end{align}
\end{proofof}
\begin{lemma}
Let $\mathbf{X},\mathbf{Y}$ and $\mathbf{Z_0}$ be pairwise disjoint node sets in a $\DAG$ $\g[D]$ such that $\mathbf{Z_0}$ satisfies the adjustment criterion relative to $(\mathbf{X,Y})$ in $\g[D]$. Let $\mathbf{Z_1} \subseteq \An(\mathbf{X} \cup \mathbf{Y},\g[D]) \setminus ( \De(\mathbf{X}, \g[D]) \cup \mathbf{Y})$ and $\mathbf{Z} = \mathbf{Z_0} \cup \mathbf{Z_1}$.
Then:
\begin{enumerate}[label = (\roman*)]
\item\label{l:zp-satisf-ac1} $\mathbf{X,Y}$ and $\mathbf{Z}$ are pairwise disjoint, and
\item\label{l:zp-satisf-ac2} $\mathbf{Z}$ satisfies the adjustment criterion relative to $(\mathbf{X,Y})$ in~$\g[D]$, and
\item\label{l:zp-satisf-ac3} $\int_\mathbf{z_0} f(\mathbf{y} \mid \mathbf{x,z_0})f(\mathbf{z_0}) d\mathbf{z_0}
= \int_{\mathbf{z_0,z_1}} f(\mathbf{y} \mid \mathbf{x,z_0,z_1}) f(\mathbf{z_0,z_1}) d\mathbf{z_0}d\mathbf{z_1},$ for any density $f$ consistent with $\g[D]$.
\end{enumerate}
\label{lemma:zprime satisfies ac}
\end{lemma}

\begin{proofof}[Lemma~\ref{lemma:zprime satisfies ac}]

\ref{l:zp-satisf-ac1} Since $\mathbf{X,Y}$ and $\mathbf{Z_0}$ are pairwise disjoint, and $\mathbf{Z} = \mathbf{Z_0} \cup \mathbf{Z_1}$, where $\mathbf{Z_1} \cap (\mathbf{X} \cup \mathbf{Y}) = \emptyset$, it follows that $\mathbf{X,Y}$ and $\mathbf{Z}$ are also pairwise disjoint.

\ref{l:zp-satisf-ac2} A set that satisfies \forb{} of the adjustment criterion relative to $(\mathbf{X,Y})$ in~$\g[D]$ satisfies \blck{} if and only if it
$d$-separates $\mathbf{X}$ and $\mathbf{Y}$ in the proper back-door graph $\gpbd[D]{XY}$ (\citet[Theorem~4.6]{vanconstructing}; see Theorem~\ref{theorem:gac-alt} in Section~\ref{sec:main}).
Thus, $\mathbf{Z_0}$ d-separates $\mathbf{X}$ and $\mathbf{Y}$ in~$\gpbd[D]{XY}$. Then by Theorem~\ref{theorem:moralization},  
 any path between $\mathbf{X}$ and $\mathbf{Y}$
in $((\gpbd[D]{XY})_{\An(\mathbf{X} \cup \bY \cup \mathbf{Z_0},\gpbd[D]{XY})})^m$ contains a node in $\mathbf{Z_0}$.

Since $\mathbf{Z_1} \cap \De(\mathbf{X},\g[D]) = \emptyset$, $\fb{\g[D]} \subseteq \De(\mathbf{X},\g[D])$ and $\mathbf{Z_0} \cap \fb{\g[D]} = \emptyset$, $\mathbf{Z} \cap \fb{\g[D]} = \emptyset$.
Additionally, since  $\mathbf{Z} \supseteq \mathbf{Z_0}$, all paths between $\mathbf{X}$ and $\mathbf{Y}$
in $((\gpbd[D]{XY})_{\An(\mathbf{X} \cup \bY \cup \mathbf{Z_0},\gpbd[D]{XY})})^m$ contain a node in $\mathbf{Z}$.
Furthermore, $\mathbf{Z_1} \subseteq \An(\mathbf{X} \cup \mathbf{Y},\g[D])$ implies that $\An(\mathbf{X} \cup \bY \cup \mathbf{Z},\g[D]) = \An(\mathbf{X} \cup \bY \cup \mathbf{Z_0},\g[D])$. Thus, all paths between $\mathbf{X}$ and $\mathbf{Y}$ in $((\gpbd[D]{XY})_{\An(\mathbf{X} \cup \bY \cup \mathbf{Z},\gpbd[D]{XY})})^m$ contain a node in $\mathbf{Z}$. Hence, $\mathbf{Z}$ satisfies \blck{} relative to $(\mathbf{X,Y})$ in~$\g[D]$ (Theorem~\ref{theorem:moralization}, Theorem~\ref{theorem:gac-alt}).

\ref{l:zp-satisf-ac3} We prove this statement by induction on the number of nodes in $\mathbf{Z_1}$. Below, we prove the base case: $|\mathbf{Z_1}| = 1$. We then assume that the result holds for $|\mathbf{Z_1}| = k$, and show that it holds for $|\mathbf{Z_1}| = k+1$. Thus, let $|\mathbf{Z_1}| = k+1$ and take an arbitrary $Z_1 \in \mathbf{Z_1}$. Let $\mathbf{Z'_0} = \mathbf{Z_0} \cup \{Z_1\}$ and $\mathbf{Z'_1} = \mathbf{Z_1} \setminus \{Z_1\}$. The base case then implies
\begin{align}
\int_\mathbf{z_0} f(\mathbf{y} \mid \mathbf{x,z_0})f(\mathbf{z_0}) d\mathbf{z_0}
&= \int_{\mathbf{z_0}} f(\mathbf{y} \mid \mathbf{x,z_0}) \int_{z_1} f(\mathbf{z_0},z_1) dz_1d\mathbf{z_0}  \nonumber \\
&= \int_{\mathbf{z_0},z_1} f(\mathbf{y} \mid \mathbf{x,z_0},z_1) f(\mathbf{z_0},z_1) d\mathbf{z_0}dz_1 \nonumber \\
&= \int_{\mathbf{z'_0}} f(\mathbf{y} \mid \mathbf{x,z'_0}) f(\mathbf{z'_0}) d\mathbf{z'_0}. \label{l1:a}
\end{align}
By~\ref{l:zp-satisf-ac2} above $\mathbf{Z'_0}$ satisfies the adjustment criterion relative to $(\mathbf{X,Y})$ in~$\g[D]$ and $ \mathbf{Z'_1} \subseteq \An(\mathbf{X} \cup \mathbf{Y}, \g[D])$. Then since $|\mathbf{Z'_1}| = k$, by the induction hypothesis
\begin{align}
\int_\mathbf{z'_0} f(\mathbf{y} \mid \mathbf{x,z'_0})f(\mathbf{z'_0}) d\mathbf{z'_0}
= \int_{\mathbf{z'_0,z'_1}} f(\mathbf{y} \mid \mathbf{x,z'_0,z'_1}) f(\mathbf{z'_0,z'_1}) d\mathbf{z'_0}d\mathbf{z'_1}. \label{l1:b}
\end{align}
Combining \ref{l1:a} and \ref{l1:b} yields
\begin{align}
\int_\mathbf{z_0} f(\mathbf{y} \mid \mathbf{x,z_0})f(\mathbf{z_0}) d\mathbf{z_0}
= \int_{\mathbf{z_0,z_1}} f(\mathbf{y} \mid \mathbf{x,z_0,z_1}) f(\mathbf{z_0,z_1}) d\mathbf{z_0}d\mathbf{z_1}. \nonumber
\end{align}

It is left to prove the base case of the induction. Hence, suppose $\mathbf{Z_1} = \{Z_1 \}$.
We show below that either (a) $\mathbf{Y} \dsepp Z_1 \given \mathbf{X} \cup \mathbf{Z_0}$ or (b) $\mathbf{X} \dsepp Z_1 \given \mathbf{Z_0}$ are satisfied in~$\g[D]$. 
(Note that (a) and (b) are very similar to the conditions U1 and U2, as well as U1* and U2*, from \citealp{greenland1999causal}. These conditions are also used in Theorem 5 from \citealp{kuroki2003covariate}, Lemma 3 from \citealp{kuroki04selection} and are the foundation for the results in \citealp{DeLunaEtAl11}. Additionally, \citealp{pearl2014confounding} give a discussion of these conditions, which they refer to as c-equivalence conditions, in their Theorem 1.)

If (a) $\mathbf{Y} \dsepp Z_1 \given \mathbf{X} \cup \mathbf{Z_0}$ in~$\g[D]$, then by the probabilistic implications of d-separation we have that for any density $f$ consistent with $\g[D]$
\begin{align}
\int_\mathbf{z_0} f(\mathbf{y} \mid \mathbf{x,z_0})f(\mathbf{z_0})d\mathbf{z_0} &= \int_\mathbf{z_0} f(\mathbf{y} \mid \mathbf{x,z_0})\int_{z_1}f(\mathbf{z_0},z_1)dz_1d\mathbf{z_0} \nonumber \\
&= \int_{\mathbf{z_0},z_1} f(\mathbf{y} \mid \mathbf{x,z_0}) f(\mathbf{z_0},z_1)d\mathbf{z_0}dz_1 \nonumber \\
&= \int_{\mathbf{z_0},z_1} f(\mathbf{y} \mid \mathbf{x,z_0},z_1)f(\mathbf{z_0},z_1)d\mathbf{z_0}dz_1. \nonumber 
\end{align}

If (b) $\mathbf{X} \dsepp Z_1 \given \mathbf{Z_0}$ in~$\g[D]$, then similarly 
\begin{align}
\int_\mathbf{z_0} f(\mathbf{y} \mid \mathbf{x,z_0})f(\mathbf{z_0})d\mathbf{z_0} &= \int_{\mathbf{z_0}}f(\mathbf{z_0})\int_{z_1} f(\mathbf{y},z_1 \mid \mathbf{x,z_0})dz_1d\mathbf{z_0} \nonumber \\
&= \int_{\mathbf{z_0},z_1} f(\mathbf{y},z_1 \mid \mathbf{x,z_0})f(\mathbf{z_0})d\mathbf{z_0}dz_1 \nonumber \\
&= \int_{\mathbf{z_0},z_1} f(\mathbf{y} \mid \mathbf{x,z_0},z_1)f(z_1 \given \mathbf{x,z_0})f(\mathbf{z_0})d\mathbf{z_0}dz_1 \nonumber \\
&= \int_{\mathbf{z_0},z_1} f(\mathbf{y} \mid \mathbf{x,z_0},z_1)f(z_1 \given \mathbf{z_0})f(\mathbf{z_0})d\mathbf{z_0}dz_1 \nonumber \\
&= \int_{\mathbf{z_0},z_1} f(\mathbf{y} \mid \mathbf{x,z_0},z_1) f(z_1, \mathbf{z_0})d\mathbf{z_0}dz_1. \nonumber 
\end{align}

We complete the proof by showing that (a) or (b) must hold.
Suppose for a contradiction that both (a) and (b) are violated. Then
there is a path from $\bX$ to $Z_1$ that is $d$-connecting given $\mathbf{Z_0}$  and a path from $\bY$ to $Z_1$
that is $d$-connecting given $\bX \cup \mathbf{Z_0}$. Let $p$ be a proper path from $X \in \mathbf{X}$ to $Z_1$ that is d-connecting given $\mathbf{Z_0}$ in $\g[D]$ and let $q$ be a path from $Z_1$ to $Y \in \mathbf{Y}$ that is $d$-connecting given $\bX \cup \mathbf{Z_0}$ in $\g[D]$.
We will show that this contradicts that $\mathbf{Z_0}$ satisfies the adjustment criterion relative to $(\mathbf{X,Y})$ in $\g[D]$.

We first show, by contradiction, that $q$ also is d-connecting given $\mathbf{Z_0}$.
Thus, assume that $q=\langle Z_1,\ldots, Y \rangle$ is blocked by $\mathbf{Z_0}$.
Since $q$ is d-connecting given $\mathbf{X} \cup \mathbf{Z_0}$ and blocked by $\mathbf{Z_0}$, it must contain at least one collider in $\An(\bX,\g[D]) \setminus \An(\mathbf{Z_0},\g[D])$.
Let $C$ be the collider closest to $Y$ on $q$ such that $C \in \An(\bX,\g[D]) \setminus \An(\mathbf{Z_0},\g[D])$. Let $r=\langle C,\ldots,X' \rangle, X' \in \mathbf{X}$ be a shortest directed path from $C$ to $\mathbf{X}$ (possibly of zero length). Then no node on $q(C,Y)$ or $r$, except possibly $C$, is in $\mathbf{X}$.
 We now concatenate the paths $-r(X',C)$ and $q(C,Y)$, while taking out possible loops. Hence, let $V$ be the node closest to $X'$ on $r$ that is also on $q(C,Y)$.
Then $-r(X',V) \oplus q(V,Y)$ is non-causal since either $-r(X',V)$ is of non-zero length, or $X'=V=C$ and $q(C,Y)$ is a path into $C$, because $C$ is a collider on $q$. By the choice of $C$ and $r$, $-r(X',V) \oplus q(V,Y)$ is a proper path that is d-connecting given $\mathbf{Z_0}$. This contradicts that $\mathbf{Z_0}$ satisfies the adjustment criterion relative to $(\mathbf{X,Y})$ in~$\g[D]$.
Thus, $q$ is also $d$-connecting given $\mathbf{Z_0}$.

Let $\tilde{p}$ and $\tilde{q}$ be paths in the proper back-door graph $\gpbd[D]{XY}$ constituted by the same sequences of nodes as $p$ and $q$ in~$\g[D]$ respectively.
We first prove that the paths $\tilde{p}$ and $\tilde{q}$ exist in~$\gpbd[D]{XY}$.
Path $q$ is d-connecting given $\bX \cup \mathbf{Z_0}$, so any node in $\mathbf{X}$ on $q$ must be a collider on $q$.
Since $\gpbd[D]{XY}$ is obtained from $\g[D]$ by removing certain edges out of $\mathbf{X}$, no edges from $q$ are removed and $\tilde{q}$ exists in~$\gpbd[D]{XY}$.

Since $p$ is proper, for $\tilde{p}$ to exist in~$\gpbd[D]{XY}$, it is enough to show that $p$ does not start with an edge of type $X \rightarrow W$ in~$\g[D]$ where $W$ lies on a proper causal path from $X$ to $\mathbf{Y}$ in~$\g[D]$.
Suppose for a contradiction that $p$ does start with $X \rightarrow W$. Then $W \in \fb{\g[D]}$. Then either $p$ is a directed path from $X$ to $Z_1$ so that $Z_1 \in \De(W,\g[D])$, or $p$ is non-causal and there is a collider $C'$ on $p$ such that $C' \in \De(W,\g[D])$. In the former case, since $\De(W,\g[D]) \subseteq \fb{\g[D]}$, it follows that $Z_1 \in \fb{\g[D]}$, which contradicts the choice of $Z_1$ because $Z_1 \notin \De(\mathbf{X},\g[D])$. In the latter case, since $p$ is d-connecting given $\mathbf{Z_0}$, we have $\De(C',\g[D]) \cap \mathbf{Z_0} \neq \emptyset$. Combining this with $\De(C',\g[D]) \subseteq \fb{\g[D]}$, it follows that $\fb{\g[D]} \cap \mathbf{Z_0} \neq \emptyset$ which contradicts that $\mathbf{Z_0}$ satisfies 
\forb{} relative to $(\mathbf{X,Y})$ in $\g[D]$.
Thus, $\tilde{p}$ exists in~$\gpbd[D]{\mathbf{XY}}$.

We now show that $\tilde{p}$ and $\tilde{q}$ are d-connecting given $\mathbf{Z_0}$ in~$\gpbd[D]{XY}$. The collider/non-collider status of any node on $\tilde{p}$ and $\tilde{q}$ in~$\gpbd[D]{XY}$ is the same as the collider/non-collider status of that same node on $p$ and $q$ in~$\g[D]$ respectively. So $\tilde{p}$ and $\tilde{q}$ are both d-connecting given $\mathbf{Z_0}$, unless every causal path from a collider on either $p$ or $q$ to $\mathbf{Z_0}$ contains a first edge on a proper causal path from $\mathbf{X}$ to $\mathbf{Y}$ in~$\g[D]$. Any such causal path also contains a node in $\fb{\g[D]}$, so a node in $\mathbf{Z_0}$ would be a descendant of $\fb{\g[D]}$ in~$\g[D]$. Since $\fb{\g[D]}$ is descendral, it follows that $\mathbf{Z_0} \cap \fb{\g[D]} \neq \emptyset$. This contradicts that $\mathbf{Z_0}$ satisfies the adjustment criterion relative to $(\mathbf{X,Y})$ in~$\g[D]$, 
specifically \forb{}.

Since $\tilde{p}$ is d-connecting given $\mathbf{Z_0}$ in~$\gpbd[D]{XY}$, by Theorem~\ref{theorem:moralization} it follows that
there is a path $a$ from $X$ to $Z_1$ that does not contain a node in $\mathbf{Z_0}$ in the moral induced subgraph of $\gpbd[D]{XY}$ on nodes $\An( \bX  \cup \mathbf{Y} \cup \mathbf{Z_0} \cup \{Z_1\},\gpbd[D]{XY})$.
 Similarly, $\tilde{q}$ is a d-connecting path from $Z_1$ to $Y$ given $\mathbf{Z_0}$ in~$\gpbd[D]{XY}$, so there is a path $b$ from $Z_1$ to $Y$ that does not contain a node in $\mathbf{Z_0}$ in the moral induced subgraph of $\gpbd[D]{XY}$ on nodes $\An( \bX  \cup \mathbf{Y} \cup \mathbf{Z_0} \cup \{Z_1\},\gpbd[D]{XY})$.
By combining paths $a$ and $b$ we get a path $c$ from $\bX$ to $\bY$ that does not contain a node in $\mathbf{Z_0}$ in the moral induced subgraph of $\gpbd[D]{XY}$ on nodes $\An( \bX  \cup \mathbf{Y} \cup \mathbf{Z_0} \cup \{Z_1\},\gpbd[D]{XY})$.
Since $Z_1 \in \An(\mathbf{X} \cup \mathbf{Y},\g[D])$ and $\gpbd[D]{XY}$ is obtained from $\g[D]$ by removing certain edges out of $\mathbf{X}$, it follows that $Z_1 \in \An(\mathbf{X} \cup \mathbf{Y},\gpbd[D]{XY})$. Then $\An(\mathbf{X} \cup \bY \cup \mathbf{Z_0} \cup \{Z_1\},\gpbd[D]{XY}) = \An(\mathbf{X} \cup \bY \cup \mathbf{Z_0},\gpbd[D]{XY})$. Hence, $c$ is a path from $\mathbf{X}$ to $\mathbf{Y}$ that does not contain a node in $\mathbf{Z_0}$ in the moral induced subgraph of $\gpbd[D]{XY}$ on nodes $\An( \bX  \cup \mathbf{Y} \cup \mathbf{Z_0},\gpbd[D]{XY})$. 
Thus, by Theorem~\ref{theorem:moralization}, $\mathbf{X}$ and $\mathbf{Y}$ are d-connected given $\mathbf{Z_0}$ in~$\gpbd[D]{XY}$.
 By Theorem~\ref{theorem:gac-alt}, this contradicts that $\mathbf{Z_0}$ satisfies the adjustment criterion relative to $(\mathbf{X,Y})$ in~$\g[D]$.
\end{proofof}

\begin{lemma} 
Let $\mathbf{X},\mathbf{Y}$ and $\mathbf{Z_0}$ be pairwise disjoint node sets in a $\DAG$ $\g[D]$ such that $\mathbf{Z_0}$ satisfies the adjustment criterion relative to $(\mathbf{X,Y})$ in $\g[D]$. Let $\mathbf{Z} = \mathbf{Z_0} \cup \An(\mathbf{X} \cup \mathbf{Y},\g[D]) \setminus ( \De(\mathbf{X}, \g[D]) \cup \mathbf{Y})$.
Additionally, let $\mathbf{Z_D} = \mathbf{Z} \cap \De(\mathbf{X}, \g[D])$, $\mathbf{Z_N} = \mathbf{Z} \setminus \De(\mathbf{X}, \g[D])$, $\mathbf{Y_D} = \mathbf{Y} \cap \De(\mathbf{X},\g[D])$ and $\mathbf{Y_N} = \mathbf{Y} \setminus \De(\mathbf{X}, \g[D])$.
Then the following statements hold:

\begin{enumerate}[label=(\roman*)]
\item\label{l:g1}  $(\mathbf{X} \cup \mathbf{Y_N} \cup \mathbf{Z}) \cap \fb{\g[D]} = \emptyset$, and

\item\label{l:g2} if $p = \langle A, \dots, Y_d \rangle$ is a non-causal path from a node $A \in \mathbf{X} \cup \mathbf{Y_N} \cup \mathbf{Z}$ to a node $Y_d \in \mathbf{Y_D}$, then $p$ is blocked by $(\mathbf{X} \cup \mathbf{Y_N} \cup \mathbf{Z_N}) \setminus \{A\}$ in~$\g[D]$, and

\item\label{l:s1} $\mathbf{Y_D} \dsepp \mathbf{Z_D} \given \mathbf{Y_N} \cup \mathbf{X} \cup \mathbf{Z_N}$ in~$\g[D]$, where $\mathbf{Y_N} = \emptyset$ is allowed, and

\item\label{l:s2} if $\mathbf{Y_N} = \emptyset$ then $\mathbf{Z_N}$ is a generalized back-door set relative to $(\mathbf{X,Y})$ in~$\g[D]$, and

\item\label{l:s3}  the empty set is a generalized back-door set relative to $((\mathbf{X} \cup \mathbf{Y_N} \cup \mathbf{Z_N}),\mathbf{Y_D})$ in~$\g[D]$, and

\item\label{l:s4} $\mathbf{Y_N} \cup \mathbf{Z_N} \dsepp \mathbf{Y_D} \given \mathbf{X}$ in~$\g[D]_{\overline{\mathbf{X}}\underline{\mathbf{Y_N} \cup  \mathbf{Z_N}}}$, and

\item\label{l:s5} $\mathbf{X} \dsepp \mathbf{Z_N} \given \mathbf{Y_N}$ in~$\g[D]_{\overline{\mathbf{X}}}$.
\end{enumerate}

\label{lemma:shpitser ynyd}
\end{lemma}

\begin{proofof}[Lemma~\ref{lemma:shpitser ynyd}]

\ref{l:g1} By Lemma~\ref{lemma:zprime satisfies ac}, $\mathbf{Z}$ satisfies the adjustment criterion relative to $(\mathbf{X,Y})$ in $\g[D]$, implying that $\mathbf{Z} \cap \fb{\g[D]} = \emptyset$. By definition $\mathbf{Y_N} \cap \De(\mathbf{X},\g[D]) = \emptyset$ and $\De(\mathbf{X},\g[D]) \supseteq \fb{\g}$. Hence, $\mathbf{Y_N} \cap \fb{\g[D]} = \emptyset$. It is only left to show that $\mathbf{X} \cap \fb{\g[D]} = \emptyset$.

Suppose for a contradiction that $\mathbf{X} \cap \fb{\g[D]} \neq \emptyset$. Let $V \notin \mathbf{X}$ be a node on a proper causal path $p$ from $\mathbf{X}$ to $\mathbf{Y}$ (possibly $V \in \mathbf{Y}$) such that $V \in \An(\mathbf{X},\g[D])$. Let $q = \langle V,\dots,X \rangle$, $X \in \mathbf{X}$, be a shortest causal path from $V$ to $\mathbf{X}$. All nodes on $p(V,Y)$ and $q$ are in $\fb{\g[D]}$. We now concatenate $-q$ and $p(V,Y)$, while taking out loops. Hence, let $W$ be the node closest to $X$ on $q$ that is also on $p(V,Y)$. Then $r = -q(X,W) \oplus p(W,Y)$ is a proper non-causal path from $\mathbf{X}$ to $\mathbf{Y}$ that cannot be blocked by $\mathbf{Z}$ since $\mathbf{Z} \cap \fb{\g[D]} = \emptyset$, which contradicts Lemma~\ref{lemma:zprime satisfies ac}.

\ref{l:g2} We distinguish the cases that $p$ is~\ref{l:sga} out of, or~\ref{l:sgb} into $Y_d$.
\begin{enumerate}[label = (\alph*)]
\item\label{l:sga} Let $p$ be out of $Y_d$. Since $Y_d \in \fb{\g[D]}$ and since by~\ref{l:g1} node $A \notin \fb{\g[D]}$, there is at least one collider on $p$. The collider closest to $Y_d$ on $p$ and all of its descendants are also in $\fb{\g[D]}$. It then follows from~\ref{l:g1} that $p$ is blocked by $(\mathbf{X} \cup \mathbf{Y_N} \cup \mathbf{Z_N}) \setminus \{A\}$.
\item\label{l:sgb} Let $p$ be into $Y_d$. Since $p$ is a non-causal path from $A$ to $Y_d$, there is at least one node on $p$ that has two edges out of it. Let $B$ be the closest such node to $Y_d$ on $p$. Then $B \in \An(Y_d,\g[D])$. If any node on $p(B,Y_d)$ is in $\mathbf{X} \cup \mathbf{Y_N} \cup \mathbf{Z_N}$, then~\ref{l:g2} holds.

Hence, assume no node on $p(B,Y_d)$ is in $\mathbf{X} \cup \mathbf{Y_N} \cup \mathbf{Z_N}$. Then $B \notin \mathbf{X} \cup \mathbf{Y_N} \cup \mathbf{Z_N}$. Since $B \notin \mathbf{Z_N}$, $\mathbf{Z_N} \supseteq \An(\mathbf{X} \cup \mathbf{Y},\g[D]) \setminus (\De(\mathbf{X},\g[D]) \cup \mathbf{Y})$, it follows that $B \notin \An(\mathbf{X} \cup \mathbf{Y},\g[D]) \setminus (\De(\mathbf{X},\g[D]) \cup \mathbf{Y})$. Additionally, since $B \in \An(Y_d,\g[D])$, $B \in \An(\mathbf{X} \cup \mathbf{Y},\g[D])$. Combining $B \notin \An(\mathbf{X} \cup \mathbf{Y},\g[D]) \setminus (\De(\mathbf{X},\g[D]) \cup \mathbf{Y})$ and $B \in \An(\mathbf{X} \cup \mathbf{Y},\g[D])$ implies $B \in \De(\mathbf{X},\g[D]) \cup \mathbf{Y}$. Furthermore, $B \notin \mathbf{X} \cup \mathbf{Y_N}$ and $\mathbf{Y_D} \subseteq \De(\mathbf{X},\g[D]) \setminus \mathbf{X}$, so $B \in \De(\mathbf{X},\g[D]) \setminus \mathbf{X}$. But $B \in \An(Y_d,\g[D])$ through $p(B,Y_d)$ on which no other node is in $\mathbf{X}$, so $B \in \fb{\g[D]}$. Now, $-p(B,A)$ is a path out of $B$, where $B \in \fb{\g[D]}$ and $A \notin \fb{\g[D]}$ by~\ref{l:g1}. Using the same reasoning as in \ref{l:sga}, $p(A,B)$ is blocked by $(\mathbf{X} \cup \mathbf{Y_N} \cup \mathbf{Z_N}) \setminus \{A\}$. Hence, $p$ is also blocked by $(\mathbf{X} \cup \mathbf{Y_N} \cup \mathbf{Z_N}) \setminus \{A\}$.
\end{enumerate}

\ref{l:s1} By~\ref{l:g2} every non-causal path from $\mathbf{Z_D}$ to $\mathbf{Y_D}$ is blocked by $\mathbf{X} \cup \mathbf{Y_N} \cup \mathbf{Z_N}$. 
Additionally, $\mathbf{Z_D} \cap \fb{\g[D]} = \emptyset$ and  $\mathbf{Z_D} \subseteq \De(\mathbf{X},\g[D])$, so no node in $\mathbf{Z_D}$ can have a proper causal path to $\mathbf{Y_D}$ in $\g[D]$. Hence, any causal path $p$ from $\mathbf{Z_D}$ to $\mathbf{Y_D}$ in~$\g[D]$, must be non-proper with respect to $\mathbf{X}$, that is, $p$ has to contain a node in $\mathbf{X}$ as a non-collider. Hence, every causal path from $\mathbf{Z_D}$ to $\mathbf{Y_D}$ is also blocked by $\mathbf{X} \cup \mathbf{Y_N} \cup \mathbf{Z_N}$ in~$\g[D]$.

\ref{l:s2} Follows directly from $\mathbf{Z_N} \cap \De(\mathbf{X},\g[D]) = \emptyset$ and~\ref{l:g2} for $\mathbf{Y_N} = \emptyset$.

\ref{l:s3} Follows directly from~\ref{l:g2}.

\ref{l:s4} Since~$\g[D]_{\overline{\mathbf{X}}\underline{\mathbf{Y_N} \cup  \mathbf{Z_N}}}$ does not contain edges into $\mathbf{X}$, all paths from $\mathbf{Y_N} \cup  \mathbf{Z_N}$ to $\mathbf{Y_D}$ in $\g[D]_{\overline{\mathbf{X}}\underline{\mathbf{Y_N} \cup  \mathbf{Z_N}}}$ that contain a collider are blocked by $\mathbf{X}$. Hence, it is enough to prove that $\mathbf{X}$ blocks any path from $\mathbf{Y_N} \cup  \mathbf{Z_N}$ to $\mathbf{Y_D}$ in $\g[D]_{\overline{\mathbf{X}}\underline{\mathbf{Y_N} \cup  \mathbf{Z_N}}}$ that does not contain a collider. Let $r = \langle A, \dots , Y_d \rangle$, $A \in \mathbf{Y_N} \cup  \mathbf{Z_N}$, $Y_d \in \mathbf{Y_D}$ be an arbitrary path from $\mathbf{Y_N} \cup  \mathbf{Z_N}$ to $\mathbf{Y_D}$ in $\g[D]_{\overline{\mathbf{X}}\underline{\mathbf{Y_N} \cup  \mathbf{Z_N}}}$ that does not contain a collider. Since there are no edges out of $\mathbf{Y_N} \cup  \mathbf{Z_N}$ in $\g[D]_{\overline{\mathbf{X}}\underline{\mathbf{Y_N} \cup  \mathbf{Z_N}}}$, it follows that $r$ is into $A$ and that $r$ does not contain non-colliders that are in $\mathbf{Y_N} \cup  \mathbf{Z_N}$.
Let $r'$ in~$\g[D]$ be the path constituted by the same sequence of nodes as $r$ in $\g[D]_{\overline{\mathbf{X}}\underline{\mathbf{Y_N} \cup  \mathbf{Z_N}}}$. Since removing edges cannot d-connect blocked paths, it is enough to prove that $r'$ is blocked by $\mathbf{X}$ in $\g[D]$.
Since $r'$ is into $A$, $r'$ is a non-causal path from $\mathbf{Y_N} \cup  \mathbf{Z_N}$ to $\mathbf{Y_D}$, so by~\ref{l:g2}, $r'$ is blocked by $(\mathbf{X} \cup \mathbf{Y_N} \cup \mathbf{Z_N}) \setminus \{A\}$ in $\g[D]$. Since $A \in \mathbf{Y_N} \cup  \mathbf{Z_N}$ and $(\mathbf{Y_N} \cup  \mathbf{Z_N}) \cap \mathbf{X}=\emptyset$, it follows that $A \notin \mathbf{X}$. Then $r'$ is blocked by $\mathbf{X} \cup ((\mathbf{Y_N} \cup \mathbf{Z_N}) \setminus\{A\})$. Furthermore, since $r$ does not contain non-colliders in $\mathbf{Y_N} \cup  \mathbf{Z_N}$, the same is true for $r'$, so $r'$ is blocked by $\mathbf{X}$.
 
\ref{l:s5} Since $\g[D]_{\overline{\mathbf{X}}}$ does not contain edges into $\mathbf{X}$, all paths from $\mathbf{X}$ to $\mathbf{Z_N}$ in~$\g[D]_{\overline{\mathbf{X}}}$ are out of $\mathbf{X}$.
Since $\mathbf{Z_N} \cap \De(\mathbf{X},\g[D]) = \emptyset$ and $\De(\mathbf{X},\g[D]) \supseteq \De(\mathbf{X},\g[D]_{\overline{\mathbf{X}}})$, all paths from $\mathbf{X}$ to $\mathbf{Z_N}$ in~$\g[D]_{\overline{\mathbf{X}}}$ contain at least one collider. The closest collider to $\mathbf{X}$ on any such path is then in $\De(\mathbf{X},\g[D]_{\overline{\mathbf{X}}})$. Since $\mathbf{Y_N} \cap \De(\mathbf{X},\g[D]) = \emptyset$ and $\De(\mathbf{X},\g[D]) \supseteq \De(\mathbf{X},\g[D]_{\overline{\mathbf{X}}})$, this path is blocked by $\mathbf{Y_N}$ in $\g[D]_{\overline{\mathbf{X}}}$.
\end{proofof}

\vskip 0.2in
\bibliography{biblioteka}

\begin{thebibliography}{}

\bibitem[Ali et~al., 2009]{ali2009markov}
Ali, A.~R., Richardson, T.~S., and Spirtes, P. (2009).
\newblock Markov equivalence for ancestral graphs.
\newblock {\em Ann. Stat.}, 37:2808--2837.

\bibitem[Bareinboim et~al., 2014]{Bareinboim2014}
Bareinboim, E., Tian, J., and Pearl, J. (2014).
\newblock Recovering from selection bias in causal and statistical inference.
\newblock In {\em Proceedings of AAAI 2014}, pages 2410--2416.

\bibitem[Chickering, 2002]{Chickering02-optimal}
Chickering, D.~M. (2002).
\newblock Optimal structure identification with greedy search.
\newblock {\em J. Mach. Learn. Res.}, 3:507--554.

\bibitem[Claassen et~al., 2013]{ClaassenMooijHeskes13}
Claassen, T., Mooij, J., and Heskes, T. (2013).
\newblock Learning sparse causal models is not {NP}-hard.
\newblock In {\em Proceedings of UAI 2013}, pages 172--181.

\bibitem[Colombo and Maathuis, 2014]{colombo2014order}
Colombo, D. and Maathuis, M.~H. (2014).
\newblock Order-independent constraint-based causal structure learning.
\newblock {\em J. Mach. Learn. Res.}, 15:3741--3782.

\bibitem[Colombo et~al., 2012]{Colombo2012}
Colombo, D., Maathuis, M.~H., Kalisch, M., and Richardson, T.~S. (2012).
\newblock Learning high-dimensional directed acyclic graphs with latent and
  selection variables.
\newblock {\em Ann. Stat.}, 40:294--321.

\bibitem[Correa and Bareinboim, 2017]{correa2017causal}
Correa, J.~D. and Bareinboim, E. (2017).
\newblock Causal effect identification by adjustment under confounding and
  selection biases.
\newblock In {\em Proceedings of AAAI 2017}, pages 3740--3746.

\bibitem[De~Luna et~al., 2011]{DeLunaEtAl11}
De~Luna, X., Waernbaum, I., and Richardson, T.~S. (2011).
\newblock Covariate selection for the nonparametric estimation of an average
  treatment effect.
\newblock {\em Biometrika}, 98(4):861--875.

\bibitem[Entner et~al., 2013]{EntnerHoyerSpirtes13}
Entner, D., Hoyer, P.~O., and Spirtes, P. (2013).
\newblock Data-driven covariate selection for nonparametric estimation of
  causal effects.
\newblock In {\em Proceedings of AISTATS 2013}, pages 256--264.

\bibitem[Frot et~al., 2018]{frot2017learning}
Frot, B., Nandy, P., and Maathuis, M.~H. (2018).
\newblock Learning directed acyclic graphs with hidden variables via latent
  {G}aussian graphical model selection.
\newblock arXiv preprint arXiv:1708.01151.

\bibitem[Greenland et~al., 1999]{greenland1999causal}
Greenland, S., Pearl, J., and Robins, J.~M. (1999).
\newblock Causal diagrams for epidemiologic research.
\newblock {\em Epidemiology}, 10(1):37--48.

\bibitem[Guo and Dawid, 2010]{guo2010sufficient}
Guo, H. and Dawid, P.~A. (2010).
\newblock Sufficient covariates and linear propensity analysis.
\newblock In {\em Proceedings of AISTATS 2010}, pages 281--288.

\bibitem[Hahn, 1998]{hahn1998role}
Hahn, J. (1998).
\newblock On the role of the propensity score in efficient semiparametric
  estimation of average treatment effects.
\newblock {\em Econometrica}, pages 315--331.

\bibitem[Hahn, 2004]{hahn2004functional}
Hahn, J. (2004).
\newblock Functional restriction and efficiency in causal inference.
\newblock {\em Rev. Econ. Stat.}, 86(1):73--76.

\bibitem[Heinze-Deml et~al., 2017]{heinze2017causal}
Heinze-Deml, C., Maathuis, M.~H., and Meinshausen, N. (2017).
\newblock Causal structure learning.
\newblock {\em Annu. Rev. Stat. Appl.}, 5:371--391.

\bibitem[Kalisch et~al., 2012]{kalischpcalg}
Kalisch, M., M\"achler, M., Colombo, D., Maathuis, M.~H., and B\"uhlmann, P.
  (2012).
\newblock Causal inference using graphical models with the {R} package {pcalg}.
\newblock {\em J. Stat. Softw.}, 47(11):1--26.

\bibitem[Koster, 2002]{koster2002marginalizing}
Koster, J.~T. (2002).
\newblock Marginalizing and conditioning in graphical models.
\newblock {\em Bernoulli}, 8(6):817--840.

\bibitem[Kuroki and Cai, 2004]{kuroki04selection}
Kuroki, M. and Cai, Z. (2004).
\newblock Selection of identifiability criteria for total effects by using path
  diagrams.
\newblock In {\em Proceedings of UAI 2004}, pages 333--340.

\bibitem[Kuroki and Miyakawa, 2003]{kuroki2003covariate}
Kuroki, M. and Miyakawa, M. (2003).
\newblock Covariate selection for estimating the causal effect of control plans
  by using causal diagrams.
\newblock {\em J. Roy. Stat. Soc. B}, 65(1):209--222.

\bibitem[Lauritzen et~al., 1990]{lauritzen1990independence}
Lauritzen, S.~L., Dawid, P.~A., Larsen, B.~N., and Leimer, H.-G. (1990).
\newblock Independence properties of directed {M}arkov fields.
\newblock {\em Networks}, 20(5):491--505.

\bibitem[Lauritzen and Spiegelhalter, 1988]{lauritzen1988local}
Lauritzen, S.~L. and Spiegelhalter, D.~J. (1988).
\newblock Local computations with probabilities on graphical structures and
  their application to expert systems.
\newblock {\em J. Roy. Stat. Soc. B}, 50(2):157--224.

\bibitem[Maathuis and Colombo, 2015]{maathuis2013generalized}
Maathuis, M.~H. and Colombo, D. (2015).
\newblock A generalized back-door criterion.
\newblock {\em Ann. Stat.}, 43:1060--1088.

\bibitem[Maathuis et~al., 2010]{MaathuisColomboKalischBuehlmann10}
Maathuis, M.~H., Colombo, D., Kalisch, M., and B\"uhlmann, P. (2010).
\newblock Predicting causal effects in large-scale systems from observational
  data.
\newblock {\em Nat. Methods}, 7:247--248.

\bibitem[Maathuis et~al., 2009]{MaathuisKalischBuehlmann09}
Maathuis, M.~H., Kalisch, M., and B\"uhlmann, P. (2009).
\newblock Estimating high-dimensional intervention effects from observational
  data.
\newblock {\em Ann. Stat.}, 37:3133--3164.

\bibitem[Malinsky and Spirtes, 2017]{malinsky2017estimating}
Malinsky, D. and Spirtes, P. (2017).
\newblock Estimating bounds on causal effects in high-dimensional and possibly
  confounded systems.
\newblock {\em Int. J. of Approx. Reason.}, 88:371--384.

\bibitem[Mardia et~al., 1980]{mardia1980multivariate}
Mardia, K.~V., Kent, J.~T., and Bibby, J.~M. (1980).
\newblock {\em Multivariate Analysis (Probability and Mathematical
  Statistics)}.
\newblock Academic Press London.

\bibitem[Meek, 1995]{meek1995causal}
Meek, C. (1995).
\newblock Causal inference and causal explanation with background knowledge.
\newblock In {\em Proceedings of UAI 1995}, pages 403--410.

\bibitem[Nandy et~al., 2018]{nandyMaathuisArges}
Nandy, P., Hauser, A., and Maathuis, M.~H. (2018).
\newblock High-dimensional consistency in score-based and hybrid structure
  learning.
\newblock {\em Ann. Stat.}
\newblock To appear.

\bibitem[Nandy et~al., 2017]{nandy2014estimating}
Nandy, P., Maathuis, M.~H., and Richardson, T.~S. (2017).
\newblock Estimating the effect of joint interventions from observational data
  in sparse high-dimensional settings.
\newblock {\em Ann. Stat.}, 45(2):647--674.

\bibitem[Pearl, 1993]{pearl1993bayesian}
Pearl, J. (1993).
\newblock Comment: {G}raphical models, causality and intervention.
\newblock {\em Stat. Sci.}, 8:266--269.

\bibitem[Pearl, 2009]{Pearl2009}
Pearl, J. (2009).
\newblock {\em Causality: {M}odels, Reasoning, and Inference}.
\newblock Cambridge University Press, New York, NY, second edition.

\bibitem[Pearl and Paz, 2014]{pearl2014confounding}
Pearl, J. and Paz, A. (2014).
\newblock Confounding equivalence in causal inference.
\newblock {\em J. Causal Infer.}, 2(1):75--93.

\bibitem[Perkovi\'c et~al., 2017]{perkovic17}
Perkovi\'c, E., Kalisch, M., and Maathuis, M.~H. (2017).
\newblock Interpreting and using {CPDAG}s with background knowledge.
\newblock In {\em Proceedings of UAI 2017}.

\bibitem[Perkovi\'c et~al., 2015]{perkovic15_uai}
Perkovi\'c, E., Textor, J., Kalisch, M., and Maathuis, M.~H. (2015).
\newblock A complete generalized adjustment criterion.
\newblock In {\em Proceedings of UAI 2015}, pages 682--691.

\bibitem[Richardson, 2003]{richardson2003markov}
Richardson, T.~S. (2003).
\newblock Markov properties for acyclic directed mixed graphs.
\newblock {\em Scand. J. Statist.}, 30:145--157.

\bibitem[Richardson and Spirtes, 2002]{richardson2002ancestral}
Richardson, T.~S. and Spirtes, P. (2002).
\newblock Ancestral graph {M}arkov models.
\newblock {\em Ann. Stat.}, 30:962--1030.

\bibitem[Robins, 1986]{robins1986new}
Robins, J.~M. (1986).
\newblock A new approach to causal inference in mortality studies with a
  sustained exposure period-application to control of the healthy worker
  survivor effect.
\newblock {\em Math. Mod.}, 7:1393--1512.

\bibitem[Rubin, 2008]{Rubin2008}
Rubin, D. (2008).
\newblock Author's reply.
\newblock {\em Stat. Med.}, 27:2741--2742.

\bibitem[Shachter, 1998]{Shachter1998}
Shachter, R.~D. (1998).
\newblock Bayes-ball: {T}he rational pastime.
\newblock In {\em Proceedings of UAI 1998}, pages 480--487.

\bibitem[Shpitser, 2012]{shpitser2012avalidity}
Shpitser, I. (2012).
\newblock Appendum to ``{O}n the validity of covariate adjustment for
  estimating causal effects''.
\newblock Personal communication.

\bibitem[Shpitser and Pearl, 2006]{shpitser2006identification}
Shpitser, I. and Pearl, J. (2006).
\newblock Identification of joint interventional distributions in recursive
  semi-{M}arkovian causal models.
\newblock In {\em Proceedings of AAAI 2006}, pages 1219--1226.

\bibitem[Shpitser et~al., 2010]{shpitser2012validity}
Shpitser, I., VanderWeele, T., and Robins, J.~M. (2010).
\newblock On the validity of covariate adjustment for estimating causal
  effects.
\newblock In {\em Proceedings of UAI 2010}, pages 527--536.

\bibitem[Shrier, 2008]{Shrier2008}
Shrier, I. (2008).
\newblock Letter to the editor.
\newblock {\em Stat. Med.}, 27:2740--2741.

\bibitem[Spirtes et~al., 2000]{spirtes2000causation}
Spirtes, P., Glymour, C., and Scheines, R. (2000).
\newblock {\em Causation, Prediction, and Search}.
\newblock MIT Press, Cambridge, MA, second edition.

\bibitem[Textor et~al., 2016]{textor2016robust}
Textor, J., van~der Zander, B., Gilthorpe, M.~S., Li{\'s}kiewicz, M., and
  Ellison, G.~T. (2016).
\newblock Robust causal inference using directed acyclic graphs: the {R}
  package 'dagitty'.
\newblock {\em Int. J. Epidemiol.}, 45(6):1887--1894.

\bibitem[Tian and Pearl, 2002]{tian2002general}
Tian, J. and Pearl, J. (2002).
\newblock A general identification condition for causal effects.
\newblock In {\em Proceedings of AAAI 2002}, pages 567--573.

\bibitem[van~der Zander and Li\'skiewicz, 2016]{vanDerZander16}
van~der Zander, B. and Li\'skiewicz, M. (2016).
\newblock Separators and adjustment sets in {M}arkov equivalent {DAG}s.
\newblock In {\em Proceedings of AAAI 2016}, pages 3315--3321.

\bibitem[van~der Zander et~al., 2014]{vanconstructing}
van~der Zander, B., Li\'skiewicz, M., and Textor, J. (2014).
\newblock Constructing separators and adjustment sets in ancestral graphs.
\newblock In {\em Proceedings of UAI 2014}, pages 907--916.

\bibitem[van~der Zander et~al., 2018]{van2018separators}
van~der Zander, B., Li{\'s}kiewicz, M., and Textor, J. (2018).
\newblock Separators and adjustment sets in causal graphs: Complete criteria
  and an algorithmic framework.
\newblock arXiv preprint arXiv:1803.00116.

\bibitem[VanderWeele and Shpitser, 2011]{vanderweele2011new}
VanderWeele, T.~J. and Shpitser, I. (2011).
\newblock A new criterion for confounder selection.
\newblock {\em Biometrics}, 67(4):1406--1413.

\bibitem[West and Koch, 2014]{West2014}
West, S.~G. and Koch, T. (2014).
\newblock Restoring causal analysis to structural equation modeling.
\newblock {\em Struct. Equ. Modeling}, 21:161--166.

\bibitem[Westreich and Greenland, 2013]{Westreich2013}
Westreich, D. and Greenland, S. (2013).
\newblock {{T}he table 2 fallacy: {P}resenting and interpreting confounder and
  modifier coefficients}.
\newblock {\em Am. J. Epidemiol.}, 177:292--298.

\bibitem[Wright, 1921]{wright1921correlation}
Wright, S. (1921).
\newblock Correlation and causation.
\newblock {\em J. Agric. Res.}, 20(7):557--585.

\bibitem[Zhang, 2006]{zhang2006causal}
Zhang, J. (2006).
\newblock {\em Causal Inference and Reasoning in Causally Insufficient
  Systems}.
\newblock PhD thesis, Carnegie Mellon University.

\bibitem[Zhang, 2008a]{zhang2008causal}
Zhang, J. (2008a).
\newblock Causal reasoning with ancestral graphs.
\newblock {\em J. Mach. Learn. Res.}, 9:1437--1474.

\bibitem[Zhang, 2008b]{zhang2008completeness}
Zhang, J. (2008b).
\newblock On the completeness of orientation rules for causal discovery in the
  presence of latent confounders and selection bias.
\newblock {\em Artif. Intell.}, 172:1873--1896.

\end{thebibliography}
\end{document}